%% file: main.tex
\documentclass{template}
\input{math_commands.tex}

\usetikzlibrary{shapes.geometric}
\usetikzlibrary{positioning,fit}
\usepackage{enumitem}
\usepackage{wrapfig}
\usepackage[normalem]{ulem}
\usepackage{xcolor}
\online
\begin{document}

\ensubject{fdsfd}

\BeginPage{1} %

\title[Learning to Optimize: A Tutorial]
{Learning to optimize: A tutorial for continuous and mixed-integer optimization}

\author[]{Xiaohan Chen}{{xiaohan.chen@alibaba-inc.com}}
\author[]{Jialin Liu}{{jialin.liu@alibaba-inc.com}}
\author[$\ast$]{Wotao Yin}{wotao.yin@alibaba-inc.com}

\AuthorMark{Chen}


\address[]{Decision Intelligence Lab, Alibaba DAMO Academy, Bellevue, WA {\rm 98004}, USA}

\abstract{Learning to Optimize (L2O) stands at the intersection of traditional optimization and machine learning, utilizing the capabilities of machine learning to enhance conventional optimization techniques. As real-world optimization problems frequently share common structures, L2O provides a tool to exploit these structures for better or faster solutions. This tutorial dives deep into L2O techniques, introducing how to accelerate optimization algorithms, promptly estimate the solutions, or even reshape the optimization problem itself, making it more adaptive to real-world applications. By considering the prerequisites for successful applications of L2O and the structure of the optimization problems at hand, this tutorial provides a comprehensive guide for practitioners and researchers alike.}


\maketitle

\section{Introduction}

\paragraph{What is learning to optimize (L2O)?}
Learning to Optimize (L2O) is an emerging area where traditional optimization techniques are enhanced or even replaced by machine learning models learned from data. Instead of solely relying on hand-crafted algorithms, L2O employs data to derive optimization strategies using data-driven approaches. The motivation stems from the realization that there are patterns in the way optimization problems are solved, and these patterns can be learned and generalized. 

\paragraph{Why bringing learning to optimization?}
A short answer is that machine learning provides a strong alternative to understand and \textit{utilize} the special structures of data for faster and better optimization, which can be difficult if done analytically. It can be conceptualized as two ends of the same stick:
\begin{itemize}[itemsep=0pt]
    \item On one end, consider a set that contains only one optimization problem, and it is the only single optimization problem that we are interested in. Then the most efficient method is to solve it with a solver and memorize the solution for retrieval when needed.

    \item On the other end, when we are interested in solving all instances in an optimization problem class, classic algorithms with \textit{worst-case} guarantee are the best choices. We simply select the algorithm with the best guarantee.
\end{itemize}
In real-world scenarios, however, the situation lies between these two ends. In a specific application, the optimization problems of interests often are not arbitrary ones. Instead, they share some \textit{common structures}. Effectively exploiting these special structures can help with accelerating the optimization process or generating better solutions. Sometimes, these common structures can be described in an accurate mathematical language and used for deriving better optimization algorithms.

But in other cases, such structures are too vague to characterize explicitly; hence it is hard to utilize them analytically to improve optimization.
Machine learning, especially deep learning, is adept at capturing hidden structures in data by learning from past experience. This capability is endowed by the great capacity of learning-based models and the data-driven learning process. Therefore, when combined properly with classic optimization, learning can help to achieve acceleration or superior performance.

\paragraph{When to consider L2O?}
The first prerequisite for the success of all learning-based methods, including L2O, is the access to abundant past experience, that is, in the context of optimization, \textit{many examples} of optimization problems and, often, also their solutions. 
According to the famous \textit{No Free Lunch} (NFL) theorem~\cite{wolpert1997no}, we cannot expect an L2O method to work superbly on all optimization problems. Therefore, the example problems must stand as a good representation of the problems to come in the future, or, more formally, are drawn from the same \textit{task distribution}.

There are two scenarios where L2O can improve optimization:
\begin{itemize}[itemsep=0pt, leftmargin=6.5em]
    \item[\textbf{Scenario 1}:] When we need to solve similar optimization problems repeatedly. In this scenario, the task distribution is concentrated. Because it is challenging to describe the distribution explicitly in a math language, L2O is used to find solution ``shortcuts.''
    \item[\textbf{Scenario 2}:] When it is hard to formulate a optimization model. A good example is the term ``natural image'', whose accurate mathematical characterization is eagerly needed but very difficult to write as a clean formula. In this case,  we can use L2O to find better solutions by learning an optimization model or method.
\end{itemize}

\paragraph{Train offline and then deploy.} 
Unlike traditional optimization techniques, L2O methods typically involve two phases: training and deployment. Training entails refining optimization algorithms or problem structures using data. While this training process might be time-consuming, its outcome is advantageous. After training, the derived optimization algorithms or problem frameworks tend to outperform their conventional counterparts. They either offer significantly quicker results or yield solutions of superior quality, especially when tailored to real-world scenarios.

\paragraph{Organization.} 
Within this tutorial, we categorize L2O into three paradigms, each reflecting a different level of integration between optimization and machine learning: 
\begin{itemize}
    \item Learning to accelerate optimization processes. In this approach, traditional optimization algorithms or solvers are still employed, but certain procedures are replaced with machine learning models to accelerate traditional algorithms, making this a more conservative application of machine learning to the optimization process. Example: Learning for mixed-integer programming in Section \ref{sec:mip}.
    \item Learning to generate optimization solutions. This paradigm involves using a machine learning model to directly generate solutions to optimization problems. The aim is to produce solutions far more quickly than traditional algorithms, while maintaining acceptable accuracy. Example: Algorithm unrolling in Section \ref{sec:algorithm-unrolling}.
    \item Learning to adapt optimization. In this innovative paradigm, the optimization problem itself can even be altered. Recognizing that solving optimization problems is often an approximation of the ideal solution in many real-world scenarios, machine learning models may be employed to tune these problems, making them closer to the desired real-world outcomes. Example: Plug-and-play methods in Section \ref{sec:pnp} and optimization as a layer in Section \ref{sec:opt-as-layer}.
\end{itemize}
First in Section \ref{sec:intro-deep-learning}, we introduce some preliminaries on machine learning and deep neural networks. Following this, we delve into L2O techniques for continuous optimization in Sections \ref{sec:algorithm-unrolling}, \ref{sec:pnp}, and \ref{sec:opt-as-layer}. The order of presentation is determined by its dependency on conventional optimization techniques. Methods in Section \ref{sec:algorithm-unrolling} aim at quickly finding solutions to certain optimization problems. Sections \ref{sec:pnp} and \ref{sec:opt-as-layer} are directed at producing solutions that better align with real-world needs, even if these solutions do not solve the original optimization problems. Section \ref{sec:mip} turns to mixed-integer programs. Given its distinct methodologies and relative independence from the other sections, we have positioned it after the other sections.

\paragraph{Remarks on notation.} The optimization and deep learning communities use two systems of notation that can conflict with each other and cause confusions in some cases. The notation can sometimes differ between different directions of L2O. For example, a vector $\vx$ is usually used to denote the optimization variable and thus the iterate in an iterative algorithm in the field of optimization. However, in the context of deep learning, $\vx$ mostly denotes the input to a neural network or to a layer of a network. 
In this work, we do not strive at creating a new system of notation that is consistent for both fields as it would read unnatural to the audience from either community. Instead, we try to follow the commonly accepted notations by each field when there is no risk of confusion. When there are potential conflicts, however, we will make it crystal clear the distinction between our notations and the usual ones. We humbly ask caution from the reader for a better and easier understanding of the content of this article.

\section{Introduction and Deep Neural Networks}
\label{sec:intro-deep-learning}

This section commences by providing an introductory overview of deep neural networks and how to train them, along with various adaptations and methodologies that are broadly utilized by the L2O community. We also examine prominent reinforcement learning methods that are broadly applicable in L2O. Our goal is to elucidate the fundamental concepts, advantages, and constraints of those common techniques. 

\subsection{Preliminaries of Machine Learning}
\label{sec:intro-supervised-learning}

From a mathematical perspective, many machine learning tasks are about creating a mapping derived from data. Consider the task of image classification. Here, the aim is to establish a mapping that links an image to a category, such as animals or plants. We work with a data set denoted as $\sD = \{(\vx_j,y_j)\}_{j=1}^n$, where $\vx_j$ from $\sX$ denotes an image (with $\sX$ being the image space) and $y_j \in \{0,1\}$ signifies its associated category. Here, $y_j=0$ indicates that $\vx_j$ is an animal image, while $y_j=1$ indicates it is a plant. These images can come from varied sources like photography. Their corresponding labels, $y_j$ usually added by humans, are seen as the definitive answers to the image's category. Machine learning methods typically assume $y_j$ to be the accurate label for $\vx_j$ and design a mapping that links $\vx_j$ to $y_j$ for all $j=1,2,\cdots,n$. Specifically,
\begin{itemize}
    \item According to the format and properties of $\vx$, we select a parameterized mapping $g(\cdot;\vtheta): \sX \to \{0,1\}$, there $\vtheta$ involves all the parameters.
    \item We look for desirable parameters $\vtheta^\ast$ such that $g(\vx_j;\vtheta^\ast) \approx y_j$ for all $j=1,2,\cdots,n$.
    \item With $g(\cdot;\vtheta^\ast)$, we can predict the category of a new image $\vx^\prime$ not in $\sD$ with $\hat{y} = g(\vx^\prime;\vtheta^\ast)$.
\end{itemize}
We clarify the following terminologies:
\begin{itemize}
    \item The parameterized mapping $g(\cdot;\vtheta)$ is termed a \emph{machine learning model} or, more simply, a \emph{model}.
    \item The process of finding desirable parameters $\vtheta^\ast$ is called \emph{training}. A model equipped with these optimal parameters, $g(\cdot;\vtheta^\ast)$, is named a \emph{trained model}.
    \item Utilizing the trained model to predict the category of a new image is referred to as \emph{inference}.
    \item The data set $\sD$ used during the training phase is termed the \emph{training set}.
    \item A collection of images used during the inference phase is called a \emph{testing set}.
\end{itemize}
This methodology is identified as \emph{supervised learning} since each $\vx_j$ in the training set is paired with a ground truth $y_j$ by an expert, resembling the expert supervising the training phase.
Such a supervised learning pipeline can be adapted for other tasks, where $\vx$ may represent entities other than images, and $y$ might symbolize more than just image categories. Subsequent sections will delve into broader scenarios.

\subsection{Deep Neural Networks}
\label{subsec:dnn}
Let us start by discussing the construction of a machine-learning model, with a focus on neural networks. Neural networks are made of units termed as \textit{layers}, where each layer is a basic function. Similar to function composition, a layer of neural network can also consist of several (sub-)functions. The simplest neural network is a \textit{single-layer perceptron}, a linear binary classifier that makes binary prediction by applying an affine transformation to its input $\vx\in\mathbb{R}^d$ and then thresholding the output, denoted as $g(\vx): \sR^d\to\sR, \vx\mapsto h=\sigma(\vw^\top\vx + b)$, where $h$ is the binary output of the network and $\sigma: \sR\to\{0,1\}$ is a thresholding function, or \textit{activation function}, yielding 1 if an input is positive or 0 otherwise. Here, the weight $\vw\in\mathbb{R}$ and bias $b\in\mathbb{R}$ are the \textit{parameters} of the single-layer perceptron; or we say the perceptron is \textit{parameterized} by $\vtheta=(\vw,b)$, which can be learned from data. To make it self-evident to distinguish between network input and learnable parameters, we adopt the following notation to represent a single-layer perceptron:
\begin{align}
    g(\vx; \vtheta) = \sigma(\vw^\top\vx + b).
    \label{eq:slp}
\end{align}

A single-layer perceptron can be easily extended to have multi-dimensional outputs by introducing a weight matrix $\mW\in\sR^{d \times d^\prime}$ and a bias vector $\vb\in\sR^{d^\prime}$ with $d^\prime$ being the output dimension. The extended mapping is defined as $\vx \mapsto \vh = \sigma(\mW^\top \vx + \vb)$. 
When viewed graphically, each coordinate of either $\vx$ or $\vh$ can be thought of as a node. Given a general matrix $\vW$, every coordinate of $\vx$ is linked to every coordinate of $\vh$. Therefore, such a mapping is usually called a \textit{fully-connected} (FC) layer. Note that the thresholding function in \eqref{eq:slp} is not the sole option for the activation function $\sigma$ in a fully-connected layer. Below are additional activation functions to consider:
\begin{itemize}
    \item ReLU $\sigma$: Defined as $\vx \mapsto \max(\vzero,\vx)$, where all operations are executed coordinate-wise.
    \item Sigmoid $\sigma$: Defined as $\vx \mapsto \vone / ( \vone + \exp(-\vx))$, where all operations are executed coordinate-wise.
    \item Softmax $\sigma$: For a vector $\vx$, its $i$-th component is transformed as $\sigma(\vx)_i = \exp(x_i) / (\sum_i \exp(x_i))$, where the sum runs over all components of the vector.
\end{itemize}
Multiple fully-connected layers can be stacked to form a \textit{multi-layer perceptron} (MLP). Given an input $\vx\in\sR^{d_{\mathrm{in}}}$ and a set of $L$ fully-connected layers, each of which is parameterized by $\vtheta_i = (\mW_i, \vb_i)$ and denoted as $\vg_i(\cdot;\vtheta_i): \sR^{d^{i-1}}\to\sR^{d^{i}}$ for any $i\in\{1,\dots,L\}$, an MLP is constructed as a composition function, formulated in an iterative manner as
\begin{align}
  \vh^{(0)} = \vx, \quad  \vh^{(i)} = \vg_i(\vh^{(i-1)}; \vtheta_i) := \sigma(\mW_i^\top \vh^{(i-1)} + \vb_i), ~~ i=1,\dots,L,
    \label{eq:mlp}
\end{align}
and finally returns $\vy = \vh^{(L)}$. Note that in \eqref{eq:mlp}, the dimensions of $\vW_1$ and $\vW_L$ should be compatible with $\vx,\vy$: $d_0 = d_{\mathrm{in}}$ and $d_L = d_{\mathrm{out}}$. 
The whole of learnable parameters, denoted as $\vtheta=\{\vtheta_i\}_{i=1}^L$, can be trained with data to take arbitrary values so that the resulting MLP is able to yield accurate approximations of desired outputs. The power of deep learning lies at the capacity of neural networks to approximate complicated mappings when it uses many layers, i.e., becomes deep~\cite{hornik1989multilayer}. MLPs are simple yet decently effective in many machine learning tasks including L2O. In the subsequent sections, we use the following condensed representation for the entire model as given in \eqref{eq:mlp}:
\begin{equation}
    \label{eq:mlp-simple}
    \vy = \mathrm{MLP}(\vx; \vtheta), ~~ \text{or} ~~ \vy = \vg(\vx;\vtheta).
\end{equation}

\subsection{Training a Neural Network}
\label{subsec:dnn-training}

Given initial parameters $\vtheta$, typically with random values, \textit{training} a neural network is to update these parameters so that the neural network can approximate the target unknown mapping. This process is also called learning. Training requires several components: \emph{(i)} a \textit{dataset} that we learn from; \emph{(ii)} a quantitative metric for the quality of the approximation on the dataset, called \textit{loss function}; \emph{(iii)} a mechanism that minimizes the loss function by updating the parameters; 
\emph{(iv)} regularization techniques that avoid overfitting. Step (iii) is often, though not always, based on gradient descent, and the gradient is typically, though not always, computed with the \textit{chain rule} and \textit{backpropagation}.

\paragraph{Datasets.} A neural network is trained on a finite set of data samples. Each sample consists of the input to the network $\vx\in\sR^{d_\mathrm{in}}$ and the desired output $\vy\in\sR^{d_\mathrm{out}}$. In this article, we refer to $\vx$ as the \textit{input} and $\vy$ as the \textit{label}, and denote the dataset as $\sD=\{(\vx_j, \vy_j)\}_{j=1}^N$ with $N$ being the size of the dataset. Note that the label $\vy$ is optional and can be absent in unsupervised and self-supervised tasks.

The format of a sample $(\vx, \vy)$ depends on the learning task of interest. For image classification, $\vx$ is an image and $\vy$ is a scalar indicating the class of $\vx$. For natural language translation, $\vx$ is a sequence of tokens in the source language and $\vy$ is the semantically equivalent sequence of tokens in the target language.
However, since L2O is aimed to learn how to solve optimization problems, an input to a network for L2O is what characterizes an optimization problem instance. For example, if a neural network is trained to solve least squares in the form of $\min_\vz \|\vd-\mA \vz\|^2$, an input sample is $\vx=(\vd, \mA)$ and its label is the solution $\vy=\vz^\ast=\argmin_\vz \|\vd-\mA \vz\|^2$ (assuming $\vz^\ast$ is the only minimizer).

\paragraph{Loss functions.} Given a data set $\sD$ of size $N$, the neural network to be trained $\vg(\cdot; \vtheta):\sR^{d_\mathrm{in}}\to\sR^{d_\mathrm{out}}$ is parameterized by $\vtheta\in\sR^p$ and maps each input $\vx_j$ from $\sD$ to the output $\hat{\vy}_j(\vtheta)=\vg(\vx_j;\vtheta)$. Here, $\vg(\cdot; \vtheta)$ denotes an arbitrary parameterized mapping, such as the MLP defined in \eqref{eq:mlp-simple}. A \textit{loss function} provides a metric to measure the quality of the outputs over the whole set $\{\hat{\vy}_j(\vtheta)\}_{j=1}^N$, which we need for updating the parameters in the neural network accordingly.

Considering the above mentioned supervised learning setting where each input $\vx_j$ comes with a label $\vy_j$, we first pick a metric function $\ell: \sR^{d_\mathrm{out}}\times\sR^{d_\mathrm{out}}\to\sR$ that maps a prediction-label pair $(\hat{\vy}_j, \vy_j)$ to a scalar $l_j=\ell(\hat{\vy}_j, \vy_j)$, which measures the \textit{loss} or \textit{regret} induced by the prediction.
For example, the simplest metric function, the mean-squared error (MSE), is defined by $\ell(\hat{\vy},\vy):=\|\hat{\vy}-\vy\|^2$. Cross-entropy (CE) is another commonly used metric function in machine learning.
A loss function can be defined as a function of the parameter $\vtheta$ over all samples in $\sD$:
\begin{align}
    \mathcal{L}(\vtheta) = \sum_{j=1}^N l_j = \sum_{j=1}^N \ell\big(\hat{\vy}_j(\vtheta), \vy_j\big).
    \label{eq:loss-supervised}
\end{align}
In unsupervised or self-supervised learning where the labels are absent, we can adopt a metric function whose value only depends on the output of the network with $l_i = \ell(\hat{\vy_j})$.

\paragraph{Backpropagation.}
To train the neural network, we try to find the parameter $\vtheta^\ast$ that nearly minimizes the loss function: $\vtheta^\ast \approx \argmin_{\vtheta}\ccL(\vtheta)$. Since neural networks are usually deep and thus overparameterized, the learnable parameter $\vtheta$ is extremely high-dimensional in most cases. To make the training feasible, it is most common to adopt first-order methods to minimize the loss function, which requires calculating the gradient of the loss function with respect to the parameter, i.e., $\partial\mathcal{L}/\partial\vtheta$. The chain rule in calculus is leveraged to do so efficiently. 

Consider the $i$-th layer $\vg^{(i)}(\cdot;\vtheta^{(i)})$ of the network parameterized by $\vtheta^{(i)}$. The input to the layer is $\vh^{(i-1)}$, which is the output of the previous layer, and the output is $\vh^{(i)}=\vg^{(i)}(\vh^{(i-1)}; \vtheta^{(i)})$, which is then fed into the next layer as input. Assuming we are given $\frac{\partial\mathcal{L}}{\partial\vh^{(i)}}$, the gradient of the loss function with respect to $\vh^{(i)}$, the gradient with respect to the parameter and input of this layer can be calculated with the chain rule as
\begin{align}
    \frac{\partial\mathcal{L}}{\partial\vtheta^{(i)}}
        & = \frac{\partial\mathcal{L}}{\partial\vh^{(i)}}
            \cdot \frac{\partial\vh^{(i)}}{\partial\vtheta^{(i)}},
        \label{eq:chain-rule-param} \\
    \frac{\partial\mathcal{L}}{\partial\vh^{(i-1)}}
        & = \frac{\partial\mathcal{L}}{\partial\vh^{(i)}}
            \cdot \frac{\partial\vh^{(i)}}{\partial\vh^{(i-1)}},%
        \label{eq:chain-rule-input}
\end{align}
where the latter is for propagating the gradient back to previous layers. 
The second terms on the right-hand sides of both equations depend solely on $\vg^{(i)}$. Therefore, this procedure can be performed sequentially from the output layer to the input layer. This process is called \textit{backpropagation} in contrast to the \textit{forward propagation} from inputs to outputs, as it essentially propagates the error signals in the backward direction from outputs to inputs. Modern deep learning frameworks, such as TensorFlow and PyTorch, provide out-of-the-box support of the backpropagation for most deep learning layers and modules. Practitioners do not need to code them except when developing new layers. 

\paragraph{Training algorithms.} Neural networks are mostly trained with first-order methods. For example, we can directly apply gradient descent to minimize the loss function \eqref{eq:loss-supervised} over the whole training set. However, since the training set for modern deep learning is very large, the most popular optimization approach is stochastic gradient descent (SGD), which samples a small subset of the entire training set, called a \textit{mini-batch} at each training step. Such a training method is defined in Algorithm \ref{algo:sgd}.

\begin{algorithm}[ht]
\renewcommand{\algorithmicrequire}{\textbf{Input:}}
\renewcommand{\algorithmicensure}{\textbf{Output:}}
\caption{Stochastic gradient descent in supervised learning}
\label{algo:sgd}
\begin{algorithmic}[1]
    \REQUIRE A data set $\sD=\{(\vx_j,\vy_j)\}_{j=1}^N$, a model $\vg(\cdot;\vtheta)$ with parameter $\vtheta \in \sR^p$.
    \STATE Determine hyperparameters: the number of epoch $E$, learning rate $\alpha$, mini-batch size $M$.
    \STATE Initialize the values of $\vtheta$ randomly.
    \FOR{$e = 1$ to $E$}
        \STATE Sample a subset $\hat{\sD}$ with size $M$: $\hat{\sD} \subset \sD$.
        \STATE Calculate the gradient through backpropagation: $\partial \hat{\ccL} / \partial\vtheta$, where $\hat{\ccL} := \sum_{(\vx,\vy) \in \hat{\sD}} \ell \big(\vg(\vx;\vtheta), \vy \big)$.
        \STATE Update the parameters $\vtheta \leftarrow \vtheta - \alpha  \partial \hat{\ccL} / \partial\vtheta$.
    \ENDFOR
    \RETURN $\vtheta$
\end{algorithmic}
\end{algorithm}

Motivations behind such a mini-batch method are two-fold. Firstly, effective training of heavily overparameterized neural networks requires so many training samples that the device memory cannot hold all of them at once. Secondly, a proper amount of stochasticity helps the training algorithm escape some local minima of $\ccL$ while too much of it can significantly increase the variance of gradients, slowing down convergence. Although SGD does not guarantee convergence to the global minimum $\vtheta^\ast$ due to the non-convexity nature of the loss function, it often yields a parameter $\vtheta$ corresponding to a reasonably low loss function value under properly chosen hyperparameters. (A hyperparameter tuning method will be presented in the subsequent paragraph ``Validation".)

Modifications to SGD have been proposed for more efficient and stable stochastic training, among which Adam \cite{adam} and its variants stand out. Adam updates the network parameters with the momentum, the moving average of the gradients over time, which is normalized with the variance of the gradients. Adam implicitly leverages the first-order and second-order information of the loss function and thus is faster than SGD in many applications.
SGD and Adam are enough for training neural networks used in most L2O works and are good starting points for new applications of L2O before special needs identified.

Popular training techniques, including SGD and Adam, have been implemented in most modern deep-learning platforms, eliminating the need for manual coding. In the rest of this paper, \emph{we will treat the minimization of loss functions as solvable problems}, shifting our emphasis to other critical subjects.

\paragraph{Validation.} Although it is common to manually determine the values of hyperparameters (e.g. $E,\alpha,M$ in Algorithm \ref{algo:sgd}), there is a more systematic approach known as \emph{validation}. This method partitions a complete dataset into two subsets: the training set $\sD_\mathrm{train}$ and the validation set $\sD_\mathrm{val}$. Algorithm \ref{algo:sgd} is employed on $\sD_\mathrm{train}$. The role of the validation set is to evaluate the model's capability to generalize to new data, enabling the selection of superior training hyperparameters. The validation loss function is defined as 
\[ \ccL_{\mathrm{val}}(\vtheta) := \sum_{(\vx,\vy) \in \sD_{\mathrm{val}}} \ell\Big( \vg(\vx;\vtheta), \vy \Big) \]
Imagine we are comparing two hyperparameter configurations, $(E,\alpha,M)$ and $(E',\alpha',M')$, and their resulting parameters $\vtheta(E,\alpha,M)$ and $\vtheta(E',\alpha',M')$, respectively. 
If $\ccL_{\mathrm{val}}(\vtheta(E,\alpha,M)) < \ccL_{\mathrm{val}}(\vtheta(E',\alpha',M'))$, then we conclude that $(E,\alpha,M)$ is better than $(E',\alpha',M')$. 

\paragraph{Regularization.} While the training of a neural network is conducted on the training set, the ultimate goal is to learn a network that can perform well on samples unseen during training. A neural network with this ability is considered to \textit{generalize} well or have a great \textit{generalization ability}. However, an inappropriate training process can hurt the generalization ability of a neural network by \textit{overfitting} the parameters only to the samples in the training set, resulting in poor performances on unseen samples. Regularization methods aim to address this issue by adding additional constraints or penalties to the learning process, encouraging the model to generalize better.

Here are some commonly used regularization techniques in deep learning and specifically for L2O:
\emph{(1)} \textit{L1 and L2 regularization (Weight Decay)}, which encourages the model to find simpler solutions by shrinking the weights, reducing their impact on the overall loss;
\emph{(2)} \textit{Dropout}, which combats overfitting by randomly disabling a proportion of the neurons during each training iteration. We clarify that no neurons are dropped out during inference though their outputs are scaled down by the dropout probability.
\emph{(3)} \textit{Early stopping}, which forces the training process to exit early when the model's performance on a validation set starts to degrade. Early stopping prevents the model from overfitting by finding a balance between training for too long (leading to overfitting) and stopping too early (resulting in underfitting);
\emph{(4)} \textit{Data augmentation}, which artificially increases the size of the training dataset by applying various transformations to the existing data samples so that the model can learn to be invariant to such variations, thereby improving generalization.

These regularization techniques can be used individually or in combination, depending on the specific problem and characteristics of the dataset. With regularization, deep neural networks can become more robust, generalize better to unseen data, and exhibit improved performance in real-world applications.

\subsection{Important Variants of Neural Networks}
\label{sec:other-neural-nets}

\paragraph{Recurrent neural networks.}
A \textit{recurrent neural network} (RNN) is a special class of neural networks designed for processing sequential data. Unlike standard neural networks that work with a single input with a fixed size, RNNs manage an array of sequential inputs, denoted as $\vx^{(1)}, \dots, \vx^{(T)}\in\sR^{d_\mathrm{in}}$, where the input's dimension is represented by $d_\mathrm{in}$. Given that sequence length $T$ differs between inputs, traditional MLPs fall short, as they are tailored to fixed-size inputs. This variability is commonly seen in applications like machine translation where each $\vx^{(i)}$ symbolizes a word, necessitating the handling of sentences of diverse lengths.

To cater to this need, RNNs were introduced, with the core idea of \textit{weight sharing}, that is to share the same set of parameters between layers, so that the networks can have arbitrary depths. Moreover, an RNN maintains a \textit{hidden state} vector $\vh^{(i)}\in\sR^{d_\mathrm{emb}}$ which is assumed to be a lossy representation or \textit{embedding} of the previous inputs $\vx^{(1)},\dots,\vx^{(i)}$ with $d_\mathrm{emb}$ being the dimension of the hidden embeddings. Denoting the output dimension as $d_\mathrm{out}$, a generic parameterized RNN $\vg:\sR^{d_\mathrm{in}}\times\sR^{d_\mathrm{emb}}\to\sR^{d_\mathrm{out}}\times\sR^{d_\mathrm{emb}}$ can be formulated as
\begin{align}
    \vy^{(i)}, \vh^{(i)} = \vg(\vx^{(i)}, \vh^{(i-1)}; \vtheta).
    \label{eq:rnn}
\end{align}
Here, the $i$-th layer or iteration of the RNN yields the \textit{output} of the current layer $\vy^{(i)}\in\sR^{d_\mathrm{out}}$, which can be used to calculate losses for training, and the hidden state $\vh^{(i)}$, which will be sent to the next iteration to pass the previous information. In practice, we set $\vh^{(0)}$ to all zero or a randomly sampled vector. Note that the operator $\vg$ and its learnable parameter $\vtheta$ is time-invariant and thus the network can be applied to sequences with arbitrary lengths.

However, a vanilla RNN suffers from the gradient vanishing or explosion issue when processing very long sequences~\cite{pascanu2013on}, in which case the backpropagation to early steps involves too many multiplications of the Jacobian matrix of $\vg$. While gradient explosion is more devastating for failing training, gradient vanishing happens more often and prevents the RNN from learning long-term memories as the result of early iterations receiving near zero gradients. The long short-term memory (LSTM) network was proposed with two core techniques to address this issue. One is the introduction of self-loops so the gradient can backpropagate better for longer sequences. The other is to make the self-loops dependent on the context, i.e., the input sequences so that the network can condition itself to memorize or forget certain parts of an input sequence. LSTM achieves empirical success in many fields such as audio analysis~\cite{wollmer2013lstm} and machine translation~\cite{gehring2016convolutional}. For more details of RNN and LSTM networks and advanced techniques, we refer the reader to~\cite{deeplearning}.

RNN and its variants have been widely adopted by the L2O community. The reason is that numerous optimization problems are solved with iterative algorithms, and the iterative formulation has a natural correspondence with the recurrent structure of an RNN.
Consider the gradient descent (GD) algorithm with a fixed step size that tries to minimize a differentiable objective function $f(\vx):\sR^d\to\sR$. One iteration of GD is formulated as $\vx^{(i)} = \vx^{(i-1)} - \alpha \nabla f(\vx^{(i-1)})$. If we think of the gradient $\vd^{(i-1)}=\nabla f(\vx^{(i-1)})$ as an input to the current iteration, GD can be written in a similar manner to that in \eqref{eq:rnn} as $\vx^{(i)} = \vg(\vd^{(i-1)}, \vx^{(i-1)}; \alpha)$. Here, the current iterate $\vx^{(i)}$ corresponds to the hidden state $\vh^{(i)}$ in \eqref{eq:rnn} that carries a summary of all gradients information in the past, and the iteration is parameterized by a single scalar $\alpha$, the step size in GD. Based on this similarity, many researchers construct RNNs that mimic iterative algorithms and then exploit data-driven learning in expectation of achieving certain levels of acceleration. We will revisit this relation for Algorithm Unrolling in Section~\ref{sec:algorithm-unrolling} with practical examples.

\paragraph{Convolutional neural networks.}
\textit{Convolutional neural networks} (CNNs) represent a category of neural networks that excel in processing grid-structured data. This includes 1-D grids like time-series data  and 2-D grids such as images.
While a vanilla CNN has a structure akin to an MLP as seen in \eqref{eq:mlp}, the core difference lies in its linear mapping. Instead of the generic linear mapping, $\vx \mapsto \vW^\top \vx + \vb$, CNNs deploy a unique linear mapping termed as \textit{convolution}, represented by $\vx \mapsto \vw\ast\vx + \vb$, with $\vw$ being the convolutional kernel.
For illustrative purposes, let us focus on the 1-D convolution. In this scenario, the input $\vx\in\sR^{d_\mathrm{in}}$ is a discretized variant of a 1-D function $x(t)$. Concurrently, the kernel $\vw\in\sR^K$ mirrors this, where the grid size $K$ is also referred to as the \emph{kernel size}. By denoting $\vx[i]$ as the $i$-th component of the vector $\vx$, the convolution operation, $\vw\ast\vx$, can be defined coordinate-wise as:
\begin{align}
    \label{eq:conv}
    (\vw\ast\vx)[i] = \sum_{k=1}^{K} \vw[k] \cdot \vx[i+k-1], \quad \text{for all } i = 1,2,\cdots,d_\mathrm{in}.
\end{align}
Note that the operation formulated in \eqref{eq:conv} is mathematically known as \textit{cross-correlation} but is referred to as convolution in most deep learning papers. Such a convolution can be easily extended to 2-D cases. Contrasted with the generic linear mapping characterized by parameters $(\vW,\vb)$, convolutional mapping, defined by parameters $(\vw,\vb)$, provides a more structured, parameter-efficient approach.

Convolution layers have three key characteristics that are fundamental to the success of CNNs in a wide array of applications. To begin with, when the kernel size $K$ is much smaller than the input dimension $d_\mathrm{in}$, the kernel $\vw$ only locally interacts with part of the input, which endows CNNs with sensitivity to local features. Secondly, the identical kernel is applied across all positions of the input. This design, referred to as \textit{parameter sharing}, serves as an implicit but strong regularization for better generalization and computational efficiency of CNNs. Lastly, convolutions inherently exhibit \textit{equivariance} to translations. This means that the output of a convolution mapping remains unchanged, up to the identical translation, when applied to a translated input. Consequently, this allows CNNs to effectively detect essential patterns even when inputs are shifted, which is considerably challenging for traditional MLPs where the parameters have one-on-one correspondence between input and output coordinates.

In the existing corpus of literature pertaining to L2O, CNNs predominantly find applicability in two principal use cases. The first is to incorporate convolutional layers to imitate certain operations in an optimization process, especially when the optimization involves convolution operations itself. Examples include convolutional compressive sensing, which we will explore in greater detail in Section~\ref{sec:algorithm-unrolling} dedicated for Algorithm Unrolling. The second use case leverages CNNs as black-box mappings, exploiting their benefits in specific domains. For examples, 2-D CNNs have been observed to excel at natural image denoising. Consequently, researchers have proposed the application of well-trained CNNs as black-box denoisers, intended to be integrated in a plug-and-play manner within a larger optimization framework. Section~\ref{sec:pnp} is dedicated to covering L2O methods in this category.

\paragraph{Graph neural networks.}
\label{sec:gnn}

Let us consider an undirected graph represented as $(\cV, \cE)$, where $\cV$ is the collection of vertices and $\cE$ is the collection of edges. Each node is indexed by $i$, and an edge connecting nodes $i$ and $j$ is labeled as the tuple $(i,j)$. The graph-structured data can be described as $\cG = \big( \{\vv_i\}_{i \in \cV}, \{\ve_{i,j}\}_{(i,j)\in\cE} \big)$, in which $\vv_i$ is the data attached to node $i$ and $\ve_{i,j}$ is the data attached to edge $(i,j)$. For instance, in a social network graph, each node might symbolize a user with the data $\vv_i$ giving information about them. An edge might signify a connection between two users, and the $\ve_{i,j}$ describes their relationship. It is worth noting that \emph{enormous optimization problems, especially discrete optimization, can be expressed by graph-structured data} \cite{khalil2017learning,gasse2019exact}. We defer more details of this type of approach to Section~\ref{sec:mip}. 

GNNs map graph-structured data $\cG$ to desired outputs $\vy$. Here we consider a simple one-layer message-passing GNN to present its formal expression. The first step is defined as:
\begin{align}
    \label{eq:gnn-aggregate}
    \vh_i = \mathrm{MLP}\left(  \vv_i, \sum_{j \in \cN(i)} \mathrm{MLP}(\vv_i, \vv_j, \ve_{i,j}; \vtheta_1) ; \vtheta_2  \right).
\end{align}
Here each node $i$ collects information from its neighbor $\cN(i)$, and then we invoke MLP for each node. After the operation of \eqref{eq:gnn-aggregate}, we obtain a \emph{hidden state} $\vh_i$ for each node $i$. The second step is mapping the hidden states to the output $\vy$. Based on the structure of $\vy$, there are two potential mappings. If $\vy$ illustrates the graph's entire properties, all the states ${\vh_i}$ are combined, followed by an MLP invocation for each graph. However, if $\vy$ consists of labels of all nodes formatted as $\vy = (\vy_1,\vy_2,\cdots)$, an MLP is invoked for each node individually. To elaborate:
\begin{align}
\text{Graph-level output:}~~ &  \vy = \mathrm{MLP}\left(  \sum_{i\in\cV} \vh_i ; \vtheta_3 \right) \label{eq:gnn-output-1}\\
 \text{Node-level output:}~~ & \vy_i = \mathrm{MLP}\big(\vh_i; \vtheta_4\big)  \label{eq:gnn-output-2}
\end{align}
For the rest of this paper, the complete GNN model is condensed to:
\[
\begin{aligned}
    \vy = & \ \mathrm{GNN}(\cG;\vtheta) \quad \text{Graph-level GNN: \eqref{eq:gnn-aggregate}, \eqref{eq:gnn-output-1}}\\
    \vy_i = & \ \mathrm{GNN}(i,\cG;\vtheta) \quad \text{Node-level GNN: \eqref{eq:gnn-aggregate}, \eqref{eq:gnn-output-2}}
\end{aligned}
\]
For a graph-level GNN, the parameters are given by $\vtheta = (\vtheta_1,\vtheta_2,\vtheta_3)$, while for a node-level GNN, they are $\vtheta = (\vtheta_1,\vtheta_2,\vtheta_4)$.

The GNN structure outlined earlier presents multiple advantages when dealing with graph-structured data. Firstly, \textbf{scalability} is inherent in its design. It is crucial to understand that every MLP employed is consistent across all nodes and edges. As a result, even when a new graph of varying size emerges, the methodology to determine hidden states remains unaltered. Moreover, the same repetitive process can be enacted on this new graph without modifying the MLPs. This implies that this particular message-passing GNN is adaptable to graphs of any dimension. Secondly, the structure ensures \textbf{permutation invariance}. Evaluating \eqref{eq:gnn-aggregate}, \eqref{eq:gnn-output-1}, and \eqref{eq:gnn-output-2}, it becomes apparent that altering the indices of two nodes does not impact the output of a graph-level GNN. On the other hand, for a node-level GNN, the outputs adjust corresponding to the permutation of inputs. Refer to \cite{jegelka2022gnn,wu2022gnn} for more properties and applications of GNN.

\subsection{Reinforcement Learning} 
\label{sec:rl}

In addition to supervised learning discussed in Section \ref{sec:intro-supervised-learning}, another important machine learning approach is \emph{reinforcement learning (RL)}. Instead of relying on vast amounts of data with manually labeled ground-truth as in supervised learning, RL focuses on learning a decision-making policy through interactions with the underlying system. This subsection introduces basic concepts of RL and provides a concise guide for beginners.

\paragraph{Markov decision process (MDP).} At the heart of reinforcement learning lies the Markov decision process. Imagine playing a game where at every step, based on your current position, you choose an action. This action leads to a specific outcome and potentially a reward. The goal is to maximize rewards over time, and this setup is described by an MDP.

Formally, an MDP is a discrete-time stochastic process during time $t=0,1,\cdots,T$. At each time $t$, we have three random variables $\vs_t,\va_t,r_t$. 
The state variable $\vs_t$ reflects the system's condition, like a chess board's layout or a car's status (position, speed, etc.). 
The action variable $\va_t$ represents the decision at time $t$.
The reward variable $r_t$ is a value (typically real) given by the system or set by individuals.
An MDP is characterized by a tuple $(\sS,\sA,p,r)$:
\begin{itemize}
    \item $\sS$ denotes the \emph{state space}, encompassing all possible values of a state variable.
    \item $\sA$ denotes the \emph{action space}, encompassing all possible values of an action variable.
    \item $p:\sS \times \sS \times \sA \to \sR$ defines the probability of \emph{transitioning} from state $\vs$ to $\vs'$ after taking action $\va$: 
    \[\mathrm{Pr}(\vs_{t+1}=\vs'|\vs_t=\vs,\va_t = \va) = p(\vs,\vs',\va)~~\text{for all }\vs,\vs^\prime \in \sS, \va \in \sA, t = 0,1,2,\cdots,T.\]
    \item $r:\sS \times \sA \to \sR$ describes the \emph{reward} received after taking action $\va$ at state $\vs$.
\end{itemize}
Once $(\sS,\sA,p,r)$ are specified, the MDP is fully defined. 
For every $(\vs,\va)$, if $p(\vs,\vs',\va) = 1$ for some $\vs'$, the transition is deterministic. 
Depending on the availability of $p$, the MDP can either be a conventional control problem (known $p$) or a reinforcement learning problem (unknown $p$).

\paragraph{Reinforcement learning.} Imagine that we cannot directly access the formula for $p$ and view the MDP as a ``black-box." Our interaction with this black-box is iterative: we observe state $\vs_t$ from it, decide to take action $\va_t$, and then receive a reward $r_t$ and a new state $\vs_{t+1}$ in response. Given these interactions, the challenge is: Can we deduce an optimal policy that gives the probabilities of taking different actions for a particular state? This is the fundamental query addressed by RL. To elucidate, let us define some RL terms:
\begin{itemize}
    \item Environment: The black-box MDP $(\sS,\sA,p,r)$ where $p$ is not directly accessible.
    \item Policy: A parameterized model $\pi(\va,\vs; \vtheta)$ that indicates the probability of selecting action $\va$ when in state $\vs$. Mathematically, $\mathrm{Pr}(\va_t=\va|\vs_t=\vs) = \pi(\va,\vs; \vtheta)$ for all $\va \in \sA,\vs \in \sS,t=0,1,\cdots,T$.
    \item Trial: With a given environment and policy, the distribution of random variables $\vs_t,\va_t$ are well-defined. Concatenating all these variables together, it is denoted as $\vtau := \{(\vs_t,\va_t)\}_{t=0}^T$. The joint probability distribution of $\vtau$, induced by the MDP $(\sS,\sA,p,r)$ and policy $\pi(\va,\vs; \vtheta)$, is denoted with $P_{\vtheta}(\vtau)$ and is given by:
    \begin{equation}
        \label{eq:joint-dist}
        P_{\vtheta}(\vtau) := \mathrm{Pr}(\vs_0,\va_0,\cdots,\vs_T,\va_T) = \mathrm{Pr}(\vs_0) \left(\prod_{t=0}^{T-1} \mathrm{Pr}(\va_t|\vs_t) \mathrm{Pr}(\vs_{t+1}|\vs_t,\va_t) \right) \mathrm{Pr}(\va_T|\vs_T)
    \end{equation}
    A realization of $\vtau$ is called a \emph{trial}. A trial manifests when the policy interacts with the environment, logging the resulting data.
\end{itemize}
The goal of RL is to minimize the loss function defined by the expected cumulative reward over time:
\begin{equation}
    \label{eq:reward-finite}
   \Min_{\vtheta}~ \ccL_T(\vtheta) := - \EE_{\vtau \sim P_{\vtheta}(\vtau)} \sum_{t=1}^T r(\vs_t,\va_t).
\end{equation}

\paragraph{Policy gradient.} The policy gradient is a technique of computing the gradient $\partial \ccL_T / \partial \vtheta$, allowing the use of first-order optimization methods like SGD or Adam, as discussed in Section \ref{subsec:dnn-training}. Specifically, the gradient can be expressed as:
\[
\begin{aligned}
   \nabla \ccL_T(\vtheta) = & \nabla \int_{\vtau \in (\sS \times \sA)^{T+1}} P_{\vtheta}(\vtau) \left( \sum_{t=1}^T r(\vs_t,\va_t) \right) d \vtau \\
    = & \int_{\vtau \in (\sS \times \sA)^{T+1}} \frac{\partial}{\partial \vtheta} P_{\vtheta}(\vtau) \left( \sum_{t=1}^T r(\vs_t,\va_t) \right) d \vtau \\
    = & \int_{\vtau \in (\sS \times \sA)^{T+1}} P_{\vtheta}(\vtau) \frac{\partial}{\partial \vtheta} \log\big( P_{\vtheta}(\vtau)\big) \left( \sum_{t=1}^T r(\vs_t,\va_t) \right) d \vtau \\
    = & \EE_{\vtau \sim P_{\vtheta}(\vtau)} \left[ \frac{\partial}{\partial \vtheta} \log\big( P_{\vtheta}(\vtau)\big) \left( \sum_{t=1}^T r(\vs_t,\va_t) \right) \right]
\end{aligned}
\]
Using the joint distribution formula:
\[
\log\big( P_{\vtheta}(\vtau)\big) = \log\big(\mathrm{Pr}(\vs_0)\big) + \left( \sum_{t=0}^{T-1} \log\big(\mathrm{Pr}(\va_t|\vs_t) \big) + \log\big(\mathrm{Pr}(\vs_{t+1}|\vs_t,\va_t)\big) \right) + \log\big( \mathrm{Pr}(\va_T|\vs_T) \big)
\]
Note that $\mathrm{Pr}(\vs_{t+1}|\vs_t,\va_t)$ and $\mathrm{Pr}(\vs_0)$ are independent of $\vtheta$. Therefore,
\[ 
\frac{\partial}{\partial \vtheta} \log\big( P_{\vtheta}(\vtau)\big) = \sum_{t=0}^{T} \frac{\partial}{\partial \vtheta} \log\big(\mathrm{Pr}(\va_t|\vs_t) \big) = \sum_{t=0}^{T} \frac{\partial}{\partial \vtheta} \log\big( \pi(\va_t,\vs_t;\vtheta) \big) 
\]
Typically, a deep neural network is selected as the policy model $\pi$. Thus, in the above equation, $\partial \log (\pi) / \partial \vtheta$ be computed using backpropagation, as mentioned in Section \ref{subsec:dnn-training}. To solve \eqref{eq:reward-finite}, one can employ Algorithm \ref{algo:sgd}: interact with the environment $M$ times, record a trial $\vtau_i$ each time, consider a set of trials $\{\vtau_i\}_{i=1}^M$ a mini-batch in Algorithm \ref{algo:sgd}, and update $\vtheta$ similarly to Algorithm \ref{algo:sgd}.

\paragraph{Q-learning.} Let us examine the scenario where $T=\infty$. To ensure that the cumulative reward is not infinite, we introduce a \emph{discount factor} $\gamma \in (0,1)$ and define a new loss function as
\[
    \ccL_\gamma(\pi) := - \EE_{\vtau \sim P_{\pi}(\vtau)} \sum_{t=1}^{\infty} \gamma^t r(\vs_t,\va_t).
\]
Instead of the $P_{\vtheta}$ in \eqref{eq:reward-finite}, here we use a more general joint distribution $P_\pi$, where $\pi$ represents any potential policy with 
$\pi(\va,\vs) = \mathrm{Pr}(\va_t=\va|\vs_t=\vs)$. Based on this, we define a Q-function $Q_\pi:\sA\times\sS \to \sR$ as
\[ Q_\pi(\va,\vs) := \EE_{\vtau \sim P_{\pi}(\vtau)} \left[ \sum_{t'=t}^{\infty} \gamma^{t'-t} r(\vs_{t'},\va_{t'}) \bigg | \vs_t = \vs, \va_t = \va \right] \]
The Q-function measures the expected reward for the action $\va_t=\va$ as $\vs_t=\vs$. Given the nature of MDPs, this function remains time-independent. Hence, the Q-function can act as an evaluation metric for policy $\pi$: the higher its value, the better the policy. The best policy's Q-function, denoted by $Q_\ast = \max_{\pi} Q_{\pi}$, should satisfy the Bellman equation:
\begin{equation}
    \label{eq:bellmann}
    Q_\ast(\va,\vs) = \EE \left[ r(\vs,\va) + \gamma \max_{\va'\in\sA} Q_\ast(\vs_{t+1},\va') \bigg | \vs_t = \vs, \va_t = \va \right]
\end{equation}
The expectation here relies only on the transition probability $\mathrm{Pr}(\vs_{t+1}=\vs'|\vs_t=\vs,\va_t=\va)$, and is independent of the policy $\pi$. If the optimal Q-function $Q_\ast$ is available, an optimal deterministic policy is easy to derive as $\va_t = \max_{\va \in \sA} Q_\ast(\va,\vs_t)$, assuming the size of $\sA$ is moderate. Consequently, Q-learning's primary objective is to \emph{identify a Q-function that aligns with the Bellman equation}. 

For Q-learning, we adopt a parameterized machine-learning model to represent the Q-function $Q(\va,\vs;\vtheta)$. 
Initially, we run the MDP, gathering data in tuples like $(\vs,\va,r,\vs')$. Here, $\vs$ represents any state during the current run, $\va$ is the action suggested by $\argmax_{\va} Q(\va,\vs;\vtheta)$, $r = r(\vs,\va)$ is the reward value, and $\vs'$ represents the next state. Multiple such data tuples can be extracted from a single run. After sufficient runs and ample data collection, the objective becomes to train the Q-function, so that it reflects the Bellman equation. This is achieved by minimizing the loss:
\begin{equation}
    \label{eq:bellmann-data}
    \Min_{\vtheta} \ccL_{Q}(\vtheta):= \sum_{(\vs,\va,r,\vs')} \big( y - Q(\va,\vs;\vtheta) \big)^2, ~~\text{where } y = r + \max_{\va' \in \sA} Q(\va',\vs';\vtheta_{\mathrm{prev}})
\end{equation}
In this scenario, we update $\vtheta$ while keeping $\vtheta_{\mathrm{prev}}$ constant. Through backpropagation, we can compute the gradient of $\ccL_{Q}$, enabling us to apply Algorithm \ref{algo:sgd} to update $\vtheta$. After adequate updates of $\vtheta$, we synchronize $\vtheta_{\mathrm{prev}}$ and $\vtheta$. This cycle of updating $\vtheta$ continues until the computation budget is reached or $\vtheta_{\mathrm{prev}}$ and $\vtheta$ converge to the same value.

\paragraph{A start-up recipe.} It is essential to understand that many of the widely adopted reinforcement learning algorithms, like the policy gradient method and Q-learning, are already accessible through modern deep learning platforms. Hence, there's no pressing need to build these algorithms from scratch. When integrating these algorithms into your projects, the primary elements one must define in their code include:
\begin{itemize}
    \item A state space
    \item An action space
    \item A reward function
\end{itemize}
Based on these components, one also configure:
\begin{itemize}
    \item A stochastic policy model $\pi(\va,\vs;\vtheta)$ for the policy gradient method, or
    \item A parameterized Q function $Q(\va,\vs;\vtheta)$ for Q-learning.
\end{itemize}

\section{Algorithm Unrolling}
\label{sec:algorithm-unrolling}

In this section, we present an important branch of L2O methods called \textit{algorithm unrolling}, which expresses a classic iterative algorithm as a neural network of a certain depth and replaces certain components of the algorithm with learnable parameters. We start with a simple example of unrolling a projected gradient descent algorithm (PGD) into a neural network. Consider a non-negative least square optimization problem
\begin{align}
    \label{eq:least-square}
    \min_\vx f(\vx) = \frac{1}{2} \| \vd - \mA \vx \|_2^2 \quad \text{subject to} \quad \vx \geq \vzero,
\end{align}
where $\mA\in\sR^{m\times n}$ is given $\vd\in\sR^{n}$ is the observation we have for the regression. We say $\vd$ is the input to the algorithm, which characterizes an instance of non-negative least squares along with $\mA$.
PGD for \eqref{eq:least-square} is the iteration
\begin{align}
    \label{eq:gd-lsq}
    \vx^{(i)} = P_+\left(\vx^{(i-1)} - \alpha \nabla f(\vx^{(i-1)})\right),\quad i=1,2,\dots,
\end{align}
where $P_+(\cdot)$ projects its input to the non-negativity orthant, $\nabla f(\vx) = \mA^\top(\mA\vx - \vd)$, and $\alpha$ is the step size. $P_+(\cdot)$ coincides with 
the RELU activation function, which is ubiquitously used in deep learning.
If we treat the iterate $\vx^{(i)}$ as a state vector of the $i$-th iteration and $(\alpha, \mA)$ as the parameters,
with some rearrangement, we can re-write \eqref{eq:gd-lsq} as a recurrent system
\begin{align}
    \label{eq:pgd-lsq-recurrent}
    \vx^{(i)}  = P_+\left((\alpha\mA^\top)\vd + (\mI_n - \alpha \mA^\top\mA) \vx^{(i-1)}\right) =: \vg\Big(\vd, \vx^{(i-1)}; (\alpha, \mA)\Big),
\end{align}
where $\mI_n\in\sR^{n\times n}$ is the identity matrix and $=:$ defines its right-hand side by its left-hand side. We can see the correspondence between PGD~\eqref{eq:pgd-lsq-recurrent} and RNN~\eqref{eq:rnn}, with $P_+(\cdot)$ being the activation function and $\vd$ the constant input to each iteration.%
\footnote{
    Note there is a difference in notation, where the state vector is denoted as $\vx^{(i)}$ in~\eqref{eq:pgd-lsq-recurrent} while it is denoted as $\vh^{(i)}$ in~\eqref{eq:rnn}.}
This similarity essentially motivates the researchers to actually ``relax'' the algorithm into an RNN by \textit{parameterizing} part of the algorithm to be learnable with data, e.g., $(\alpha, \mA)$ in \eqref{eq:pgd-lsq-recurrent}.
Note that there are many ways of parameterization other than simply plain conversion. We will discuss in more details later in this section.

\begin{figure}[t]
    \centering
    \begin{tabular}{cc}
        \includegraphics[width=0.46\textwidth]{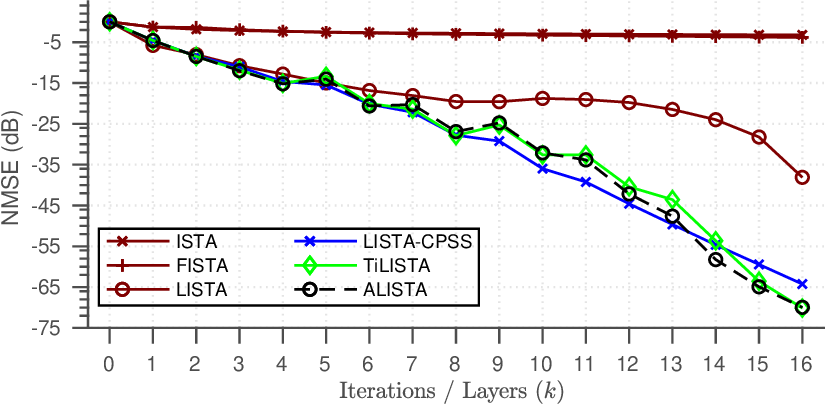}
        &
        \includegraphics[width=0.46\textwidth]{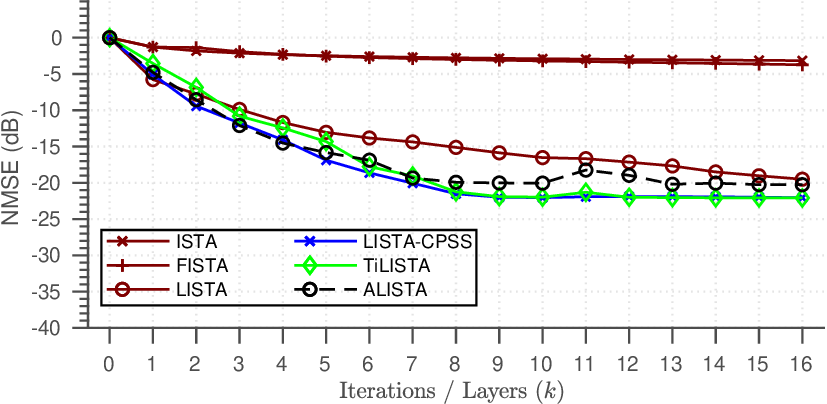}
    \end{tabular}
    \caption{Impressive acceleration can be achieved by algorithm unrolling methods (LISTA, LISTA-CPSS, TiLISTA and ALISTA) compared to the classic iterative algorithm ISTA and its variant FISTA accelerated with momentum. Algorithm unrolling uses orders of magnitude fewer iterations than ISTA/FISTA to achieve the same precision.
    Left figure: noiseless case. Right figure: noisy case (SNR = 20). X-axis is the number of iterations; Y-axis is the normalized mean squared error (lower is better). \textit{Plot source: Figure~1 of~\cite{alista}.}}
    \label{fig:trailer-unrolling}
\end{figure}

Once the parameterization is done, we are left with a neural network that is expected to achieve acceleration on a certain type of optimization problems after proper data-driven training.
The original motivation of algorithm unrolling~\cite{lista} is that the conversion from an iterative algorithm to a neural network enables data-driven learning for solving a specific type of optimization problem. At the cost of convergence guarantee for worst cases, the resulting networks are observed to accelerate when applied to \textbf{unseen} problems of the same type, as shown in Figure~\ref{fig:trailer-unrolling}.
Since its emergence, algorithm unrolling has found broad applications, including but not limited to (natural) image restoration and enhancement~\cite{zhang2017leanring}, medical and biological imaging~\cite{li2021deep}, and wireless communication~\cite{balatsoukas2019deep}. Algorithm unrolling is also referred to as \textit{algorithm unfolding} in the literature.

In this section, we strive to present a comprehensive guideline of how to unroll an algorithm in hand, the major design elements to consider and useful techniques in practice. We first set up the basics and notations for algorithm unrolling and introduce the key motivations behind it. Then we use learned ISTA (LISTA)~\cite{lista} for a case study to lay out the general routine of unrolling an iterative algorithm into a neural network. After discussions on the key design elements and advanced techniques for algorithm unrolling, we will conclude this section with the introduction to mathematical analyses behind algorithm unrolling in the existing literature.

\subsection{Basics and motivations}
\label{subsec:unrolling-basics}

\paragraph{Generic formulation.}
We consider optimization problems in the form of  
\[ \Min_{\vx \in \sX(\vd)} f(\vx; \vd) \]
where $\vx$ is the  the variable to optimize and $\vd$ characterizes an instance of such type of problems. We first select a general iterative algorithm that targets regressing the optimal solution from the observation $\vd$ with the updating form
\begin{equation}
    \label{eq:generic-algorithm}
    \vx^{(i)} = \vg(\vd, \vx^{(i-1)}; \vtheta), \quad i = 1,2,\cdots,\infty.
\end{equation}
In this case, we say $\vd$ is the \textit{input} to the algorithm
and $\vtheta$ is the parameter of the algorithm that consists of the knowledge of the optimization problem and the hyperparameters selected by users. We hope that $\vx^{(i)}(\vd)\to\vx^\star(\vd)$ where $\vx^\star(\vd):=\argmin_{\vx \in \sX(\vd)}f(\vx;\vd)$ for all $\vd$ of interest.

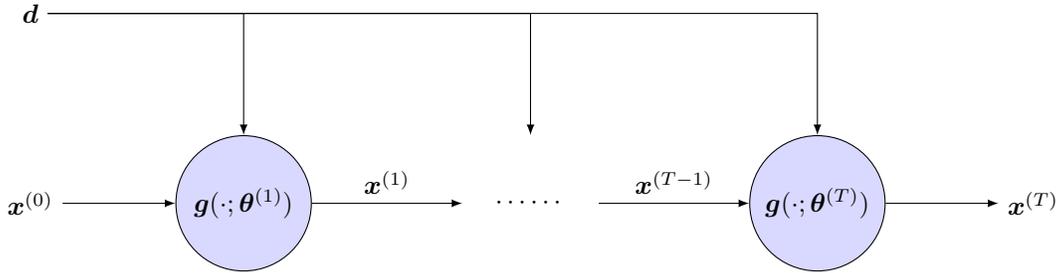
\begin{figure}[t]
    \centering
    \input{unrolling-general.tikz}
    \caption{An illustration of a generic unrolling process that turns an iterative algorithm into truncated neural network. The observation $\vd$ can be seen as the constant input to all iterations. The parameter $\vtheta^{(i)}$ can be either temporally varying or constant, corresponding to a recurrent or feedforward system, respectively.}
    \label{fig:unrolling-general}
\end{figure}

Algorithm unrolling unfolds the optimization iteration~\eqref{eq:generic-algorithm} into a neural network of a certain depth and replacing certain parameters of the algorithm, i.e., $\vtheta$ in~\eqref{eq:generic-algorithm}, into learnable parameters. The new updates take the form
\begin{equation}
    \label{eq:model-unroll}
    \vx^{(i)} = \vg(\vd, \vx^{(i-1)}; \vtheta^{(i)}), \quad i = 1,2,\cdots,T,
\end{equation}
where $T$ is the depth of the neural network, and $\vtheta^{(i)}$ are the learnable parameters in the $i$-th layer.
If we use depth-invariant parameters, that is $\vtheta^{(i)}\equiv\vtheta$ for all $i$, we will get an RNN, which is the case of the first network derived using algorithm unrolling in~\cite{lista}. We also say that the parameter $\vtheta$ is \textit{shared} or \textit{tied} across layers in this case.
Many works that follow~\cite{lista} propose to use time-variant parameters for a larger model capacity and update them independently during training~\cite{borgerding2017amp,chen2018theoretical}. Other works also show the benefits of a hybrid parameterization scheme that partially share parameters over time~\cite{alista,aberdam2021adalista}. We will discuss about this design choice in Section~\ref{subsec:unrolling-design-factors}.
The process of unrolling an iterative algorithm into a neural network is illustrated in Figure~\ref{fig:unrolling-general}.

The next step is to train the resulting neural network on a training set $\sD_\mathrm{train}$ of size $N$. Consider a simple supervised learning setting where the target signal is accessible by users. Each sample in the training set is an input-label pair $(\vd_j, \vx^\ast_j)\in \sD_\mathrm{train}$. Denote the resulting network as a end-to-end mapping $G:\sR^m\to\sR^n$ that maps input $\vd_j$ to output $\vx_j^{(T)}=G(\vd_j;\vTheta)$ where $\vTheta =\{\vtheta^{(i)}\}_{i=1}^T$.
With a proper metric function $\ell:\sR^n\times\sR^n\to\sR, (\vx_j^{(T)}, \vx_j^\ast)\mapsto l_j=\ell(\vx_j^{(T)}, \vx_j^\ast)$, we train $\vtheta$ by minimizing a loss function
\begin{equation}
    \label{eq:unrolling-train-loss}
    \min_\vtheta \mathcal{L}(\vTheta) = \sum_{j=1}^{N} c_j \cdot l_j = \sum_{j=1}^{N} c_j \cdot \ell\left( G(\vd_j;\vTheta),  \vx^\ast_j \right),
\end{equation}
where $c_j\geq0$ is the coefficient of the weighted sum with $\sum_{j=1}^N c_j = 1$. When $\vx^\ast_j$ is not available, we can adopt an unsupervised metric function such as the objective function itself.

\paragraph{Motivations.}
Introduced to bridge the gap between classic optimization methods and black-box deep learning approaches, algorithm unrolling combines benefits from both domains. By unrolling an iterative optimization algorithm into a fixed-depth neural network, it can be trained to solve a specific distribution of optimization problems, often achieving significant acceleration over the original algorithm on similar problems. This learning ability allows the unrolled algorithm to adapt to the problem class and learn optimal parameters for faster convergence.

At the same time, algorithm unrolling inherits the mathematical structures and domain knowledge embedded in the classic optimization methods developed by experts. This inherited knowledge acts as a strong regularization for the neural networks, constraining the search space and guiding the learning process. Hence, the unrolled algorithms can be trained more efficiently, requiring fewer training samples compared to generic deep neural networks.

\paragraph{Why can algorithm unrolling be effective?}
We consider algorithm unrolling to be effective if the resulting network can, compared to the original algorithm, accelerate convergence to a desired solution or generate a solution with better quality. The key motivations for the effectiveness are two-fold.

In one scenario, the optimization problems in a specific application are usually drawn from a concentrated distribution.
For example, the authors of~\cite{lista} studied sparse-coding-based image reconstruction on the MNIST dataset~\cite{mnist}. A common set of basis vectors is used for all images from MNIST; hence the reconstruction problems are the same but have different input images. The classic algorithms for these problems such as gradient-projection iterations work for a much larger class of problems and are not tailored to the specific reconstruction problems where the images to reconstruct share a common linear basis.
Algorithm unrolling can learn from the training data about the commonality and yield a neural network that takes a fast ``shortcut'' to the solution not known by human experts.

The other scenario arises where the optimization formulations derived by human experts fail to fully capture the intricacies and complexities of real-world applications. Consequently, the exact minimizer of such an optimization problem may yield suboptimal solutions. A prime example of this limitation is the use of total variation regularization in image processing, which removes noise in flat regions while preserving sharp edges but also introduces undesirable staircase artifacts in the resulting image. Total variation regularization and subsequent improvements only partially model the inherent properties of natural images. By applying algorithm unrolling, which significantly relaxes the original algorithm, we can discover effective corrections to the original formulation's algorithms and thus capture the nuances of data-specific information, which in turn yield high-quality solutions in practical applications, overcoming the limitations of traditional, analytic formulations.

\subsection{A recipe for algorithm unrolling --- Learned ISTA as an example}
\label{subsec:unrolling-recipe}

Given a specific task, we lay out a general outline of unrolling an iterative algorithm into a neural network for solving a specific type of optimization problems for the task:

\begin{itemize}[leftmargin=1.8cm]
    \item[\textbf{Step 1:}] Set up a proper modeling of the task as an optimization problem and identify a classic iterative algorithm for solving the modeled optimization.
    \item[\textbf{Step 2:}] Unroll the iterative algorithm over time into the form of a recurrent system as in~\eqref{eq:pgd-lsq-recurrent} and truncate it a finite number of steps.
    \item[\textbf{Step 3:}] Parameterize the unrolled recurrent system, that is to make part of the recurrent update form learnable by directly replacing them with learnable parameters or generating them with a neural network with learnable parameters.
    \item[\textbf{Step 4:}] Train the unrolled and parameterized network with data.
\end{itemize}

The first and second steps are intuitive in most cases, especially when the original goal is to accelerate a specified algorithm for a specific task, such as accelerating ISTA for sparse coding problems in~\cite{lista} and projected gradient descent for Multiple-Input-Multiple-Output (MIMO) detection in~\cite{samuel2019learning}.

The fourth step can play an important role in the empirical success of an neural network derived with algorithm unrolling. To train the network derived effectively, one needs a proper training objective or loss function that is suitable for the task of interest, an effective training algorithm with necessary regularization, and other training techniques that helps to stabilize or accelerate the training process. We will present the basic training scheme and part of the techniques via the case study in this subsection and have a more detailed discussion in Section~\ref{subsec:unrolling-design-factors}.

The third step, i.e., the parameterization step, is the essential yet the most challenging step in the methodology of algorithm unrolling. A proper parameterization is fundamental to the easiness of training later and to the empirical performance and acceleration that can be expected during inference, but at the same time, requires careful consideration of a number of design elements. These design elements directly influence the architecture of the resulting network and thus its empirical performance. There is no one ``correct'' way of parameterization, but it should be adaptive to the task of interest, 
including its objective, data format, and common structures shared among the data. 

For example, \textit{\textbf{i}terative \textbf{s}hrinkage \textbf{t}hresholding \textbf{a}lgorithm} (ISTA), which was first unrolled in~\cite{lista}, now has more than ten versions with different parameterizations and the resulting neural networks. These variants have been proposed for varying purposes with adaptation for specific tasks. Some variants have special parameter representation for improved trainability and interpretability~\cite{moreau2017understanding,chen2018theoretical,ablin2019learning}; some variants introduce additional structures or modules so that the resulting network can generalize to a larger range of data without re-training~\cite{behrens2020neurally,aberdam2021adalista,chen2021hyperparameter}; and other variants replace certain parts of the unrolled algorithms with more advanced deep learning modules~\cite{zhang2018istanet}.
For its richness in the literature, ISTA is a perfect example for algorithm unrolling to explain the multiple aspects of network design.

\paragraph{A case study with Learned ISTA.} We adopt the original \textit{learned ISTA} (LISTA) network~\cite{lista} as a study case to go over the four steps laid out in the recipe.
We choose LISTA~\cite{lista} as the study case example because it is the first algorithm unrolling work and the most impactful one in the sense that the algorithm unrolled in~\cite{lista}, ISTA, has been investigated most frequently in the literature.

\paragraph{Step 1 - Set up the optimization model and algorithm.}
The first step is to set up the optimization model and identify a classic iterative algorithm for solving it. In~\cite{lista} the authors strove at learning a fast approximation to sparse coding problems by unrolling ISTA. Given an observation $\vd\in\sR^m$ as input, the sparse coding problem is usually formulated as finding the optimal sparse vector $\vx^\ast\in\sR^n$ that minimizes the \textit{LASSO} objective function, which is a weighted sum of the squared reconstruction error and an $\ell_1$ sparsity regularization, formulated as
\begin{align}
    \label{eq:lasso}
    \mathrm{Find} ~~ \vx^\ast = \argmin_\vx \frac{1}{2}\|\vd-\mA\vx\|_2^2 + \lambda\|\vx\|_1.
\end{align}
Here, $\mA\in\sR^{m\times n}$ is the given dictionary that is often overcomplete ($m<n$) and $\lambda$ is a hyperparameter that controls the regularization strength, which is usually tuned by hand on a set of sparse coding problems. ISTA algorithm is an iterative algorithm for solving~\eqref{eq:lasso}. In each iteration, ISTA essentially takes a gradient descent step with respect to the squared reconstruction error term with a step size of $\alpha$ and then shrink the non-zero entries towards zero through a so-called soft-thresholding function $\eta(\cdot; \rho)$ with a positive threshold $\rho$
\begin{align}
    \label{eq:ista}
    \text{Classic ISTA~\cite{fista}:} \quad
    & \vx^{(i)} = \eta\left( \vx^{(i-1)} - \alpha\mA^\top(\mA\vx^{(i-1)} - \vd); \rho \right), ~ i = 1, \dots, \infty \\
    & \vx^{(0)} = \bm{0}, \nonumber
\end{align}
where $\eta(z;\rho)=\mathrm{sign}(z)(|z|-\rho)_+$ is called \textit{soft-thresholding function} which is applied to the input vector in a coordinatewise manner.
With a proper selection of the step size and the threshold, ISTA algorithm is provably converging to $\vx^\ast$ at a sublinear speed of $\mathcal{O}(1/i)$.
\footnote{
    Typically for ISTA, we set the step size $\alpha=1/L$ where $L$ is the largest eigenvalue of $\mA\mA^\top$ and threshold $\rho=\lambda/L$.
}

\paragraph{Step 2 - Unroll the iterative algorithm and truncate.}
In~\cite{lista}, the authors rearranged and merged terms in~\eqref{eq:ista} and truncate the iteration to a fixed number steps, say $T$, turning the algorithm into the following form
\begin{align}
    \label{eq:ista-unrolled-truncated}
    \text{LISTA~\cite{lista}:} \quad
    & \vx^{(i)} = \eta\left( \mW_1\vx^{(i-1)} + \mW_2\vd; \rho \right), ~ i = 1, \dots, T \\
    & \mW_1 = \mI_n - \alpha \mA^\top\mA, \quad \mW_2 = \alpha\mA^\top, \quad \vx^{(0)} = \bm{0}. \nonumber
\end{align}

\paragraph{Step 3 - Parameterize the unrolled system into a neural network.}
With the unrolled and truncated system in~\eqref{eq:ista-unrolled-truncated}, the authors in~\cite{lista} directly converted the $\mW_1$, $\mW_2$ matrices and the threshold $\rho$ into learnable parameters that are shared across all layers, thus turning ISTA into a recurrent neural network. We denote the learnable parameters of the resulting RNN as the concatenation of the two matrices and the threshold $\vtheta=(\mW_1, \mW_2, \rho)$. Note that while we allow the parameters to be updated through data-driven training, we can still initialize them in the same way as we set them in standard ISTA. This initialization scheme provides a good starting point that is already a good approximation to ISTA in the first place, which can help the training process in the next step.

Besides directly converting $\mW_1$, $\mW_2$ and $\rho$ into learnable parameters and sharing them across layers, there are many other ways of parameterization resulting in different variants of LISTA. For example, a simple extension to the original recurrent version of LISTA in~\cite{lista} is to \textit{untie} $\mW_1$, $\mW_2$ and $\rho$ and instead learn layer-dependent parameters $(\mW_1^{(i)}, \mW_2^{(i)}, \rho^{(i)})$ for all $ i\in\{1,\dots,T\}$. This untied parameterization strategy is widely adopted in following works such as improved LISTA for sparse inverse recovery~\cite{chen2018theoretical} and algorithm unrolling for other algorithms~\cite{borgerding2017amp}. The main benefit of this strategy is that it radically increases the number of learnable parameters in the resulting network, lifting the \textit{complexity of parameterization} from $\mathcal{C}(m(m+n))$ in~\cite{lista} to $\mathcal{C}(Tm(m+n))$.%
\footnote{
    We use a different notation of $\mathcal{C}(\cdot)$ to represent the complexity of parameterization to distinguish it from the time complexity of computation for an algorithm or a neural network, which is denoted with $\mathcal{O}(\cdot)$. 
}
The increased complexity in parameterization enables the neural network to learn more complicated mappings and thus expect potentially better acceleration, but it also makes the training process more difficult and requires special techniques for a stabilized training process to actually achieve that expected acceleration. We will have a thorough discussion in the next subsection about the trade-offs between various parameterization strategies.

\paragraph{Step 4 - Train the resulting neural network.}
To set up a proper training process we first need to generate a training set. In~\cite{lista}, the authors assumed that all sparse coding problems of their interest are \textit{drawn from the same distribution}. To be specific, the inputs, $\{\vd_j\}_{j=1}^N$, are $N$ image patches extracted from one dataset. All sparse coding problems share the same dictionary matrix $\mA$, which is learned beforehand on the same dataset with an online dictionary learning algorithm.

As we elaborated in \textbf{Step~1}, the objective of the original LISTA is to learn a fast approximation to optimal solutions of sparse coding problems, i.e., the minimizers of the LASSO objective $\vx^\ast$ in~\eqref{eq:lasso}. Therefore, the loss function for training targets at regressing towards $\vx^\ast$. The authors of~\cite{lista} trained LISTA networks in a supervised setting, where they adopted a coordinate descent algorithm to solve the sparse coding problem~\eqref{eq:lasso} exactly to generate the optimal sparse code $\vx^\ast_j$ given each input $\vd_j$. In this way, a training set of $N$ input-label pairs $\sD_{\mathrm{train}}=\{(\vd_j, \vx^\ast_j)\}_{j=1}^N$ is generated. LISTA networks are then trained on this training set to minimize the mean squared error (MSE) between the approximated sparse codes and the ground-truth sparse codes
%
\begin{equation}
    \label{eq:lista-train-loss}
    \min_\vTheta \mathcal{L}(\vTheta) = \sum_{j=1}^{N} \| \vx^{(T)}_j - \vx^\ast_j \|_2^2,
\end{equation}
where $\vTheta = \{(\mW_1^{(i)}, \mW_2^{(i)}, \rho^{(i)})\}_{i=1}^T$ involves all the parameters to learn and $\vx^{(T)}_j$ is dependent on $\vTheta$ since it is computed using~\eqref{eq:ista-unrolled-truncated} with $\vd_j$ being the input. In practice, the network is trained with stochastic optimization algorithms such as SGD and Adam, where training data are fed in mini-batches randomly sampled from the whole training set to calculate the loss function values and the gradients with respect to the parameters through backpropagation.

\begin{figure}
    \centering
    \includegraphics[width=0.43\textwidth]{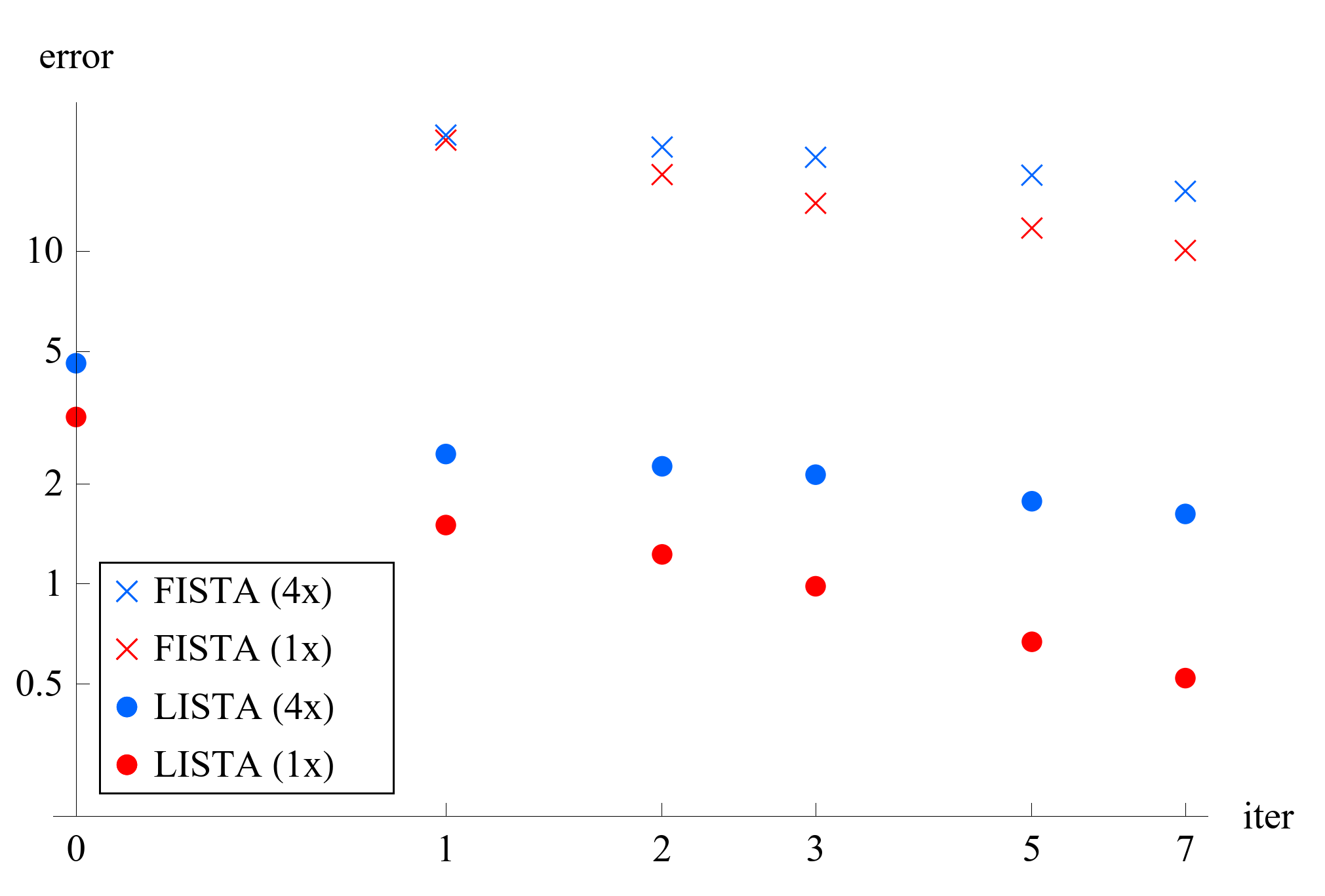}
    \caption{Comparison of LISTA an FISTA for solving sparse coding problems on the testing set. It takes 18 iterations of FISTA to reach the error for LISTA with just one iteration for $n=100$, and 35 iteration for $n=400$. \textit{Plot source: Figure~3 of~\cite{lista}}.}
    \label{fig:lista-figure3}
\end{figure}

\paragraph{Deploy and evaluate.}
After the offline training is finished, the final act is the deployment of the trained network for solving new optimization problems.
To evaluate the empirical performance of a trained LISTA network from the last step, we first prepare a testing set in the same way as we do for the training set --- all problems in the testing set share the same dictionary but the input vectors are different yet unseen during training.
Figure~\ref{fig:lista-figure3} (from \textit{Figure~3} in~\cite{lista}) shows the evaluation of LISTA on the Berkeley image database, where the input vectors are 10$\times$10 image patches and two dictionary matrices are learned: one dictionary is 1$\times$ complete, i.e., $m=n=100$, presented in red in the figure; the other is 4$\times$ complete, i.e., $n=4m=400$, presented in blue. LISTA is compared with FISTA~\cite{fista}, an extension to ISTA accelerated with momentum, in terms of how fast the prediction error decreases as iterations go on. It is reported in~\cite{lista} and we quote here that \textit{it takes 18 iterations of FISTA to reach the error for LISTA with just one iteration for $n=100$, and 35 iteration for $n=400$. Hence one can say that LISTA is roughly 20 times faster than FISTA for approximate solutions.}

\subsection{Design elements, techniques and applications}
\label{subsec:unrolling-design-factors}
In this subsection, we overview the various design elements, techniques, and applications that were developed to meet specific practical requirements. Researchers have created a rich body of literature by developing diverse design elements. We will explore the key design elements featured in the literature, discussing their origins, the impact of different design choices, and their advantages and disadvantages. To help improve the understanding of these design elements, we will establish connections with practical examples from the literature. 
We will also discuss the training techniques known to be effective for algorithm unrolling.

\paragraph{Share parameters across unrolled iterations?}
The first design element to consider when unrolling an iterative algorithm is: \textit{do we share the parameters across layers in the resulting network?}
Although it is trained and tested with a fixed number of iterations, the original LISTA~\cite{lista} unrolls ISTA into an RNN where all iterations share the same set of learnable parameters. In this case, we also say that the network has \textit{tied weights} in the literature. As the authors stated in~\cite{lista}, the recurrent system derived with unrolling, if truncated to a fixed number of iterations for all inputs, can also be viewed as a feedforward neural network of a fixed number of layers with tied weights. Based on this perspective, most work that followed~\cite{lista} proposed to unroll algorithms into feedforward networks by breaking the parameter sharing design, which means the learnable parameters in different layers are \textit{untied} and learned independently~\cite{borgerding2017amp,moreau2017understanding,chen2018theoretical}.
The untying process converts the recurrent formulation~\eqref{eq:ista-unrolled-truncated} into an MLP-type architecture:
\begin{align}
    \label{eq:lista-untied}
    \text{LISTA untied:} \quad
    & \vx^{(i)} = \eta\left( \mW_1^{(i)}\vx^{(i-1)} + \mW_2^{(i)}\vd; \rho^{(i)} \right), ~ i = 1, \dots, T, \\
    & \vx^{(0)} = \bm{0}. \nonumber
\end{align}
In this feedforward formulation, we denote the set of learnable parameters as $\vTheta=\{(\mW_1^{(i)}, \mW_2^{(i)}, \rho^{(i)})\}_{i=1}^T$.
The two types of networks are both valid choices for algorithm unrolling but they have their pros and cons depending on the scenarios and should be decided accordingly to make the most out of L2O.

It is worth noting that networks derived with algorithm unrolling actually process non-sequential inputs in most L2O applications despite them being RNNs. For example, in the case of unrolling gradient descent that we formulated in~\eqref{eq:pgd-lsq-recurrent} and the original LISTA~\cite{lista}, the input to the network is a single vector $\vd$, which is the input to the corresponding optimization problem. This means that the resulting RNNs only exploit the idea of parameter sharing and the recurrent structure, and there is no temporal correlation in inputs for them to capture. In some applications, however, the recurrent parameterization still has certain advantages over the non-recurrent counterpart.
To be more specific, an RNN and a feedforward network derived with algorithm unrolling differ in different aspects:

\begin{itemize}
    \item \textit{Capacity.} An RNN uses shared parameters across time and thus is compact in terms of the overall number of learnable parameters. In contrast, a feedforward network can have different parameters in each layer, independently learned from data. The total number of learnable parameters is thus proportionally increasing as the network goes deeper. This discrepancy in the compactness of parameterization results in differences in the \textit{capacity} of the resulting networks, by which we mean the ability of the network to approximate complex mappings when the learnable parameters take proper values~\cite{hornik1989multilayer}. With a larger of amount of parameters, feedforward networks have greater capacity to learn more complex mappings when compared with RNNs.
    On the other hand, however, the compact parameterization of an RNN brings benefits at the hardware level, e.g., lower storage and memory requirements for containing the trained network and for the computations. The hardware efficiency allows RNNs to be deployed on resource-constrained platforms such as edge devices~\cite{zhao2011online,han2013online} and distributed or Internet-of-Things (IoT) systems~\cite{zou2022proximal}.
    \item \textit{Flexibility.} An essential benefit of the idea of parameter sharing in RNNs is the \textit{flexibility} of applying them to inputs of varying lengths and generalize among them. While an RNN is truncated to a fixed number of iterations during training in~\cite{lista} and most algorithm unrolling works, it is straightforward to apply it for an arbitrary number of iterations for inference on new optimization problems. In fact, there is a natural need for such flexibility in L2O networks because some optimization problems are intrinsically more difficult than others and thus require more iterations to converge well. On the other hand, feedforward networks obviously lack such flexibility. Once the number of layers is decided during unrolling, the number of parameters is also fixed, and there is not a trivial way to properly apply parameters in previous layers to later layers. An simple extension that allows a trained feedforward network to process longer iterations is to append iterations after the last layer using a classic iterative algorithm~\cite{heaton2023safeguarded}.
    \item \textit{Trainability.} RNNs are generally considered harder to train compared to feedforward networks, especially when they are trained with long sequences or a large number of iterations in the case of algorithm unrolling~\cite{pascanu2013on}. This difficulty of training is caused by gradient vanishing (the more common case in practice) and explosion issues. Consider an RNN whose recurrent layer has an ill-conditioned Jacobian matrix. Gradients will be scaled downwards or upwards when backpropagated to earlier time steps with the Jacobian matrix multiplied to the gradients repeatedly. Practically, the gradient explosion issue can be mitigated with gradient clipping. However, alleviating the gradient vanishing issue for an RNN usually requires modifications to the architecture, which is impractical in algorithm unrolling because the architecture is decided by the original algorithm. In contrast, feedforward networks are generally easier to train but certain training tricks are also needed in some cases, especially when the resulting networks are deep and thus have too many learnable parameters. We will have a further discussion on training techniques later.
\end{itemize}

\paragraph{How to parameterize the resulting network?}
This is the core question that practitioners need to answer when adopting the algorithm unrolling methodology. The answer to this question decides the architecture of the resulting network, thus its capacity, trainability and empirical performance achievable in the end.
To find the better way of parameterization for the task of interest, we need to carefully consider of a number of factors. While it is impractical and overwhelming to list every piece of work in the literature and describe how its parameterization is done, we list three major factors that have driven different designs of parameterization: \emph{(i)} data structure, \emph{(ii)} task goal, and \emph{(iii)} interpretability.
Here, we limit the discussion to algorithm unrolling methods that strictly follow the formula of the original iterative algorithms. Variants with more advanced and more heavily parameterized deep learning extensions are elaborated in discussions of other design elements below.

\begin{itemize}[leftmargin=1.5em]
    \item[$-$] \textit{Data structure} describes the inherent complexity of the data distribution of the task that we strive to tackle using algorithm unrolling. It has direct influence on the decision of parameterization because the resulting network should at least provide the capacity to cope with that data complexity.
    The parameterization design of the original LISTA~\cite{lista} already reflects such consideration, where the authors proposed to use the same learnable matrices $\mW_1$ and $\mW_2$ for all sparse coding problems in the dataset. The rationale behind this parameterization strategy is rooted at the assumption that all sparse coding problems share the same dictionary $\mA$, which decides the value of $\mW_1$ and $\mW_2$ in classic ISTA iterations as formulated in~\eqref{eq:ista-unrolled-truncated}. Therefore, it is reasonable to learn the same $\mW_1$ and $\mW_2$ for all problems so that some \textit{common shortcut} can be utilized through data-driven learning from those optimization problems.

    However, if the data structure is more complex, e.g., the optimization problems use varying dictionaries, the original LISTA might fail to generalize. To deal with this increased data complexity, the authors of~\cite{aberdam2021adalista} proposed Ada-LISTA (standing for Adaptive LISTA) which \textit{augments} ISTA by replacing the input dictionary $\mA$ in each iteration with two matrices $\mG_1$ and $\mG_2$ that are linearly transformed from $\mA$ with two learnable matrices $\mW_1$ and $\mW_2$:
    \begin{align}
        \label{eq:ada-lista}
        \text{Ada-LISTA~\cite{aberdam2021adalista}:} \quad
        & \vx^{(i)} = \eta\left( (\mI_n - \alpha^{(i)} \mG_1^\top\mG_1)\vx^{(i-1)} + \alpha^{(i)}\mG_2^\top\vd; ~ \rho^{(i)} \right), ~ i = 1, \dots, T \\
        & \mG_1 = \mW_1\mA, \quad \mG_2 = \mW_2\mA, \quad \vx^{(0)} = \bm{0}. \nonumber
    \end{align}
    Here, the parameter matrices $\mW_1$ and $\mW_2$ are shared across layers but the layerwise learnable step sizes and thresholds are untied, similar to~\eqref{eq:lista-untied}. It is worth noting that in the formulation of Ada-LISTA the input to the network now is an observation-dictionary pair $(\vd, \mA)$.
    This augmentation endows Ada-LISTA enough capacity to learn from and generalize among sparse coding problems with varying dictionaries. But this capacity is still limited. The authors showed theoretically and empirically that Ada-LISTA can adapt to variants of a clean dictionary perturbed with column permutations or (small) random noises.
    If we expect the dataset contains more general sparse coding problems, further improvement to the network is necessary.

    \item[$-$] \textit{Task goal} is the target that we expect the resulting network to yield. In the original LISTA~\cite{lista} and many following works~\cite{moreau2017understanding,ablin2019learning,aberdam2021adalista}, 
    the task goal is to learn fast approximation to the minimizer of the LASSO objective $\vx^\ast$ as we denoted in~\eqref{eq:lasso}. In some applications, the goal is to recover an unknown signal, which we denote as $\bar{\vx}^\ast\in\sR^n$ to distinguish it from the minimizer of an objective function. The recovery is based on its observation generated with a forward model. A simple example of a forward model is a linear mapping perturbed with noise:
    \begin{align}
        \label{eq:linear-noisy-forward-model}
        \vd = \mA \bar{\vx}^\ast + \rvepsilon,
    \end{align}
    where $\mA$ represents the forward linear mapping and $\rvepsilon$ is the random noise. When $\bar{\vx}^\ast$ is sparse, the minimizer of LASSO $\vx^\ast$ is a decent approximation to $\bar{\vx}^\ast$ when the $\ell_1$ penalty level $\lambda$ is properly selected. However, $\vx^\ast$ is not a perfect recovery of $\bar{\vx}^\ast$ due to the bias introduced by the uniform $\ell_1$ penalty on all coordinates~\cite{fan2001variable}. Therefore, the two task goals are subtly different.

    \newcounter{tempCounter}
    The difference in task goal sometimes requires a change of the parameterization strategy. For example, instead of minimizing LASSO, the authors of~\cite{chen2018theoretical} unrolled ISTA into a feedforward network, aiming at solving the linear sparse recovery problems formulated in~\eqref{eq:linear-noisy-forward-model} with the sparsity assumption on $\bar{\vx}^\ast$. They theoretically showed that the $\mW_1^{(i)}$ and $\mW_2^{(i)}$ matrices in LISTA~\eqref{eq:lista-untied} need to asymptotically satisfy a coupling relationship, that is $\mW^{(i)}_1 - (I-\mW^{(i)}_2 A) \to \vzero $ when $i\to\infty$, so that the output of the network converges to $\bar{\vx}^\ast$ as the number of layers increases. With this result, the authors proposed to sustain this coupling structure throughout all layers, leading to the following parameterization:~\footnote{%
        We omit the usage of the advanced soft-thresholding technique called \textit{support selection} in~\cite{chen2018theoretical} for simplicity.%
        \setcounter{tempCounter}{\value{footnote}}%
    }%
    \begin{align}
        \label{eq:lista-cpss}
        \text{LISTA-CP~\cite{chen2018theoretical}:} \quad
        & \vx^{(i)} = \eta\left( \vx^{(i-1)} + {\mW^{(i)}}^\top(\vd - \mA \vx^{(i-1)}); ~ \rho^{(i)} \right), ~ i = 1, \dots, T, \\
        & \vx^{(0)} = \bm{0}. \nonumber
    \end{align}
    In contrast, it is shown in~\cite{ablin2019learning} that, however, if the resulting network is expected to converge to the LASSO minimizer, it must \textit{tend towards ISTA iterations} up to a learnable step size $\alpha^{(i)}$ in each layer. The threshold of the soft-thresholding function is set to $\rho^{(i)}=\alpha^{(i)}\lambda$ accordingly:
    \begin{align}
        \label{eq:step-lista}
        \text{Step-LISTA~\cite{ablin2019learning}:} \quad
        & \vx^{(i)} = \eta\left( \vx^{(i-1)} + \alpha^{(i)}\mA^\top(\vd - \mA \vx^{(i-1)}); ~ \alpha^{(i)}\lambda \right), ~ i = 1, \dots, T, \\
        & \vx^{(0)} = \bm{0}. \nonumber
    \end{align}
    We refer the reader to~\cite{chen2022learning} for more works in the literature targeted at the two different task goals.

    \item[$-$] \textit{Interpretability} is demanded in the field of algorithm unrolling for understanding how the acceleration is achieved via data-driven learning or for a more transparent computation process than black-box deep learning. The derivation of an interpretable network might bring up an architecture with special structures as a by-product and thus a new way of parameterization.
    
    For example, the authors of~\cite{moreau2017understanding} explained the acceleration of LISTA from the perspective of matrix factorization on the Gram matrix of the dictionary $\mA$. This perspective can derive a re-parameterized variant of the original LISTA which achieves similar convergence rate. In the case of linear sparse recovery, the authors of~\cite{alista} theoretically showed that all layers in the coupled LISTA formulated in~\eqref{eq:lista-cpss} can use the same matrix $\tilde{\mW}$ up to a learnable step size. Instead of being trained with data, $\tilde{\mW}$ is analytically generated by solving a normalized mutual coherence minimization with respect to the dictionary $\mA$. This analytical formulation leaves only the layerwise step sizes and thresholds for learning:~\footnotemark[\value{tempCounter}]
    \begin{align}
        \label{eq:alista}
        \text{ALISTA~\cite{alista}:} \quad
        & \vx^{(i)} = \eta\left( \vx^{(i-1)} + \alpha^{(i)}\tilde{\mW}^\top(\vd - \mA \vx^{(i-1)}); ~ \rho^{(i)} \right), ~ i = 1, \dots, T, \\
        & \vx^{(0)} = \bm{0}. \nonumber
    \end{align}
\end{itemize}

\paragraph{When to stop the iteration?}
While it is most standard to have an algorithm unrolled and truncated to a fixed number of iterations, it is usually needed to have a more sophisticated and adaptive strategy to decide how many layers or iterations are needed for a specific input. This need is well-motivated by the fact that optimization problems are not difficult to the same level. For problems that are easier to solve, it is possible to approximate the solution with only a few layers, while more iterations are needed for other more difficult problem instances. Similar motivation stirred the dynamic inference direction in the field of deep learning, where the researchers also want the network learn to \textit{exit early} if the input is easy~\cite{teerapittayanon2016branchynet}. This direction is also related to the optimal stopping time in control theory.

Another approach is to design a dynamic mechanism for algorithm unrolling to achieve this goal. For example, the authors of~\cite{zhang2018dynamically} designed a stopping policy based on reinforcement learning (RL) that dynamically decides the stopping time of a recurrent process for image restoration. The stopping policy is learned with Deep Q-Learning~\cite{mnih2015human}.
A similar RL-based approach is adopted in a plug-and-play method~\cite{wei2020tuning} for solving inverse problems in the image domain, where the authors learn a CNN-based policy network that simultaneously predicts the stopping probability and the parameters of the outer iterative framework. The network is trained using the actor-critic framework~\cite{sutton2000policy}. While it is designed for a plug-and-play method, we believe the dynamic stopping policy in~\cite{wei2020tuning} is still valuable for the practitioners in algorithm unrolling to develop their own dynamic strategy.

Besides RL-based methods, the authors of~\cite{chen2020learning} proposed a dynamic stopping strategy implemented with a \textit{variantional auto-encoder} (VAE) that learns when to terminate iterations early. The authors interpreted their proposed framework from a variantional Bayes perspective and connected it with RL. However, their method is still based on a neural network truncated to fixed \textit{max length} and only learns how to exit early, resulting in the lack of flexibility to extend to longer iterations.

\paragraph{Domain-specific designs.}
When algorithm unrolling comes to certain domains, we can do more than simply strictly following the original formulation of the iterative algorithm when unrolling and parameterizing it. Instead, we can customize parts of the resulting network to satisfy specific needs emerging in those domains. While the customization comes with the cost of losing the strict correspondence between the original algorithm and the neural network derived, which slightly undermines its interpretability, we can expect much better empirical performance.

One scenario where special designs are needed is in the field of wireless communication systems and scientific imaging where we need processing real-valued data and optimization problems. However, most modern deep learning frameworks do not support processing complex-valued input out of the box. To make algorithm unrolling still applicable in such applications, for example, the authors of~\cite{balatsoukas2019deep} used real-value decomposition to transform the original problem formulation so that a project gradient descent algorithm can be unrolled for solving MIMO detection problems in wireless communication systems.

Another important and ubiquitous example is the exploitation of convolutional neural networks in algorithm unrolling. CNNs have been proved to be extremely effective for image-based tasks thanks to the great properties of convolutional layers that we discussed see Section~\ref{subsec:dnn}. Therefore, when algorithm unrolling is applied to such tasks, it is well-motivated to combine it with CNNs. For example, in~\cite{zhang2018istanet}, the authors unrolled ISTA algorithm for natural image compressive sensing. The optimization-based compressive sensing minimizes a LASSO-type objective but assumes the image signal has a sparse representation under \textit{some linear transformation}, say $\Psi$, such as DCT and wavelet transformations. However, this assumption only approximately holds and barely captures the complexity of natural images. Therefore, the authors of~\cite{zhang2018istanet} proposed to insert two separate two-layer CNNs, $\mathcal{F}^{(i)}$ and $\Tilde{\mathcal{F}}^{(i)}$, into each layer of the network, with $\mathcal{F}^{(i)}$ replacing the $\Psi$ transformation and $\Tilde{\mathcal{F}}^{(i)}$ replacing the backward transformation, which is usually $\Psi^\top$ in classic algorithms.
For better regularization and interpretability, the authors also proposed to penalize $\|\Tilde{\mathcal{F}}^{(i)}\circ\mathcal{F}^{(i)}(\vx^{(i)})-\vx^{(i)}\|$ during training so that the learned $\mathcal{F}^{(i)}$ and $\Tilde{\mathcal{F}}^{(i)}$ are encouraged to follow the symmetry constraint.

This combination of algorithm unrolling and more advanced CNN modules has been empirically proved to be successful on the task of natural image compressive sensing~\cite{zhang2018istanet}, and been later applied to many other image-domain tasks, especially the field of biomedical imaging, e.g., CT imaging~\cite{kang2018deep,gupta2018cnn}, photo-acoustic tomography (PAT)~\cite{hauptmann2018model}, magnetic resonance imaging (MRI)~\cite{pramanik2020deep}, ultrasound imaging~\cite{solomon2019deep} and Shepp-Logan phantom and human phantoms~\cite{adler2018learned}. Another case of introducing CNNs to algorithm unrolling is where the optimization problems are themselves modeled with convolutional operators, such as convolutional sparse coding~\cite{sreter2018learned}, which is later extended in~\cite{alista}.

\paragraph{Instance adaptivity.}
Neural networks derived from algorithm unrolling typically learn parameters corresponding to components in the original algorithm, such as step sizes in gradient descent and thresholds in ISTA. Once the training is completed, the set of learned parameters is fixed. However, the optimal choices for these components are often instance-specific. To address this, some works propose adaptively generating these components for each instance using information from the current iteration or past optimization trajectory, a property known as \emph{instance adaptivity}. Two main methodologies exist for instance-specific generation.

One way is to adopt neural network-based augmentations. The authors of~\cite{behrens2020neurally} proposed an extension to ALISTA~\eqref{eq:alista} that generates step sizes and thresholds using an LSTM module. The LSTM network takes the $\ell_1$ norm of the residual of the current iterate and the $\ell_1$ norm of the residual transformed with the $\tilde{\mW}$ matrix analytically generated as in~\cite{alista}. The outputs of the LSTM network are two scalars, This method is based on the authors' observation on the correlation between the $\ell_1$ norm of the residual and the $\ell_1$ norm of the difference between the current iterate to the ground truth signal in each iteration. A recent work~\cite{liu2023towards} also utilizes an LSTM network for generating the diagonal matrix that preconditions the gradient step, the additive bias term, and the step size of the momentum term. The authors derived the basic mathematical conditions that successful update rules commonly satisfy. The overall framework derived in this way shows great adaptivity and even the ability to generalize cross datasets and tasks.

Another method for generating the instance-specific components for iterations is via analytical derivation which guarantees theoretical convergence. The authors of~\cite{chen2021hyperparameter} introduced a new ALISTA parameterization augmented with momentum acceleration for the sparse recovery task, and then analytically derived the formula of the instance-optimal key components of the ALISTA algorithm~\eqref{eq:alista}, i.e., the step size, soft-thresholding threshold, step size of the momentum term and size of the selected support. This formulation leaves reduces the number of learnable parameters in ALISTA to only tuning three scalars but endows the resulting network with great adaptivity.

\paragraph{Training techniques.}
%
Here we have a brief introduction to the selection of training objective and regularization techniques in algorithm unrolling training.

\begin{itemize}[leftmargin=1.5em]

    \item[$-$] \textit{Training objective} mainly depends on the task goal, which we have discussed about its relationship with the design of parameterization. Task goal also decides the training objective and the loss function utilized for training. We still take LISTA as an example, which can be trained to either \emph{(i)} approximate the LASSO minimizers; or \emph{(ii)} recover unknown signals.
    
    We generally have two ways to achieve the first goal. One approach is to generate optimal solutions of the optimization problems in the training set using a third-party solver, which for sure induces computation overheads during dataset generation and requires the precision of the solver. The neural network is then trained to regress the generated solutions using a regression loss function, e.g., the mean squared error (MSE).
    The second approach is to directly adopt the LASSO objective function as the loss function and to train the network to minimize the objective value at the output of the whole network. In this case, however, extra caution is required for potential numerical stability issues, especially when the objective includes logarithmic operations. In the case of LASSO, we have observed numerical issues when the $\ell_1$ penalty level $\lambda$ is large (around or larger than $10^{-1}$). Our conjecture is that the absolute value function consistently yield gradients with magnitude $1$, which is likely to be scaled upwards and to cause gradient explosion during backpropagation.
    
    For the second goal, it is straightforward to train the neural network in a supervised regression manner as long as the ground truth signals are accessible for training. However, this accessibility is forbidden in certain real-world applications. Two workarounds for this limitation are: (1) we can change the formulation of the optimization problems so that the optimal solutions approximate the desired signals better, and then fall back to the previous case; or (2) we can train the neural network on synthesized problems whose ground truth signals are generated instead of solved. The latter approach requires the neural network to have good generalization ability to generalize from synthesized problems to those in real-world scenarios.
    
    \item[$-$] \textit{Regularization techniques} are designed for more stabilized training process and better generalization ability of deep neural networks. Firstly, most of standard deep learning regularization techniques introduced in Section~\ref{subsec:dnn-training} can also be used for algorithm unrolling training, such as early stopping when the network performance stops improving on a hold-out validation set. Secondly, customized regularization terms can be added to the training loss function to encourage sustaining special structure in the neural network derived. A good example is the regularization term added for training ISTA-Net~\cite{zhang2018istanet} which enforces the composition of the forward and backward convolutional modules to approximate identity mapping, as we discussed earlier in the \textbf{Domain-specific designs} paragraph. Moreover, special training schemes can also serve as implicit regularization on the training process. For example, a stage-wise training strategy with learning rate decaying was adopted in~\cite{borgerding2017amp,chen2018theoretical} for training resulting networks with untied parameters, which was empirically shown to help with stabilizing the training and achieving better performance.
\end{itemize}

\paragraph{Applications.}
The application of algorithm unrolling extends across a diverse array of fields, where its implementation has been markedly successful. Particularly in the area of image restoration and reconstruction, a burgeoning literature has emerged, reflecting its widespread applicability. Within this field, algorithm unrolling has been utilized in classical low-level restoration tasks, such as denoising~\cite{zhang2017leanring}, deblurring~\cite{corbineau2019learned}, inpainting~\cite{aberdam2021adalista}, super-resolution~\cite{giryes2018tradeoffs} and video background subtraction~\cite{cai2021learned}. It also plays a crucial role in image artifacts removal tasks like JPEG artifact reduction~\cite{wang2016d3}. Furthermore, reconstruction tasks, e.g., natural image compressive sensing~\cite{zhang2018istanet,mardani2018neural}, have also benefited from the application of this method.
The scope of algorithm unrolling is not restricted to the aforementioned areas but also encompasses scientific imaging domains, including seismic, medical, and biological imaging. An exhaustive and well-integrated review of algorithm unrolling in the context of biomedical imaging can be found in~\cite{li2021deep}.
Additionally, this technique has found resonance in wireless communication systems, offering solutions to various challenges such as resource management~\cite{takabe2020complex}, channel estimation and signal detection~\cite{he2020model}, and LDPC coding~\cite{wadayama2019deep}. A seminal survey that elaborates on the methods of algorithm unrolling within communication systems can be referred to in~\cite{balatsoukas2019deep}.

\subsection{Mathematics behind algorithm unrolling}
\label{subsec:unrolling-math}

In algorithm unrolling techniques, iterative optimization algorithms, or their generalized variants, are treated as deep neural networks. This raises many interesting mathematical problems, some of which are still open or not well-studied. In this subsection, we present some of them. For more aspects on related topics, readers can refer to \cite{chen2022learning,monga2019algorithm,scarlett2022theoretical}.

\paragraph{Model selection.} Let us recall the recipe introduced above. The first step is to set up the optimization model and its associated iterative scheme. The optimization problem is naturally determined by the user's interest, while determining the specific formula of the iterative scheme is an art. Within the general framework presented in Section \ref{subsec:unrolling-basics}, the challenge lies in specifying the operator $\vg$ in \eqref{eq:model-unroll}. Clearly, not all operators are suitable. Some must be dismissed, such as:
\begin{itemize}
    \item $\vg(\vd,\vx;\vtheta)=\vx + \vone$. This operator does not yield a fixed point.
    \item $\vg(\vd,\vx;\vtheta)=2\vx$. Iteration through this operator cannot converge unless the starting point is $\vx=\vzero$.
\end{itemize}
A very recent study \cite{liu2023towards} marks an initial step in this direction. The authors study some basic conditions that such an operator $\vg$ should satisfy. Based on these conditions, they proposed a structured iterative scheme that exhibits superior performance on convex problems.

\paragraph{Expressivity.} Assuming the optimization model and associated iterative algorithm are determined, given an operator $\vg$, consider the following two sequences:
\begin{itemize}
    \item Sequence 1: $\{\vx^{(i)}\}_{i=1}^T$ generated by \eqref{eq:generic-algorithm}.
    \item Sequence 2: $\{\vx^{(i)}\}_{i=1}^T$ generated by \eqref{eq:model-unroll}.
\end{itemize}
This raises a question: are there parameters $\{\vtheta^{(i)}\}_{i=1}^T$ in \eqref{eq:model-unroll} such that Sequence 2 is significantly better than Sequence 1, provided that the parameters $\vtheta$ in \eqref{eq:generic-algorithm} is properly chosen? If not, algorithm unrolling may not be a wise option, as it induces extra cost of training without surpassing traditional algorithms. Fortunately, many works have demonstrated and quantified the potential advantages of algorithm unrolling over traditional methods. For example, in \cite{chen2018theoretical}, the authors show in the context of LISTA that, a conventional algorithm typically converges in a sub-linear rate, while there exist parameters such that algorithm unrolling yields a linear rate. For more results, readers may refer to \cite{alista,ablin2019learning,xie2019diff,wu2020sparse,zarka2020deep,yang2020learning,takabe2020theoretical,aberdam2021adalista,chen2021hyperparameter}.
Note that these results only show the existence of desirable parameters, but do not guarantee such parameters can be obtained by training. In machine learning, such results are categorized as the expressive power of machine learning models, regarding unrolled algorithms as a deep neural network.

\paragraph{Training.} Given the existence of parameters such that algorithm unrolling surpasses conventional methods, the question arises: Are there algorithms to compute such parameters with theoretical assurance? Specifically, one needs to solve minimization problems like \eqref{eq:unrolling-train-loss}, which is essentially a bi-level optimization. The lower-level optimization is defined by \eqref{eq:model-unroll} and the upper-level optimization is given by \eqref{eq:unrolling-train-loss}. Although complete results with the optimality guarantee still lack, some studies like \cite{brauer2022learning} identify some theoretical properties of the gradient of the loss function. Further references to related topics can be found in Section \ref{sec:diff-through-opt}.

\paragraph{Generalization.} Note that, even if problem \eqref{eq:unrolling-train-loss} can be perfectly solved, there might still be a gap between the trained and ideal models. The gap results from the differences between the training set $\sD_{\mathrm{train}}$ (defined before \eqref{eq:unrolling-train-loss}) and the upcoming new instances not in $\sD_{\mathrm{train}}$. To quantify this gap, one needs to create a testing set $\sD_{\mathrm{test}} = \{(\vd_j,\vx^\star_j)\}_j$, where all instances are independent of those in the training set. Then one can define the loss function on that testing set as:
\[ \ccL_{\mathrm{test}}(\vTheta) = \sum_{(\vd,\vx^\star) \in \sD_{\mathrm{test}}} \ell\Big( G(\vd;\vTheta), \vx^\star \Big) \]
If a relation like
\[ \ccL_{\mathrm{test}}(\vTheta_\star) \leq \ccL(\vTheta_\star) + \delta,\quad \text{where }\vTheta_\star := \argmin_{\vTheta} \ccL(\vTheta) \]
can be established with a reasonable $\delta > 0$, the trained model is said to have robust generalization performance on instances not found in the training set. Luckily, tools in statistical learning such as Rademacher complexity analysis can confirm such results.  Refer to \cite{chen2020understanding,behboodi2020generalization,schnoor2021generlization,joukovsky2021generalization,kouni2022generalization} for more details.

\section{Plug-and-Play Methods: DNN-Assisted Optimization Model}
\label{sec:pnp}

In this section, we present the \textit{Plug-and-Play} (\textbf{PnP}) framework with DNN assistance, which incorporates deep learning-based image denoisers into classic iterative algorithms.
More precisely, this is achieved by substituting an analytic expression with a trained neural network (usually an image denoiser), followed by the immediate deployment of the modified algorithm to optimize the subjects derived from the identical task distribution. Additional training of the modified algorithm is optional.

\subsection{An introduction to Plug-and-Play methods}

Here, we provide a demonstration of the PnP framework, integrating a neural network with a classical optimization algorithm, the alternating direction method of multipliers (ADMM). 
Consider the optimization problem: $\Min_{\vx \in \sR^{d}} f(\vx) + g(\vx)$. While a direct minimization of $f+g$ can be challenging, minimizing $f$ and $g$ separately might be considerably simpler. Hence, one can alternately tackle $f$ and $g$ in an iterative manner. This is facilitated by introducing an auxiliary variable $\vy$:
\begin{align}
    \label{eq:admm-objective}
    \Min_{\vx,\vy} f(\vy) + g(\vx), \quad \text{subject to} ~ \vx = \vy.
\end{align}
Subsequently, the ADMM is employed to solve \eqref{eq:admm-objective} in an iterative fashion:
\begin{subequations}
    \begin{align}
        \vx^{(i)} &= \mathrm{Prox}_{\beta g}  (\vy^{(i-1)} - \vu^{(i-1)}) := \argmin_{\vx} \left\{ \beta g(\vx) + \frac{1}{2} \left\| \vx - \left(\vy^{(i-1)} - \vu^{(i-1)}\right) \right\|^2  \right\} \label{eq:admm-1} \\
        \vy^{(i)} &= \mathrm{Prox}_{\alpha f} (\vx^{(i)}   + \vu^{(i-1)}) := \argmin_{\vy} \left\{ \alpha f(\vy) + \frac{1}{2} \left\|\vy - \left(\vx^{(i)}   + \vu^{(i-1)}\right) \right\|^2  \right\} \label{eq:admm-2} \\
        \vu^{(i)} &= \vu^{(i-1)} + \vx^{(i)} - \vy^{(i)} \label{eq:admm-3}
    \end{align}
\end{subequations}
In the above,~\eqref{eq:admm-1} and~\eqref{eq:admm-2} are the proximal operators for $f(\vy)$ and $g(\vx)$, where $\alpha$ and $\beta$ are step sizes. \eqref{eq:admm-3}~updates the dual variable $\vu$. Under mild assumptions on $f$ and $g$, the global convergence of ADMM can be proved. Specifically, from any starting point $\vx^{(0)},\vy^{(0)},\vu^{(0)}$, the sequences $\vx^{(i)}$ and $\vy^{(i)}$ will converge to the optimal solutions of \eqref{eq:admm-objective}. One can refer to~\cite{boyd2011distributed} for more details of ADMM.

Now let us revisit the optimization problem \eqref{eq:admm-objective} and ADMM in the context of image processing. For instance, imagine a scenario where the goal is to recover an image, denoted as $\vx^\star$, from its observations $\vb$. The connection between $\vb$ and $\vx^\star$ can be expressed as:
\[ \vb = \mA \vx^\star + \varepsilon, \]
where $\mA$ models the observation physics, such as CT, and $\varepsilon$ represents the noise during observation. For simplicity, we assume $\mA$ is a known linear operator. To recover $\vx^\star$ from $\vb$, one can solve an optimization problem 
\begin{equation}
    \label{eq:tv-quadratic}
\end{equation}
\[ \Min_{\vx \in \sR^d} \underbrace{(1/2)\|\mA \vx - \vb\|^2}_{f(\vx)} + \underbrace{\lambda \mathrm{TV}(\vx)}_{g(\vx)} \]
The function $f(\vx)$ measures the consistency of its input $\vx$ with the data we observe. The function $g(\vx)$ serves as a regularization term that measures the noise level in $\vx$ and the ``TV" within $g(\vx)$ stands for total-variation\footnote{For a 1D signal represented as $\vx \in \sR^{d}$, the total variation is given by $\sum_{i=1}^{d-1} |\vx[i+1]-\vx[i]|$. In the case of a 2D image denoted by $\vx \in \sR^{d_1 \times d_2}$, the TV is formulated as $\sum_{i=1}^{d-1}\sum_{j=1}^{d-1} \sqrt{(\vx[i+1,j]-\vx[i,j])^2 + (\vx[i,j+1]-\vx[i,j])^2}$.} \cite{rudin1992nonlinear}. Therefore, the roles of the proximal operators can be summarized as
\begin{align}
    \begin{array}{rrcl}
      \mathrm{Prox}_{\beta g} :  & \quad \text{noisy image}           & \quad \mapsto & \quad \text{less noisy image} \\
      \mathrm{Prox}_{\alpha f} : & \quad \text{less consistent image} & \quad \mapsto & \quad \text{more consistent image with observation}
    \end{array}
\end{align}
Note that the function $f$ is determined by data, while the function $g$ is designed manually. 
The design of the regularization function $g$ reflects our understanding on the structural property we wish an ideal solution $\vx^\star$ to possess. Given that $g$ only provides a simple approximation to the complexity of real-world image structures, it raises the question: 
\emph{why not supersede $g$ or $\mathrm{Prox}_{\beta g}$ with a more effective denoiser} that might not originate from optimization? This leads us to the essence of the ``plug-and-play" strategy: \textbf{plug} a stand-alone denoising operator into \eqref{eq:admm-1}:
\begin{align}
    \label{eq:pnp-admm-1}
    \vx^{(i)} = \mathrm{Denoiser}_\sigma (\vy^{(i-1)} - \vu^{(i-1)}),
    \tag{\ref{eq:admm-1}$^\prime$}
\end{align}
\textbf{and play} as \eqref{eq:pnp-admm-1},\eqref{eq:admm-2},\eqref{eq:admm-3}. In equation \eqref{eq:pnp-admm-1}, $\sigma$ plays a similar role of $\lambda$ in \eqref{eq:tv-quadratic}, denoting the estimated noise level in the input to the denoiser, which is accepted by most modern denoisers as an hyperparameter. The higher the noise level $\sigma$, the more aggressive denoising is performed. 

Specialized denoisers, such as BM3D~\cite{dabov2017image}, is more sophisticated than TV, and thus is more promising to yield visually better results than solving \eqref{eq:tv-quadratic}. Today, one can consider a denoiser based on a well-trained deep neural network, such as DnCNN \cite{zhang2017beyond}, within \eqref{eq:admm-1}, and witness remarkable performance in image-related tasks.

\subsection{Histories, motivations and a recipe for DNN-assisted Plug-and-Play}
\label{subsec:pnp-basics}

\paragraph{History.}
The PnP methodology dates back to the work~\cite{venkatakrishnan2013plug}, where the authors creatively replaced one of the two proximal operators in ADMM with a \emph{manually-designed denoiser} (such as non-local means~\cite{buades2005nonlocal}, K-SVD~\cite{aharon2006ksvd}, and BM3D~\cite{dabov2017image}). The resulting ADMM, detailed in equations \eqref{eq:pnp-admm-1},\eqref{eq:admm-2}, and \eqref{eq:admm-3}, exhibited greater performance than the original ADMM.

At beginning, PnP had not been explicitly linked to deep learning or deep neural networks. It was not until the emergence of~\cite{meinhardt2017learning,rick2017one,zhang2017leanring} that the notion of ``learning operators" was integrated into the PnP framework. Instead of using a manually designed denoiser, these studies replaced the proximal operator by a deep neural network that is trained with data. The empirical efficacy of this approach significantly surpassed previous implementations of PnP.

\paragraph{Motivations.} Focusing narrowly on image denoising, deep neural networks, especially convolutional networks, have been the state of the art. The concept of residual learning~\cite{he2016deep} and techniques such as batch normalization~\cite{ioffe2015batch} make network training significantly more effective, leading to trained networks with greater capacity and flexibility to process more real-world images.
\textit{DnCNN}~\cite{zhang2017beyond}
and more advanced variants~\cite{mao2016image,zhang2018ffdnet,zhang2021plug} integrate these techniques with convolutional neural networks and are vastly superior than analytical denoisers such as BM3D.

Therefore, it is well-motivated to use CNN denoisers in the PnP framework. On top of the great network capacity brought by stacking up more layers, the intuitions behind the success of CNN denoisers might also lie at the convolutional operation itself. Paper \cite{ulyanov2018deep} shows convolutional operations naturally generate \textit{self-similarity} in its output due to its inherent parameter-sharing structure, which makes CNNs very suitable for image restoration tasks, even if they have random parameters. See more discussion at the end of the next subsection.

\paragraph{A recipe for DNN-assisted plug-and-play methods.}
Procedures for deriving a DNN-assisted PnP method are more straightforward than other L2O frameworks as the derivation involves a minimal modification to the original optimization scheme:
\begin{itemize}[leftmargin=1.8cm]
    \item[\textbf{Step 1:}] Based on an optimization formulation, select an iterative solution method.
    \item[\textbf{Step 2:}] Identify one or more parts of the method that can be replaced with DNN denoisers. 
    \item[\textbf{Step 3:}] Train the DNN denoisers on proper datasets.
    \item[\textbf{Step 4:}] \textit{Plug} the well-trained DNN operators into the original optimization algorithm, and \textit{play}!
\end{itemize}

\subsection{Design elements and techniques of Plug-and-Play algorithms}
\label{subsec:pnp-design-factors}

\paragraph{Options and selection of plugged-in denoisers.}
Among the initial attempts of CNN-based denoising, \textit{DnCNN}~\cite{zhang2017beyond} stands as a significant development, which is grounded in VGG networks~\cite{vgg}. DnCNN set the new state of the art for image denoising through the incorporation of residual learning concepts~\cite{he2016deep} and the employment of batch normalization layers~\cite{he2016deep}. DnCNN is used for PnP in~\cite{ryu2019plug}. Another representative CNN denoiser was proposed in~\cite{mao2016image}, which devised an encoder-decoder architectural complemented by skip connections.

It is noteworthy, however, that these denoisers are configured to address inputs with a singular, fixed noise level. In scenarios with the presence of multiple noise levels within images, the approach necessitates the training of distinct CNNs, each tailored to the specific noise levels encountered. Consequently, the application of these models becomes contingent on the precise noise characteristics of the inputs, thus underscoring a limitation in their generalizability across diverse noise conditions.

In an effort to overcome the aforementioned limitation pertaining to flexibility and adaptivity, a novel CNN denoiser, termed FFDNet, was introduced by the authors of~\cite{zhang2018ffdnet}. It enables the handling of inputs contaminated with noises across a spectrum of levels through the utilization of a single network. The unique implementation accepts a noise level map as an auxiliary input, thereby fostering adaptivity to varying noise intensities.

Such adaptability is vital for PnP. As the iterative process of PnP advances, the quality of the denoiser's input also improves (i.e., less noisy). Consequently, an optimal strategy entails a gradual decrease of the noise level corresponding with these improvements. This is a capability uniquely afforded by FFDNet~\cite{yuan2020plug}, making it a strong candidate in contexts demanding adaptive responsiveness to varying noise conditions.

In a more recent PnP work~\cite{zhang2021plug}, the authors proposed a new CNN-based denoiser that essentially combines FFDNet and the U-Net architecture~\cite{unet}. U-Net is effective in image-to-image translation, because The skip connections introduced in U-Net provide the CNN-based denoisers the ability to synthesize features at different scales, thus generating more visually vivid features for better image restoration performances.

\paragraph{Parameter tuning in PnP methods.}
Using learning in Plug-and-Play (PnP) methods is not just about finding the right denoising tool or helper. It can also help figure out the best settings for the PnP process itself.
Although denoisers (including the DNN-based), are integrated into PnP frameworks in a hot-plugged fashion, each iterative step necessitates the determination of particular parameters, such as the step size of the other proximal step and the noise level specific to the current iteration. The calibration of these elements demands coordination in a sense for optimal performance of PnP. Intuitively, the step size of the other proximal step influences the ``noisiness'' of the input to the denoiser.

Recognizing the complexity of this interaction, the authors of~\cite{wei2020tuning} proposes to learn a policy network using reinforcement learning, to meticulously tune these parameters and even decide the appropriate stopping time of the PnP iterative process. The new optimization method can get results as good as those from using perfect settings, which are usually based on ground truth we cannot normally access. However, the introduction of reinforcement induces an extra step of learning a policy network besides pre-training the plugged-in denoiser.

\paragraph{Discussion on denoisers as image prior.}
A key motivation of adopting CNN-based denoisers in the PnP frameworks besides the network capacity is that convolutional neural networks serve as better image prior. But why do CNNs capture the features of natural images so well? In Section 2.2 we have discussed the three key characteristics of a CNN that are fundamental to its empirical success: local interaction with inputs, parameter sharing and the equivariance to translations. Here we provide another perspective that is slightly out of the context of PnP but related.
The authors of~\cite{ulyanov2018deep} showed that even a randomly-initialized CNN can be used as a handcrafted prior with excellent results in standard inverse problems including denoising. The network is trained to fit a degraded measurement by mapping a randomly generated but fixed 3D tensor without any other explicit regularization. They called this phenomenon deep image prior. In~\cite{ulyanov2018deep}, the authors attribute it to the convolutional filters themselves, which naturally impose self-similarity in the outputs. When combined with skip connections such as in U-Net, such self-similarity structures can be sustained at multiple scales in the image, which is essential for the visual quality of the restored images.

\subsection{Mathematical analysis on the convergence of Plug-and-Play methods}

Concerning the convergence of PnP-ADMM, theoretical guarantees were initially attempted in~\cite{sreehari2016plug}, based on the unrealistic assumption that the denoiser's derivative was doubly stochastic. Next, \cite{chan2016plug} proves the convergence of PnP-ADMM through the implementation of a more pragmatic assumption pertaining to ``bounded denoisers.''

Enforceable convergence conditions for PnP-FBS and PnP-ADMM were given in~\cite{ryu2019plug}, using the fixed-point analysis. Their assumption is a particular Lipschitz condition imposed on the denoisers, which can be achieved by applying normalization when the denoisers are trained. A more recent study is~\cite{matthieu2021enhanced}.

The aforementioned fixed-point framework was extended in~\cite{cohen2021regularization} to the RED-PRO framework.
Independently, \cite{heaton2022wasserstein} introduces a technique of learning projection operators onto the compact manifolds of true data. Under specific assumptions regarding the manifold and the adequate representation by the sampled data, they establish that their method for the synthesis of a PnP operator can, with statistical assurance, approximate the projection onto a low-dimensional manifold of data.

\paragraph{Enforcing special structures in denoisers.}
The integration of PnP with L2O has been revealed to confer theoretical advantages. Most of the existing theoretical insights into the convergence of PnP towards the desired solution depend on specific conditions associated with the denoisers. Despite the intrinsic difficulty in furnishing such guarantees for manually crafted denoisers, such as K-SVD and BM3D, the adoption of learning-based methods affords a flexibility that facilitates the convergence conditions.

Normalization or regularization strategies on the Lipschitz continuity and its variants of DNNs are useful tools such purpose. Notably, spectral constraints have been applied to both DNNs~\cite{bartlett2017spectrally,bansal2018can,oberman2018lipschitz} and Generative Adversarial Networks (GANs)~\cite{miyato2018spectral,brock2019large}. Regularizing Lipschitz continuity serves to stabilize the training process and increase the network's robustness to adversarial attacks~\cite{weng2018evaluating,qian2018l2}. 

In a specific illustration, the researchers in~\cite{miyato2018spectral} innovated a strategy to normalize all weights to be congruent with unit spectral norms, thereby circumscribing the Lipschitz constant of the entire network to a maximum value of one. Concurrently, the authors of~\cite{ryu2019plug} introduced a normalization technique explicitly tailored for the training of deep learning-based denoisers, ensuring compliance with the Lipschitz condition, further substantiating the theoretical strengths of the amalgamation of PnP with L2O.

\section{End-to-end learning and optimization as a layer}
\label{sec:opt-as-layer}
``Predict then optimize'' is a process used in decision-making processes, especially when predictions about future data are necessary for defining how to optimize and their constraints. The technique starts with training a machine learning model, then applies to model to predict future data, and finally optimizes the decision based on these predictions.
Typical examples include predicting travel times for route optimization, predicting sales for inventory optimization, and predicting returns and volatility for portfolio optimization.

Suppose we have a set of decision variables $\vx$ and a random variable $\vY$, which represents some future events so is currently unknown. The relationship between the decision $\vx$ and the random variable $\vY$ is represented by a loss function $f(\vx, \vy)$, which measures the performance of decision $\vx$ given the realization $\vy$ of $\vY$. In the predict-then-optimize technique, the first step involves building a machine learning model $\vs_\vtheta$ to predict $\vY$ based on some features $\vbeta$. This yields a prediction $\hat{\vy} = \vs_\vtheta(\vbeta)$. The next step is to solve the optimization problem according to the model prediction $\hat{\vy}$: $\min_{\vx}f(\vx,\hat{\vy})$. 
In a conventional predict-then-optimize pipeline, the learning and optimization are separated. In the training phase, we collect data pairs in the form of $\{\vbeta_i,\vy_i^\star\}_{i=1}^n$, and train a machine learning model $\vs_\vtheta$ that maps each $\vbeta_i$ to $\vy_i^\star$. Here $\vy^\star_i$ is an ``ideal  prediction" under the input $\vbeta$, which should be built upon domain knowledge. In the prediction phase, we make decisions $\vx$ according to the model prediction under the upcoming input $\vbeta$. Specifically, it is:
\begin{equation}
    \label{eq:predict-then-opt}
    \mathrm{Learning:}~ \Min_{\vtheta} \sum_{i=1}^n \ell(\vs_\vtheta(\vbeta_i), \vy_i^\star);
    \quad 
    \text{Predict-then-Optimize:}~ \Min_\vx f(\vx, \vs_\vtheta(\vbeta)).
\end{equation}
Here $\ell$ measures the difference between $\vs_\vtheta(\vbeta_i)$ and $\vy^\star_i$, which might be taken as the squared error, cross-entropy, etc. 

Note that the quality of the decision $\vx$ obtained from this optimization problem is typically measured by its actual loss $f(\vx, \vy)$, rather than the predicted loss $f(\vx, \hat{\vy})$. This introduces a discrepancy between the process of making the decision, which is based on the predicted loss, and the evaluation of the decision, which is based on the actual loss. This discrepancy is a key challenge in the predict-then-optimize technique. To mitigate this discrepancy, different approaches have been proposed and explored:
\begin{itemize}
    \item \textbf{Stochastic optimization.} This approach uses the samples and other prior knowledge of $\vY$, given the contextual information $\vbeta$, to model its randomness by a distribution $D(\vbeta)$. Then, the decision is made by minimizing the expected loss:
    \begin{align*}
        \Min_\vx \EE_{\vy\sim D(\vbeta)}\left[f(\vx,\vy)\right].
    \end{align*}
    \item \textbf{(Distributional) robust optimization.} This approach takes into consideration the uncertainty in the predictions. It finds the decision that performs best under the worst-case scenario within a certain uncertainty set around $\hat{\vy}$, say $B(\hat{\vy})$, leading to the following optimization problem:
    \begin{align*}
        \Min_\vx \max_{\vy\in B(\hat{\vy})} f(\vx,\vy).
    \end{align*}
    Distributional robust optimization combines stochastic optimization and robust optimization by introducing an ambiguity set $\gD$:
    \begin{align*}
        \Min_\vx \max_{D\in \gD(\vbeta)}\EE_{\vy\sim D}\left[f(\vx,\vy)\right].
    \end{align*}
\end{itemize}

\subsection{End-to-end learning}

More recently, the approach known as {end-to-end learning} or {decision-focus learning}~\cite{Kao2009directed,Donti2017taskbased,Wilder2019Melding,Bertsimas2020predictive,Elmachtoub2022smart} integrates the optimization step into the learning process, which serves as a straightforward approach to address the discrepancy induced by \eqref{eq:predict-then-opt}. It does not rely on an intermediate prediction step in the traditional approaches but instead attempts to learn a direct mapping from features $\vbeta$ to decisions $\vx^\star$. Specifically, in the training phase, we collect data pairs in the form of $\{\vbeta_i,\vx_i^\star\}_{i=1}^n$, and search parameters such that the final decision $\vx_i$ matches the ideal solution $\vx_i^\star$:
\[ \Min_{\vtheta} \sum_{i=1}^n \Bar{\ell}(\vx_i, \vx_i^\star),~~ \text{where}~~ \vx_i = \argmin_{\vx} f(\vx,\vs_\vtheta(\vbeta_i)). \]
Here, $\Bar{\ell}$ directly measures the difference between the model-based decision $\vx_i$ and the ideal decision $\vx^\star_i$. Once such a training process is done, the predict-then-optimize format remains the same with \eqref{eq:predict-then-opt}.

\paragraph{A general form.} In the context of end-to-end learning, we no longer depend on the intermediate variable $\vy$. Therefore, in this subsection, we can rewrite the function $f(\vx,\vs_\vtheta(\vbeta))$ into a more general form $f_\vtheta(\vx,\vbeta)$ and introduce some explicit constraints that the decisions $\vx$ must satisfy for reasons such as safety, regulations, and physical feasibility:
\begin{equation}
    \label{eq:amortized-opt}
    \Min_{\vx} f_\vtheta(\vx,\vbeta) ~~ \St ~ \vg_\vtheta(\vx,\vbeta) = \vzero, ~  \vh_\vtheta(\vx,\vbeta) \leq \vzero.
\end{equation}
Finding the optimal parameter $\vtheta$ in \eqref{eq:amortized-opt} can be expressed explicitly in the bi-level optimization problem:
\begin{equation}
    \label{eq:bi-level-opt}
    \Min_{\vtheta}~ \frac{1}{n}\sum_{i=1}^n \Bar{\ell}(\vx_i,\vx^\star_i), 
    \St~ \vx_i = \argmin_{\vx}\{f_\theta(\vx,\vbeta_i): g_\vtheta(\vx,\vbeta_i) = \vzero, h_\vtheta(\vx,\vbeta_i)\le \vzero\}.
\end{equation}
For simplicity, we ignore the case when the lower-level problem may have multiple solutions. This approach yields the mapping
\begin{align*}
   \vx_\vtheta:~\vbeta \mapsto \argmin_{\vx}\{f_\vtheta(\vx,\vbeta): \vg_\vtheta(\vx,\vbeta) = \vzero, \vh_\vtheta(\vx,\vbeta) \leq \vzero\}.
\end{align*}

\paragraph{Pros and cons.} In the traditional approaches \eqref{eq:predict-then-opt}, errors from the standalone prediction step $\delta \vy := \vy^\star - \hat{\vy}$ can propagate to the optimization solution, and the error can be difficult to control since, even in linear programming, arbitrary small errors in data can lead to an arbitrarily large error in the solution. End-to-end learning avoids this issue by optimizing the decision variables directly with respect to the actual loss. 
Additionally, end-to-end learning does not explicitly model the randomness in the data, but instead tries to directly learn a good decision mapping from data. So, it may perform better than stochastic optimization when the distribution of the data is complex and difficult to model

On the other hand, end-to-end learning does not explicitly consider the worst-case scenario, making it less suitable than robust optimization for cases where robustness is crucial. 

In addition, the solution mappings of most optimization are generally not differentiable, so they prevent the gradients easily flow back through in the backpropagation algorithm. We address this challenge below. 

\subsection{Differentiation through solution and KKT system}
\label{sec:diff-through-opt}
In order to learn a mapping $\vx_\vtheta$, a common approach is to approximate the gradient 
\begin{align*}
    \frac{d}{d\vtheta} \Bar{\ell}(\vx_\vtheta(\vbeta_i),\vx^\star_i)
\end{align*}
for every training data point $i$, and then employ an SGD-like algorithm, such as Adam \cite{adam}, to find a desirable $\vtheta$. 
The key obstacle is to ``differentiating through'' the solution, that is, to compute
\begin{align}\label{eq:sol_diff}
    \frac{\partial}{\partial \vtheta}\vx_\vtheta(\vbeta_i)
\end{align}
with $\vx_\vtheta(\vbeta_i)$ being the solution of an optimization problem. Even when the underlying optimization problem is a linear program, the differentiation is ill-defined. 
Specifically, a linear program may have multiple solutions, and even when the solution is unique, it is locally piece-wise constant with respect to the coefficients of the objective and those in the constraints. In this case, \eqref{eq:sol_diff} is either zero or undefined. 
To obtain an informative gradient with respect to $\vtheta$, modifications to the optimization formulation is necessary.

The Karush-Kuhn-Tucker (KKT) system is a set of equalities and inequalities that describe what optimal primal and dual solutions should satisfy. Under proper constraint qualifications, the KKT system is a necessary and sufficient optimality condition. Specifically, we introduce dual variables $\vlambda$ and $\vnu$ for constraints $\vg_\vtheta$ and $\vh_\vtheta$ respectively, and the KKT system is then defined for the primal-dual solution tuple $(\vx,\vlambda,\vnu)$ as follows:
\begin{equation}
    \label{eq:kkt}
    \cG_\vtheta(\vx,\vlambda,\vnu,\vbeta) = \vzero,~~\vnu \geq \vzero,~~ \vh_\vtheta(\vx,\vbeta) \leq \vzero,
\end{equation}
where the operator $\cG_\vtheta$ is defined by
\[ \cG_\vtheta(\vx,\vlambda,\vnu,\vbeta) = \begin{bmatrix}
    \frac{\partial f_{\theta}}{\partial \vx}(\vx,\vbeta) +  \left(\frac{\partial \vg_{\vtheta}}{\partial \vx}(\vx,\vbeta)\right)^\top \vlambda + \left(\frac{\partial \vh_{\vtheta}}{\partial \vx}(\vx,\vbeta)\right)^\top \vnu\\
    \vg_\vtheta(\vx,\vbeta)\\
    \vnu \circ \vh_\vtheta(\vx,\vbeta)
\end{bmatrix} \]
Since the above equation \eqref{eq:kkt} describes the optimality condition of $(\vx,\vlambda,\vnu)$, the decision $\vx_\vtheta(\vbeta)$ in \eqref{eq:sol_diff} must satisfy the KKT system above, serving as the $\vx$ in \eqref{eq:kkt}. This method enables one to determine the differential of $\vx_\vtheta(\vbeta)$ with respect to $\vtheta$. 
Furthermore, to illustrates the dependency of $\vlambda,\vnu$ on $\vbeta$ and $\vtheta$, we use the notions $\vlambda_\vtheta(\vbeta),\vnu_\vtheta(\vbeta)$. By assuming the smoothness of $f,\vg,\vh$ and applying the implicit function theorem, one can derive $\partial \vx_\vtheta / \partial \vtheta$ from \eqref{eq:kkt} as:
\begin{align}
    \frac{d}{d \vtheta} \cG_\vtheta \Big( \vx_\vtheta(\vbeta),\vlambda_\vtheta(\vbeta),\vnu_\vtheta(\vbeta),\vbeta \Big) = \vzero \label{eq:implicit-1}\\
    \implies \frac{\partial \cG_\vtheta}{\partial \vx} \frac{\partial \vx_\vtheta}{\partial \vtheta}
    + \frac{\partial \cG_\vtheta}{\partial \vlambda} \frac{\partial \vlambda_\vtheta}{\partial \vtheta}
    + \frac{\partial \cG_\vtheta}{\partial \vnu} \frac{\partial \vnu_\vtheta}{\partial \vtheta}
    + \frac{\partial \cG_\vtheta}{\partial \vtheta}  
    = \vzero \label{eq:implicit-2}\\
    \implies  \Big[
        \Big(\frac{\partial \vx_\vtheta}{\partial \vtheta}\Big)^\top~
        \Big(\frac{\partial \vlambda_\vtheta}{\partial \vtheta}\Big)^\top~
        \Big(\frac{\partial \vnu_\vtheta}{\partial \vtheta}\Big)^\top
    \Big]^\top  
    = - \Big[
        \frac{\partial \cG_\vtheta}{\partial \vx}~
        \frac{\partial \cG_\vtheta}{\partial \vlambda}~
        \frac{\partial \cG_\vtheta}{\partial \vnu}
    \Big]^{-1} \frac{\partial \cG_\vtheta}{\partial \vtheta} \label{eq:implicit-3}
\end{align}

When the system is not differentiable at a certain solution, the Jacobian matrix 
$J\cG_{\vtheta} := [\frac{\partial \cG_\vtheta}{\partial \vx}~
\frac{\partial \cG_\vtheta}{\partial \vlambda}~
\frac{\partial \cG_\vtheta}{\partial \vnu}]$ 
is not invertible and, thus, the above step \eqref{eq:implicit-3} fails to hold. In subsequent paragraphs, we will review different means to address this problem. 

Additionally, it is worth noting that solving the system defined in \eqref{eq:implicit-2} or calculating the inverse of the Jacobian matrix $J\cG_\vtheta$ is computationally expensive. Suppose $\vx \in \sR^{d_1}, \vlambda \in \sR^{d_2}, \vnu \in \sR^{d_3}$. The complexity of calculating $\partial \vx_\vtheta / \partial \vtheta$ via \eqref{eq:implicit-3} will be $\cO((d_1+d_2+d_3)^3)$. In the entire training process (i.e., solving \eqref{eq:bi-level-opt}), typically thousands of such gradients have to be evaluated, resulting in considerable computational demands. To address this issue, various techniques, such as Neumann approximation \cite{geng2021training} and Jacobian-free training \cite{samy2022jfb}, have been proposed and explored. 
An earlier work~\cite{amos2017optnet} applies the KKT approach to certain quadratic programs.

\paragraph{Heuristics.}
When the solution mapping is nondifferentiable and even discontinuous, one often applies heuristics \cite{griewank2008evaluating}. Take, for instance, a situation with a singular Jacobian matrix. Instead of directly solving the system as in \eqref{eq:implicit-3}, one can compute the least squares solution for system \eqref{eq:implicit-2}. Specifically, we solve the following quadratic optimization problem as an approximation to \eqref{eq:implicit-2}:
\[ \Min_{\vd}  \left\| \Big[
        \frac{\partial \cG_\vtheta}{\partial \vx}~
        \frac{\partial \cG_\vtheta}{\partial \vlambda}~
        \frac{\partial \cG_\vtheta}{\partial \vnu}
    \Big] \vd + \frac{\partial \cG_\vtheta}{\partial \vtheta} \right\|^2 + \lambda \|\vd\|^2, \]
where the last term $\lambda \|\vd\|^2$ serves as a regularizer to ensure the magnitude of $\vd$ does not become excessively large. The scalar $\lambda \in \sR_{+}$ controls the intensity of this regularization. 
Iterative methods, like the conjugate gradient approach, can be employed to solve this optimization problem. The number of iterations is determined by the desired accuracy specified by the user.

\paragraph{Adding regularizer.} 
To address the non-differentiable issue of the solution mapping, when the underlying optimization problem \eqref{eq:amortized-opt} is convex, a common approach from \cite{wilder2019melding} is to add a small quadratic function, $\frac{\alpha}{2}\|\vx\|_2^2$ and, therefore, yielding the new objective function $f_{\theta}(\vx,\beta,\alpha)$. This makes the problem \emph{strongly convex} and the solution not only unique but also differentiable. 
Another choice of regularizer is the logarithmic barrier function~\cite{mandi2020interior}, which is often used on conic optimization.

\paragraph{Perturbation.}
Gaussian smoothing or Gaussian perturbation is a standard technique used to make the differentiation of nonsmooth or discontinuous functions possible. This technique can be particularly useful in optimization problems where standard gradient-based methods are unable to operate due to the non-differentiability of the objective function. We leverage the methodology illustrated in \cite{berthet2020learning} as a prime example: perturb the parameter in the solution mapping $\vx_\vtheta$ and take expectation as follows:
\begin{align*}
   \vx_{\vtheta,\varepsilon}:~\vbeta \mapsto \EE_{\vn \sim \cN(\vzero, \vI)}\left[\vx_{\vtheta + \varepsilon \vn}(\vbeta)\right].
\end{align*}
Assuming the dimension of $\vtheta$ is $n$, the smoothed solution mapping can be explicitly expressed as:
\[
\begin{aligned}
    \vx_{\vtheta,\varepsilon}(\vbeta) = & \frac{1}{\sqrt{(2\pi \varepsilon^2)^n}} \int_{\vn \in \sR^n} \exp\left( -\frac{1}{2} \| \vn\|^2 \right) \vx_{\vtheta + \varepsilon \vn}(\vbeta) d\vn \\
    = & \frac{1}{\sqrt{(2\pi \varepsilon^2)^n} \varepsilon} \int_{\vtheta^\prime \in \sR^n} \exp\left( -\frac{1}{2\varepsilon^2} \| \vtheta^\prime - \vtheta\|^2 \right) \vx_{\vtheta^\prime}(\vbeta) d\vtheta^\prime 
\end{aligned}
\]
According to the standard real analysis, the mapping $\vx_{\vtheta,\varepsilon}$, as a convolution of $\vx_\vtheta$ with a smooth mollifier, must be smooth as long as $\vx_\vtheta$ is measurable. Its gradient can be derived by:
  \[
\begin{aligned}
 \frac{\partial}{\partial \vtheta} \vx_{\vtheta,\varepsilon}(\vbeta) = & \frac{1}{\sqrt{(2\pi \varepsilon^2)^n} \varepsilon} \int_{\vtheta^\prime \in \sR^n} \frac{\partial}{\partial \vtheta} \exp\left( -\frac{1}{2\varepsilon^2} \| \vtheta^\prime - \vtheta\|^2 \right) 
 \otimes \vx_{\vtheta^\prime}(\vbeta) d\vtheta^\prime \\
 = & \frac{1}{\sqrt{(2\pi \varepsilon^2)^n} \varepsilon} \int_{\vtheta^\prime \in \sR^n}  \exp\left( -\frac{1}{2\varepsilon^2} \| \vtheta^\prime - \vtheta\|^2 \right)  \frac{\partial}{\partial \vtheta} \left( -\frac{1}{2\varepsilon^2} \| \vtheta^\prime - \vtheta\|^2 \right) \otimes \vx_{\vtheta^\prime}(\vbeta) d\vtheta^\prime \\
 = & \frac{1}{\sqrt{(2\pi \varepsilon^2)^n} \varepsilon} \int_{\vtheta^\prime \in \sR^n}  \exp\left( -\frac{1}{2\varepsilon^2} \| \vtheta^\prime - \vtheta\|^2 \right) \left( \frac{1}{\varepsilon^2} (\vtheta^\prime - \vtheta) \right) \otimes \vx_{\vtheta^\prime}(\vbeta) d\vtheta^\prime\\
 = & \EE_{\vn \sim \cN(\vzero, \vI)} \left[ \frac{1}{\varepsilon} \vn \otimes \vx_{\vtheta + \varepsilon\vn}(\vbeta) \right].
\end{aligned}
\]
Here, $\otimes$ means the outer product. With proper assumptions on $f_\vtheta,\vg_\vtheta,\vh_\vtheta$, it can be shown that $\partial \vx_{\vtheta,\varepsilon} / \partial \vtheta \to \partial \vx_{\vtheta} / \partial \vtheta$ as $\varepsilon \to 0$. To calculate the smoothed gradient in practice, we can use the Monte-Carlo method: sample enough points $\{\vn_i\}_{i=1}^M$ that yields $\vn_i \sim \cN(\vzero,\vI)$ and estimate the gradient via
\[  \frac{\partial}{\partial \vtheta} \vx_{\vtheta,\varepsilon}(\vbeta) \approx \frac{1}{M} \sum_{i=1}^M \frac{1}{\varepsilon} \vn_i \otimes \vx_{\vtheta + \varepsilon\vn_i}(\vbeta). \]

\subsection{Fixed-point layer}
The KKT system is not the only optimality condition. Convex optimization problems with two or more components, as well as certain nonconvex optimization problems, also have another form of optimality condition of the fixed-point equation.
To quickly introduce it, we use the following problem as an example
\begin{align}\label{eq:box_tv_denoise}
    \Min_{\vx\in\RR^n} ~ & \frac{1}{2}\|\vA\vx - \vb\|^2 + \|\vD \vx\|_1 \\
    \St~ & \vl \le \vx\le \vu, \nonumber
\end{align}
where 
\begin{align*}
    \vD = \begin{bmatrix}
        -1 &  1 &        & & \\
           & -1 & 1      & & \\
           &    & \ddots & & \\
           &    &        &-1 & 1
    \end{bmatrix}.
\end{align*}
The function $\|\vD \vx\|_1$ is also known as one-dimensional total variation. 
This problem includes the following components:
\begin{itemize}
    \item the least squares term, $\frac{1}{2}\|\vA\vx - \vb\|^2$, which is  Lipschitz differentiable,
    \item the $\ell_1$-norm, $\|\cdot\|_1$, whose proximal operator is simple and closed-form, 
    \item the finite-difference operator, $\vD$, is a simple convolution-like matrix, and so is its adjoint $\vD^*$,
    \item the box constraint, $\vl\le \vx\le \vu$, is easy to project onto.
\end{itemize}
By applying the technique of operator splitting~Condat-Vu splitting~\cite{condat,vu}, $\vx^\star$ is a solution of problem \eqref{eq:box_tv_denoise} if and only if there exists $\vy^\star$ such as $\vz^\star = [\vx^\star; \vy^\star]$ satisfies
\begin{align*}
    \vz^\star = \cT \vz^\star,
\end{align*}
where $\cT$ is called a fixed-point operator constructed with the following individual operators: $\nabla_\vx \frac{1}{2}\|\vA\vx - \vb\|^2=\vA^\top(\vA\vx-\vb)$, $\vD$, $\vD^*$, $\prox_{\|\cdot\|_1}$, and $\proj_{[\vl,\vu]}$. With $\vz =[\vx^\top \vu^\top]^\top$, $\vz^{k+1}=\cT \vz^k$ is given by
\begin{align*}
    \vx^{k+1} & = \proj_{[\vl,\vu]}(\vx^k - \alpha \vD^* \vu^k - \alpha \vA^\top(\vA\vx^k - \vb))\\
    \vu^{k+1} & = \prox_{\beta\|\cdot\|_1^*}(\vu^k+ \beta \vA(2\vx^{k+1}-\vx^k)),
\end{align*}
where $\prox_{\beta\|\cdot\|_1^*}$ can be computed via the Moreau identity  $\prox_{\beta\|\cdot\|_1^*}(\vx)+\beta\prox_{\beta^{-1}\|\cdot\|_1}(\beta^{-1}\vx) = \vx$.
When the solution is not unique, every solution $\vx^\star$ must satisfy this fixed-point equation along with some $\vy^\star$.

In general, let $\vd$ denote all the data that define an optimization problem, and consider the operator $\cT(\cdot;\vd)$. We call $\vz^\star$ a fixed point of $\cT(\cdot;\vd)$ if $\vz^\star = \cT(\vz^\star;\vd)$. We say that the fixed-point problem
\begin{align*}
    \text{solve}_{\vz}\quad\vz = \cT(\vz;\vd)
\end{align*}
is equivalent to the optimization problem if every optimization solution $\vx^\star$ forms a fixed point $\vz^\star$ (possibly along with another subvector in $\vz^\star$) of $\cT(\cdot;\vd)$ and the converse also holds. Below we briefly review two common fixed-point operators. More fixed-point operators are covered in~\cite{ryuyin2022}. 

A \textbf{fixed-point layer}, introduced also as the deep equilibrium model~\cite{bai2019deep}, given by $\cT$ is the mapping:
\begin{align*}
    \vd \mapsto \vz^\star\quad{or}\quad \vd \mapsto \vx^\star,
\end{align*}
where $\vz^\star$ is a fixed point of $\cT(\cdot;\vd)$ and $\vx^\star$ is an optimization solution that can be obtained directly from $\vz^\star$.
In other words, $\vd$ is the input and $\vz^\star$ or $\vx^\star$ is the output. Computing the output can be done by either iterating the fixed-point operator or otherwise solving the original optimization problem.  

When the input $\vd$ depends on $\theta$, the fixed-point $\vz^\star$ is also a function of $\theta$. In addition, it is possible for $\cT$ to also include tunable parameters $\theta$; when this happens, we write $\cT_\theta$. 

\paragraph{Relationship with RNN.}
A fixed-point iteration, especially when the operator $\cT_\theta$ involves trainable parameters $\theta$, is reminiscent of RNN in that the iteration maintains a memory and gradually refines its estimate of the solution. However, fixed-point iteration is equipped with a rich, solid mathematical theory regarding its convergence behavior. As we will see below, it gives rise to cheaper approaches to obtain a gradient or a descent direction than RNN.

\paragraph{Example: Davis-Yin fixed-point operator~\cite{davis2017three}.} Consider
\begin{align*}
    \Min_{\vx\in\RR^n}\quad f(\vx)+ g(\vx) + h(\vx),
\end{align*}
where $f,g,h$ are proper, closed, convex functions and $h$ is differentiable.
Assume total duality~\cite{ryuyin2022}. Then $\vx^\star$ is a solution if and only if it satisfies
\begin{align*}
    0\in \partial f(\vx^\star) + \partial g(\vx^\star) + \nabla h(\vx^\star),
\end{align*}
where $\partial f,\partial g$ are subdifferentiables of $f$ and $g$, respectively. They are possibly set-valued.
Hence, solving for $\vx^\star$ is equivalent to the problem:
\begin{align*}
    \text{find}_{\vx\in\RR^n}\quad 0\in (\mathbb{A} + \mathbb{B} + \CC)\vx
\end{align*}
with $\mathbb{A}=\partial f$, $\mathbb{B}=\partial g$, and $\CC=\grad h$.
The Davis-Yin splitting refers to the equivalent problem: for $\alpha > 0$
\begin{align*}
     \text{solve}_{\vv\in\RR^n}\quad 
     \vz = \underbrace{\left(I - \mathbb{J}_{\alpha \mathbb{B}}+\mathbb{J}_{\alpha \mathbb{A}}(2\mathbb{J}_{\alpha \mathbb{B}}-I-\alpha\CC \mathbb{J}_{\alpha \mathbb{B}})\right)}_{=:\cT_{\text{DYS}}}\vz,
\end{align*}
where $\mathbb{J}_{\alpha \mathbb{B}}:=(I+\alpha \mathbb{B})$ is the resolvent of $\mathbb{B}$ and, with $\mathbb{B}=\partial g$, equals $\prox_{\alpha g}$. 
The last two problems are related through
\begin{align*}
    \vx^\star = \mathbb{J}_{\alpha \mathbb{B}} \vz^\star.
\end{align*}
When $\nabla h$ is $L_h$-Lipschitz and $\cT$ has at least one fixed point, the iteration 
\begin{align*}
    \vz^{k+1} = \cT_{\text{DYS}}\vz^{k}
\end{align*}
converges to a fixed point $\vz^\star= \cT_{\text{DYS}}\vz^\star$, and $\vx^k =\mathbb{J}_{\alpha \mathbb{B}} \vz^k$ converges to a solution $\vx^\star$.

When $g=0$, we have $\mathbb{J}_{\alpha \mathbb{B}}=I$ and $\cT_{\text{DYS}}$ degenerating to the forward-backward splitting. When $h=0$, we have $\CC=0$ and $\cT_{\text{DYS}}$ degenerating to the Douglas-Rachford splitting. Their corresponding optimization iterations are widely known as the proximal-gradient method and the ADMM.

\paragraph{Differentiation through a fixed-point layer.} To run a backpropagation algorithm through a fixed-point layer, there are three approaches: traditional, implicit, and Jacobian-free. 

In the traditional approach, the forward computation iterates  
$\vz^{k+1}(\theta) = \cT(\vz^k(\theta);\vd(\theta))$ for a finite number of times, $k = 0,1,\dots, K-1$. (The choice of $K$ is either fixed and adaptive to the need of approximation accuracy.) Then, backpropagation is based on the differentiation of that formula
\begin{align}\label{eq:fixed_back}
    \frac{\partial}{\partial\theta}\vz^{k+1}(\theta) 
    = J_\vz\cT(\vz^k(\theta);\vd(\theta)) \frac{\partial}{\partial\theta}\vz^{k}(\theta) + J_\vd\cT(\vz^k(\theta);\vd(\theta))\frac{\partial}{\partial\theta}\vd(\theta),
\end{align}
which establishes the relationships among successive $\frac{\partial}{\partial\theta}\vz^{k}(\theta)$, $k=0,\dots,K-1$, and $\frac{\partial}{\partial\theta}\vd(\theta)$.

The implicit differentiation approach is based on the fixed-point property:
\begin{align*}
    & \vz^\star(\theta) = \cT\left(\vz^\star(\theta),\vd(\theta)\right)\\
    \Rightarrow\quad &  \frac{\partial}{\partial\theta}\vz^\star(\theta) = J_\vz\cT(\vz^\star(\theta);\vd(\theta)) \frac{\partial}{\partial\theta}\vz^\star(\theta) + J_\vd\cT(\vz^\star(\theta);\vd(\theta))\frac{\partial}{\partial\theta}\vd(\theta) \\
    \Rightarrow\quad &  \frac{\partial}{\partial\theta}\vz^\star(\theta) = \left(I - J_\vz\cT(\vz^\star;\vd)\right)^{-1}J_\vd\cT(\vz^\star;\vd)\frac{\partial}{\partial\theta}\vd(\theta),
\end{align*}
which establishes the relationship between $\frac{\partial}{\partial\theta}\vz^{\star}(\theta)$ and $\frac{\partial}{\partial\theta}\vd(\theta)$. There is no restriction on the choice of the forward computation method. Suppose the forward method generates $\hat{\vz}$, an approximate of $\vz^{\star}$. Since $\vz^\star$ is unavailable, the backprapogation algorithm is left to use the approximate formula
\begin{align*}
    \frac{\partial}{\partial\theta}\hat{\vz}(\theta) = \left(I - J_\vz\cT(\hat{\vz};\vd)\right)^{-1}J_\vd\cT(\hat{\vz};\vd)\frac{\partial}{\partial\theta}\vd(\theta).
\end{align*}

The Jacobian-free approach~\cite{samy2022jfb} allows any forward computation method to generate $\hat{\vz}(\theta)$, and, starting from $\vz^0:=\hat{\vz}(\theta)$, it performs one fixed-point iteration $\vz^{1}(\theta) = \cT(\vz^0(\theta);\vd(\theta))$. It uses $\vz^{1}(\theta)$ as the output. To apply backpropagation, the Jacobian-free approach uses
\begin{align*}
    \frac{\partial}{\partial\theta}\vz^{1}(\theta) =  J_\vd\cT(\vz^0(\theta);\vd(\theta))\frac{\partial}{\partial\theta}\vd(\theta),
\end{align*}
which is free of computing the Jacobian, $J_\vz\cT(\vz^\star;\vd)$, and any matrix inversion.

A variant of the Jacobian-free approach~\cite{Bolte2023onestep} performs more than one fixed-point iteration. Specifically, from $\vz^0=\hat{\vz}(\theta)$, it performs $K\ge 2$ fixed-point iterations: $\vz^{k+1}(\theta) = \cT(\vz^k(\theta);\vd(\theta))$, $k=0,1,\dots,K-1$ to return $\vz^{K}(\theta)$. For backpropagation, use \eqref{eq:fixed_back}.

\paragraph{Comparisons of different approaches.}
The forward computations of all the approaches return an approximate fixed-point $\vz^\star$ or an approximate solution $\vx^\star$. The traditional approach, as its backpropagation uses \eqref{eq:fixed_back}, must record all the fixed-point iterates $\vz^0,\vz^1,\dots,\vz^{K-1}$. The other approaches do not need to record the fixed-point iterates. The variant of the Jacobian-free approach needs to record only the last few iterates, and the number .

For gradient backpropagation, all the approaches need to compute $J_\vd\cT$. The main differences lie in the Jacobian $J_{\vd}\cT$. In the traditional approach, the Jacobian is computed $K$ times. In the implicit differentiation approach, it is computed only once, but the inversion cost is generally cubic of the dimension $n$. The inversion can be approximated by the Neumann approximation \cite{geng2021training} at a lower cost, which grows with the number of terms used in the approximation. The Jacobian-free approach~\cite{samy2022jfb}, though completely avoids the Jacobian and inversion, does not produce a gradient. The authors~\cite{samy2022jfb,Bolte2023onestep,McKenzie2023faster}, under different assumptions such as contraction and linear independence constraint qualification, show that the obtained directions are descent direction despite not equaling gradients.

\section{AI-Guided Mixed-Integer Optimization}
\label{sec:mip}

In this section, we turn to Mixed-Integer Programming (MIP), an important class of optimization problems where some or all of the variables are constrained to be integers. MIP finds extensive use across diverse fields such as supply chain management, transportation, energy systems, finance, and more. Unfortunately,  MIP is generally considered to be NP-hard due to its integer constraints. It might be quite challenging to solve a fairly large MIP instance. To efficiently solve a specific type of MIP, one has to dive deeply into the theory and practice related to that MIP type. Only then is it possible to design a customized solver that surpasses general solvers on that particular problem class.

Given the multitudinous types of MIP, this process can be quite laborious. However, thanks to recent progress in machine learning, it has become possible to develop customized solvers built on a general solver and a dataset, without requiring too much expertise. This section will present some typical uses of machine learning for mixed-integer linear programming (MILP) and some associated practical and theoretical open questions.

\textbf{Remarks on notation.} In Section \ref{sec:mip}, we adopt the following conventions for notation. Boldface lowercase symbols, such as $\vx$, denote vectors. Typically, $x_j$ signifies the $j$-th element of vector $\vx$. It is important to distinguish between $\vx_j$ and $x_j$; the former refers to a vector indexed by $j$, the specifics of which are determined by context. On the other hand, bold uppercase symbols, like $\vA$, indicate matrices. By default, $A_{i,j}$ is the value at the matrix's $i$-th row and $j$-th column, while $A_{i,:}$ signifies the entire row.

\subsection{Preliminaries}
A general instance of an MILP problem is represented as:
\begin{equation}\label{eq:milp}
    \min_{\vx\in\RR^n} ~~ \vc^\top \vx,\quad \text{s.t.} ~~ \vx \in \sX_{\mathrm{MILP}} = \left\{\vx: \vA\vx\circ \vb,~~ \vl\leq \vx\leq \vu,~~ x_j\in \ZZ,\ \forall~j \in \sI \right\}.
\end{equation}
Here, $\vA\in\QQ^{m\times n}$, $\vc\in\QQ^n$, and $\vb\in\QQ^m$. The symbols $\vl$ and $\vu$ signify the lower and upper bounds of the variables, respectively, and are coordinate-wise selected from $l_j \in \QQ\cup\{-\infty\}, u_j \in \QQ\cup\{+\infty\}$. The dot $\circ$ represents the type of constraints, chosen element-wise from $\circ_j \in \{\leq, = ,\geq\}$. The index set $\sI\subset\{1,2,\dots,n\}$ includes those indices $j$ where the variable $x_j$ is constrained to be an integer. Each MILP instance can be described as a tuple $(\vc,\sX_{\mathrm{MILP}})$.

\paragraph{Basic Concepts.} The set $\sX_{\mathrm{MILP}}$ that describes all the constraints in \eqref{eq:milp} is termed the feasible set. An MILP instance is considered infeasible if $\sX_{\mathrm{MILP}} = \emptyset$ and feasible otherwise. For feasible MILPs, $\inf\{\vc^\top \vx : \vx\in \sX_{\mathrm{MILP}} \}$ is referred to as the optimal objective value. If there exists $\vx^\ast \in \sX_{\mathrm{MILP}}$ such that $\vc^\top \vx^\ast \leq \vc^\top \vx,\ \forall~\vx \in \sX_{\mathrm{MILP}}$, then $\vx^\ast$ is designated as an optimal solution. There can be a situation where the objective value is arbitrarily low, meaning that for any $R>0$, the inequality $\vc^\top \vx^\ast < -R$ is satisfied for some $\vx^\ast \in \sX_{\mathrm{MILP}}$. In such scenarios, the MILP is said to be unbounded, or its optimal objective value is $-\infty$.

\paragraph{Approaches and challenges via ML.} Addressing an MILP involves determining its feasibility and boundedness, as well as providing an optimal solution when possible. A direct approach to applying ML techniques on MILP, and indeed the one adopted in the early research literature, involves training a neural network to map an MILP to its results. For instance, suppose we consider the scenario where MILPs are feasible and bounded. In this context, we attempt to approximate the optimal solution of an MILP through a trained neural network:
\[\vx^\ast \approx \mathrm{NeuralNetwork}(\vc,\sX_{\mathrm{MILP}}).\]
Assuming the availability of sufficient data, this approach is theoretically doable. Given a collection of MILP instances and their corresponding solution $\{(\vc_k, \sX_{\mathrm{MILP}_{k}}, \vx^*_k)\}_{k=1}^K$, the neural network with parameters $\vtheta$ can then be trained using the following minimization problem:
\[\min_{\vtheta} \sum_k \|\vx^*_k - \mathrm{NeuralNetwork}(\vc_k,\sX_{\mathrm{MILP}_{k}}; \vtheta)\|\]
Here, the objective is to minimize the sum of the normed differences between the optimal solutions and the outputs of the neural network for each given MILP instance.
However, this approach presents two fundamental challenges:
\begin{itemize}
    \item Ensuring Feasibility: The outcomes of an NN are not inherently constrained by the stringent requirements of an MILP problem. Unlike in the case of continuous optimizations where constraint violations may be bounded by a small upper limit, there are no guarantees that the results from an NN will be feasible in the context of an MILP problem.
    \item Verifying Optimality: Even if the NN's output is an optimal solution, validating its optimality is nontrivial. Unlike in continuous optimizations where conditions like the KKT can validate optimality, there is no straightforward method to confirm the optimality of a solution generated by a neural network.
\end{itemize}
These issues highlight the limitations of a purely data-driven approach and point towards the necessity of incorporating algorithmic approaches. The conventional MILP solver is instrumental in overcoming these limitations by serving as a safeguard, ensuring the optimality of the solutions. Concurrently, the NN can be integrated into the solver to leverage its strengths. The intuition offered by the NN can facilitate adaptive and efficient problem-solving, providing initial estimates, or ``warm starts," to guide the solver. 

To summarize, a more effective practice of machine learning for MILP problems is utilizing neural networks to \emph{guide traditional solvers, rather than replacing these solvers entirely with neural networks.} In the remaining subsections, we will dive into various components of MILP solvers, outlining how machine learning methodologies can be incorporated to enhance their performance.

\subsection{Improving Branch and Bound by ML}
\label{sec:bnb-ml}
The Branch and Bound (BnB) algorithm is a widely employed method for solving general MILP problems and forms the foundation of modern MILP solvers. It operates on the principle of ``divide and conquer," and is guaranteed to provide an optimal solution when it exists. Throughout the BnB process, numerous decisions need to be made, such as variable selection and node selection. While these choices do not affect the final result, they substantially impact the algorithm's efficiency. Unfortunately, there are no theoretical guidelines on making the best decisions for general MILP. Therefore, these decision-making strategies often rely on practitioners' expertise and experience. Machine learning (ML) can be a helpful tool in developing these strategies.
Before introducing ML for branching, we will first present some key components in BnB: LP relaxation, branching, and bounding. 

\subsubsection{Branch and Bound}

\paragraph{LP relaxation.} Removing all the integer constraints in \eqref{eq:milp} (i.e., $x_j,j \in \sI$ is no longer required to be an integer) results in a linear program (LP), known as the LP relaxation of \eqref{eq:milp}. The resulting feasible set is denoted as $\sX_{\mathrm{LP}}$, and the LP relaxation is formulated as $\min_{\vx \in \sX_{\mathrm{LP}}}\vc^\top \vx$. The LP relaxation plays a pivotal role in solving the MILP. If the LP relaxation is infeasible, one can confidently assert the infeasibility of the original MILP. If the LP relaxation is lower unbounded, one can deduce the unboundedness of the MILP, given $\vA,\vb,\vc,\vl,\vu$ are all rational data. If the LP relaxation produces an optimal solution denoted as $\underline{\vx}$, this solution provides a lower bound for the MILP as:
\[\vc^\top\underline{\vx} \leq \vc^\top \vx, ~~\text{for all } \vx \in \sX_{\mathrm{MILP}}.\]
If $\underline{\vx}$ fulfills all the integer constraints, then it can be concluded that $\underline{\vx}$ is the optimal solution to the MILP. However, it might be quite challenging to identify such a $\underline{\vx}$ if the LP relaxation has multiple optimal solutions. More commonly, the LP relaxation may not present a solution that aligns with all the integer constraints. To measure this, we define the set of indices that $\underline{\vx}$ violates the integer constraints: \[\mathrm{frac}(\underline{\vx}) := \{j \in \sI: \underline{x}_{j}\not\in \ZZ\}\]

\paragraph{Branching.} Suppose $\mathrm{frac}(\underline{\vx}) \neq \emptyset$.  To obtain a solution that has less violation than $\underline{\vx}$, we select an index $j \in \mathrm{frac}(\underline{\vx})$ and generate two sub-MILPs:
\[\begin{aligned}
    \mathrm{MILP}_{\mathrm{up}}:&~~\min_{\vx}\vc^\top \vx ~~\mathrm{s.t.}~~ \vx \in \sX_{\mathrm{MILP}_{\mathrm{up}}} = \{\vx: \vx \in \sX_{\mathrm{MILP}},  x_j \geq \lceil\underline{x}_{j}\rceil \}\\
    \mathrm{MILP}_{\mathrm{down}}:&~~ \min_{\vx}\vc^\top \vx ~~\mathrm{s.t.}~~ \vx \in \sX_{\mathrm{MILP}_{\mathrm{down}}} = \{\vx: \vx \in \sX_{\mathrm{MILP}},  x_j \leq \lfloor\underline{x}_{j}\rfloor\}
\end{aligned} \]
Since $\sX_{\mathrm{MILP}} = \sX_{\mathrm{MILP}_{\mathrm{up}}} \cup \sX_{\mathrm{MILP}_{\mathrm{down}}}$, the feasibility of the original MILP can be determined by the feasibilities of the two sub-MILPs, and the optimal solution, if it exists, must be the best among the two optimal solutions of $\mathrm{MILP}_{\mathrm{up}}$ and $\mathrm{MILP}_{\mathrm{down}}$. Thus, solving the original problem can be reduced to solving two sub-MILPs. To address $\mathrm{MILP}_{\mathrm{up}}$ and $\mathrm{MILP}_{\mathrm{down}}$, we relax their integer constraints, yielding LP relaxations denoted as $\mathrm{LP}_{\mathrm{up}}$ and $\mathrm{LP}_{\mathrm{down}}$, respectively. The LP optimal solutions are denoted as $\underline{\vx}_{\mathrm{up}}$ and $\underline{\vx}_{\mathrm{down}}$, respectively. Compared to the original LP solution $\underline{\vx}$, $\underline{\vx}_{\mathrm{up}}$ and $\underline{\vx}_{\mathrm{down}}$ are closer to the MILP's optimal solution in two respects:
\begin{itemize}
    \item (Feasibility) Both $\underline{\vx}_{\mathrm{up}}$ and $\underline{\vx}_{\mathrm{down}}$ are more likely to yield integer solutions. Typically, in cases where integer variables are binary, $\mathrm{LP}_{\mathrm{up}}$ fixes $x_j$ as $1$ and $\mathrm{LP}_{\mathrm{down}}$ fixes $x_j$ as $0$, satisfying $x_j \in \ZZ$. Even in general cases, $\mathrm{LP}_{\mathrm{up}}$ and $\mathrm{LP}_{\mathrm{down}}$ help narrow down the search space by eliminating non-integer values in $(\lfloor\underline{x}_{j}\rfloor, \lceil\underline{x}_{j}\rceil)$.
    \item (Lower bounds) The solutions $\underline{\vx}_{\mathrm{up}}$ and $\underline{\vx}_{\mathrm{down}}$ provide tighter lower bounds than $\underline{\vx}$:
    \begin{equation}
        \label{eq:branch-tighter-lb}
        \vc^\top\underline{\vx} \leq \min\left( \vc^\top\underline{\vx}_{\mathrm{up}}, \vc^\top\underline{\vx}_{\mathrm{down}} \right) \leq \vc^\top \vx, ~~\text{for all } \vx \in \sX_{\mathrm{MILP}}.
    \end{equation}
    This is because the union of feasible sets of $\mathrm{LP}_{\mathrm{up}}$ and $\mathrm{LP}_{\mathrm{down}}$ forms a strict subset of the original LP relaxation:
    \[\sX_{\mathrm{LP}} \subsetneq
\big( \sX_{\mathrm{LP}} \cap \{\vx: x_j \geq \lceil\underline{x}_{j}\rceil\} \big) 
\cup 
\big( \sX_{\mathrm{LP}} \cap \{\vx: x_j \leq \lfloor\underline{x}_{j}\rfloor\} \big).\]
\end{itemize}
If $\mathrm{frac}(\underline{\vx}_{\mathrm{up}})=\emptyset$, we conclude that $\underline{\vx}_{\mathrm{up}}$ must be an optimal solution to $\mathrm{MILP}_{\mathrm{up}}$, and \emph{at least a feasible solution to the original MILP}. Its optimality will be determined after solving $\mathrm{MILP}_{\mathrm{down}}$. Otherwise, we conclude that $\underline{\vx}_{\mathrm{up}}$ is still not a feasible solution to the original MILP, and we repeat the above process and branch $\mathrm{MILP}_{\mathrm{up}}$ into two sub-problems. We handle $\mathrm{MILP}_{\mathrm{down}}$ with a similar approach.

\paragraph{Bounding.} If we view each sub-MILP as a node, these nodes together make up a binary tree known as \emph{the BnB tree}. In this context, the original MILP is regarded as the \emph{root node} of the entire BnB tree, and $\mathrm{MILP}_{\mathrm{up}}$ and $\mathrm{MILP}_{\mathrm{down}}$ are seen as the two \emph{children} of the root, and subsequent sub-problems of $\mathrm{MILP}_{\mathrm{up}}$ are considered children of $\mathrm{MILP}_{\mathrm{up}}$, and so on. Expanding this tree fully can be time-consuming. While in practice, some nodes might be pruned. To discuss the pruning method, we first need to understand two basic concepts:
\begin{itemize}
    \item (Upper bound) During the branching process, as soon as we find a feasible solution for the original MILP, we store it in memory. It remains there until a solution with a lower objective value is found and replaces it. This stored solution, denoted as $\overline{\vx}$, is known as the \emph{incumbent} solution and serves as an upper limit for the optimal objective value:
    \[ \vc^\top \vx^\ast \leq \vc^\top \overline{\vx}, \]
    where $\vx^\ast$ represents any optimal solutions of the original MILP. The upper bound of the optimal objective value, represented as $\overline{f}:=\vc^\top \overline{\vx}$, is known as the \emph{global upper bound} or the \emph{primal bound}.
    \item (Lower bound) As mentioned before, the LP relaxations of the MILP or sub-MILPs provide lower bounds of the objective value. After obtaining $\underline{\vx}$, we get a lower bound $\underline{f}:=\vc^\top\underline{\vx}$, denoted as the \emph{global lower bound} or the \emph{dual bound}. After $\mathrm{LP}_{\mathrm{up}}$ and $\mathrm{LP}_{\mathrm{down}}$ are solved, the dual bound will be updated with a higher value $\min\left( \vc^\top\underline{\vx}_{\mathrm{up}}, \vc^\top\underline{\vx}_{\mathrm{down}} \right)$. This process continues throughout the BnB procedure.
\end{itemize}
As we progress through this procedure, the primal bound decreases while the dual bound increases until they converge. Based on these two bounds, we can prune certain nodes in the BnB tree. If the optimal objective value of the LP relaxation of a sub-MILP exceeds the global upper bound, we can prune this sub-MILP (or node) as neither it nor its children can produce a solution better than the incumbent solution. It is important to note that this includes instances where the LP relaxation is infeasible since any children of an infeasible sub-MILP must also be infeasible. A full description of BnB is provided in Algorithm \ref{algo:bnb}.

\begin{algorithm}
\renewcommand{\algorithmicrequire}{\textbf{Input:}}
\renewcommand{\algorithmicensure}{\textbf{Output:}}
\caption{Branch and Bound Algorithm}
\label{algo:bnb}
\begin{algorithmic}[1]
    \REQUIRE An MILP instance $(\vc, \sX_{\mathrm{MILP}})$; Variable selection rule; Node selection rule
    \ENSURE Label ``Infeasible", Label ``Infeasible or Unbounded", or an optimal solution $\overline{\vx}$;
    \STATE Active node set $\sS \Leftarrow \{\}$. \COMMENT{The set $\sS$ consists of tuples structured as $(d,f,\sX)$. Each of these tuples symbolizes a node within the BnB tree. Here, $d$ stands for the depth of the node, $f$ signifies a \emph{local lower bound} of the node - equivalent to the node's parent's optimal LP objective value, while $\sX \subset \sX_{\mathrm{MILP}}$ characterizes the feasible set of the sub-MILP corresponding to the node.}
    \STATE The global lower bound $\underline{f} \Leftarrow -\infty$, the global upper bound $\overline{f} \Leftarrow +\infty$.
    \STATE Push a node with data $(0,\underline{f}, \sX_{\mathrm{MILP}})$ into $\sS$.
    \WHILE{$\sS$ is not empty}
        \STATE Pop a node $(d, f, \sX)$ from $\sS$ with the node selection rule.
        \IF{The parent's objective exceeds the upper bound: $f \geq \overline{f}$}
            \STATE \textbf{continue}
        \ENDIF
        \STATE Update the global lower bound: $\underline{f} \Leftarrow \min_{(d, f,\sX) \in \sS} f$.
            \IF{The lower and upper bound converge $\overline{f} = \underline{f}$}
                \RETURN $\overline{\vx}$.
            \ENDIF
        \STATE Solve $(\vc, \mathrm{Relax}(\sX))$ and denote its optimal solution as $\underline{\vx}$ if feasible. \COMMENT{Here $(\vc, \mathrm{Relax}(\sX))$ represents the LP relaxation of $(\vc, \sX)$.}
        \IF{$(\vc, \mathrm{Relax}(\sX))$ is infeasible or its objective exceeds the upper bound: $\vc^\top \underline{\vx} \geq \overline{f}$}
            \STATE \textbf{continue}
        \ELSIF{$(\vc, \mathrm{Relax}(\sX))$ is lower unbounded}
            \RETURN Infeasible or Unbounded.
        \ELSE
            \IF{$\underline{\vx}$ satisfies integer constraints: $\mathrm{frac}(\underline{\vx}) = \emptyset$}
                \IF{$\underline{\vx}$ has better objective $\vc^\top \underline{\vx} < \overline{f}$}
                    \STATE Update the upper bound and the incumbent: $\overline{f} \Leftarrow \vc^\top \underline{\vx}$, $\overline{\vx} \Leftarrow \underline{\vx}$.
                    \IF{The lower and upper bound converge $\overline{f} = \underline{f}$}
                        \RETURN $\overline{\vx}$.
                    \ENDIF
                \ENDIF
            \ELSE
                \STATE Pick an index $j \in \mathrm{frac}(\underline{\vx})$ with the variable selection rule.
                \STATE Create two branches: $\sX_{\mathrm{up}} = \sX \cap \{\vx: x_j \geq \lceil \underline{x}_j \rceil\}$, $\sX_{\mathrm{down}} = \sX \cap \{\vx: x_j \leq \lfloor \underline{x}_j \rfloor\}$.
                \STATE Add $(d+1, \vc^\top \underline{\vx}, \sX_{\mathrm{up}})$ and $(d+1, \vc^\top \underline{\vx}, \sX_{\mathrm{down}})$ into $\sS$.
            \ENDIF
        \ENDIF
    \ENDWHILE
    \RETURN Infeasible ~\textbf{if} $\overline{f} = +\infty$ ~\textbf{else} ~\textbf{return} $\overline{\vx}$
\end{algorithmic}
\end{algorithm}

\paragraph{Efficiency analysis.} The BnB approach is known to terminate in finite steps when applied to MILP problems with rational parameters, regardless of the sequence in which variables and nodes are chosen. However, the effectiveness of the BnB method significantly depends on the strategies used for variable selection and node selection. As demonstrated in Figure~\ref{fig:bnb-var-example}, certain examples show a case of how varying variable selection rules result in final BnB trees of different sizes. The same conclusion is observed in Figure \ref{fig:bnb-node-example} with various node selection rules. 

These examples prompt us to question: \emph{What are the optimal strategies for variable and node selection for a specific type of MILP instances?} Unfortunately, there are no definitive answers to this inquiry. It is still typical for practitioners to develop these strategies based on a \emph{trial and error} approach. By testing many different strategies on a specific type of MILP (e.g., scheduling, bin-packing, vehicle routing), practitioners can gain domain-specific knowledge that informs the development of new strategies. However, this iterative process is problematic due to:
\begin{itemize}
\item Its labor-intensive nature.
\item The absence of a systematic methodology. The acquired domain knowledge is specific to the type of MILP and does not easily translate to other domains.
\end{itemize}
Hence, it would be advantageous to develop a systematic approach to autonomously create or discover effective variable and node selection rules. The remarkable advancements in machine learning techniques in recent years lead us to ask: \emph{can we leverage machine learning to extract powerful strategies from data?} In the following subsections, we will show how to learn variable and node selection rules with machine learning.

\begin{figure}
\centering
\input{var-select-1.tikz}
\hfill
\input{var-select-2.tikz}
    \caption{Variable selection can influence the size of the BnB tree. Consider the MILP: $\min 11x_1 + 12x_2 +13x_3 + 14x_4$ subject to $x_1+x_2+x_3+x_4\geq2.5$, $x_2+2x_3 \geq 2.1$ and all the variables must be binary. The root node represents an optimal solution $\underline{\vx} = (0.95,1,0.55,0)$ to the LP relaxation and its corresponding objective $\vc^\top\underline{\vx} = 29.6$. In $\underline{\vx}$, two elements are fractional: $x_1$ and $x_3$. We have to decide which one to branch on. If we adopt the first approach—branching based on \emph{the fractional element with the smallest index}—we will choose $x_1$ to branch on, leading to the tree shown in the left figure. On the other hand, if we employ a different strategy—branching based on \emph{the element with the largest fractionality}—the choice is determined by $\argmax_{j} \min(\underline{x}_j - \lfloor \underline{x}_j \rfloor, \lceil \underline{x}_j \rceil - \underline{x}_j)$. In this case, $x_3$ is selected, resulting in the tree in the right figure. Continuing this process in either case will yield different BnB trees. Note that, in this case, the node selection strategy will not influence the BnB tree, and you can verify this with any node selection method.}
    \label{fig:bnb-var-example}

\vspace{5mm}

\input{node-select-1.tikz}
\hfill
\input{node-select-2.tikz}
    \caption{Node selection can influence the size of the BnB tree. Consider the MILP: $\min x_1 + 2x_2 +3x_3 + 4x_4$ subject to $x_1+x_2+x_3+x_4\geq2.5$, $x_2+2x_3 \geq 2.1$ and all the variables must be binary. We consistently branch based on the element with the highest fractionality but employ different strategies for node selection. In the figure on the left, we employ a \emph{depth-first} order for node selection, which results in the need to solve a total of $9$ LP problems. The order to access nodes are: \textcircled{0},\textcircled{1},\textcircled{2},\textcircled{3},\textcircled{4},\textcircled{5},\textcircled{6},\textcircled{7},\textcircled{8}. Conversely, the figure on the right demonstrates the use of a \emph{breadth-first} order for node selection, which leads to pruning $2$ LPs. The order to access nodes are: \textcircled{0},\textcircled{1},\textcircled{2},\textcircled{3},\textcircled{4},\textcircled{5},\textcircled{6}.}
    \label{fig:bnb-node-example}
\end{figure}

\subsubsection{Variable selection}
\label{sec:branching-rule}

 The variable selection rule, also known as the \emph{branching rule}, determines which variable to branch on when multiple fractional values are present in the current LP relaxation solution. To identify the ideal variable to branch on, we collect information related to all potential variables and make decisions based on the collected data, which may include:
\begin{itemize}
    \item The original MILP: $\vc,\sX_{\mathrm{MILP}}$
    \item The current MILP and its LP relaxation: $\sX$, $\mathrm{Relax(\sX)}$ and $\underline{\vx}$ (defined in Algorithm \ref{algo:bnb}, Line 13)
    \item Data generated during BnB: $\underline{f},\overline{f},\overline{\vx},\sS$.
\end{itemize}
All the aforementioned information is accessible during the BnB procedure. Now, let us present the formal formulation of a branching rule:
\begin{equation}
    \label{eq:branch-with-score}
    s_j = \cV(j,\cdot),\quad j^* = \argmax_{j \in \mathrm{frac}(\underline{\vx})}s_j.
\end{equation}
In this formulation, $\cdot$ represents the collected information, and $\cV$ is a mapping function that assigns a real number, which is named as a score, to a variable. Subsequently, the index with the top score is selected for branching. In case of a tie, lexicographical order is usually adopted. The mapping $\cV$ forms the cornerstone of this rule and is undoubtedly the main component we aim to design or learn. Some typical handcrafted branching rule is listed here:
\begin{itemize}
    \item (Most Infeasible Branching). We use the fractionality as the score: $s_j = \min(\underline{x}_j - \lfloor \underline{x}_j \rfloor, \lceil \underline{x}_j \rceil - \underline{x}_j)$.
    \item (Strong Branching). A primary aim of branching is to improve the lower bound. The Strong Branching (SB) rule seeks to achieve this by testing each variable $j \in \mathrm{frac}(\underline{\vx})$ to see which one significantly tightens the lower bound the most. Suppose there are $n_{\mathrm{frac}}$ elements in $\mathrm{frac}(\vx)$, a full strong branching constructs $2n_{\mathrm{frac}}$ sub-MILPs with
    \[ \sX_{j,\mathrm{up}} = \sX \cap \{\vx:x_j \geq \lceil \underline{x}_j \rceil\},\quad \sX_{j,\mathrm{down}} = \sX \cap \{\vx:x_j \leq \lfloor \underline{x}_j \rfloor\},\quad j \in \mathrm{frac}(\vx).\]
    Subsequently, we solve the LP relaxations of these MILPs $\{(\vc,\sX_{j,\mathrm{up}}),(\vc,\sX_{j,\mathrm{down}})\}_{j\in \mathrm{frac}(\vx)}$. The solutions of these are represented as $\{\underline{\vx}_{j,\mathrm{up}}, ~\underline{\vx}_{j,\mathrm{down}}\}_{j\in \mathrm{frac}(\underline{\vx})}$ respectively. The improvement in the objective of the two sub-MILPs, obtained by branching along $x_j$, is computed as $\delta_{j,\mathrm{up}} = \vc^\top \underline{\vx}_{j,\mathrm{up}} - \vc^\top \underline{\vx}$ and $\delta_{j,\mathrm{down}} = \vc^\top \underline{\vx}_{j,\mathrm{down}} - \vc^\top \underline{\vx}$. Given a small number $\varepsilon$ (say $\varepsilon=10^{-8}$), the SB score of $x_j$ is derived by:
    \begin{equation}
        \label{eq:sb}
        s_{j,\mathrm{SB}} = \cV_{\mathrm{SB}}(j,\vc,\sX,\underline{\vx},\varepsilon) = \max(\delta_{j,\mathrm{up}},\varepsilon) \times \max(\delta_{j,\mathrm{down}}, \varepsilon).
    \end{equation}
    \item (Pseudocost Branching). The Pseudocost (PC) branching rule is essentially a SB approximation that relies on historical data. It observes that a single variable is often branched multiple times. For example, in the left tree in Figure \ref{fig:bnb-var-example}, $x_2$ and $x_3$ are branched two times. Therefore, historical objective improvements can be stored and averaged for each variable, providing an estimate for predicting the SB score. Specifically, we maintain two quantities $\delta f_{j,\mathrm{up}}$ and $\delta f_{j,\mathrm{down}}$ for each variable $x_j$ and they are initialized by zero. Every time we branch along $x_j$, we calculate the ``objective improvement per unit" and store its historical average in $\delta f_{j,\mathrm{up}}$ and $\delta f_{j,\mathrm{down}}$ respectively:
    \[ \frac{\vc^\top \underline{\vx}_{j,\mathrm{up}} - \vc^\top \underline{\vx}}{\lceil \underline{x}_j\rceil - \underline{x}_j}, \quad \frac{\vc^\top \underline{\vx}_{j,\mathrm{down}} - \vc^\top \underline{\vx}}{\underline{x}_j - \lfloor \underline{x}_j\rfloor} \]
    With the two maintained quantities, we can calculate estimates for objective improvements: $\delta_{j,\mathrm{up}} = (\lceil \underline{x}_j\rceil - \underline{x}_j ) \delta f_{j,\mathrm{up}}$ and $\delta_{j,\mathrm{down}} = (\underline{x}_j - \lfloor \underline{x}_j\rfloor ) \delta f_{j,\mathrm{down}}$. Plugging these estimates into the right-hand side of \eqref{eq:sb}, we will obtain an approximated SB score, named the PC score $s_{j,\mathrm{PC}}$.
\end{itemize}

\textbf{(Imitating SB).} Simple branching rules, such as the first fractional rule and the most fractional rule presented in Figure \ref{fig:bnb-var-example}, are straightforward but may lack effectiveness due to their short-sighted nature. 
The strong branching rule \eqref{eq:sb}, while leading to potentially higher lower bound increases, requires solving multiple LPs, and hence is computationally expensive. 
While the pseudo-cost rule reduces the overhead compared to SB but is often found to be less effective in improving node efficiency. \emph{A powerful yet efficient approximation of SB without solving large amounts of LPs would be greatly beneficial.} Motivated by this point, we are developing efficient branching rules that imitate SB while the complexity is controlled by a budget $B$ (informal form):
\begin{equation}
    \label{eq:imitate-sb}
    \min_{\cV} \sum_{(\vc,\sX_{\mathrm{MILP}})\in \sM}~ \sum_{\mathrm{BnB~node}} \sum_{j} \left\| s_{j,\mathrm{SB}} - \cV(j,\cdot) \right\| ~~\text{subject to }~\mathrm{Complexity}(\cV) \leq B.
\end{equation}

\textbf{(Evaluating branching rules).} To evaluate the performance of a mapping $\cV$ for a given MILP instance, it is necessary to use a consistent node selection rule, for example, Depth-First Search (DFS), commonly adopted in contemporary MILP solvers. This enables the evaluation of the BnB performance given the mapping $\cV$ and its resulting variable selection strategy. In particular, it is recommended to embed a timer in the BnB algorithm and track the upper bound $\overline{f}$ and the lower bound $\underline{f}$ in real-time. These quantities at a specific moment $t$ can be denoted as $\overline{f}(t), \underline{f}(t)$. Assuming the BnB begins solving $(\vc, \sX_{\mathrm{MILP}})$ at time $t=0$ and finishes at $T_{\mathrm{end}}$, a typical metric can be expressed as follows:
\begin{equation}
    \label{eq:primal-dual-gap}
    I(\cV, (\vc, \sX_{\mathrm{MILP}}), \alpha) = \alpha \int_{0}^{T_{\mathrm{end}}} \left(\overline{f}(t) - \underline{f}(t)\right) \mathrm{d}t + T_{\mathrm{end}}.
\end{equation}
In BnB, the upper bound $\overline{f}$ decreases over time $t$ while the lower bound $\underline{f}$ increases. When the gap $\overline{f}(t) - \underline{f}(t)$ vanishes, we obtain an optimal solution of the MILP. Hence, we aim for the primal-dual gap to decrease as rapidly as possible.
This drives the first term in \eqref{eq:primal-dual-gap}. In addition, a scalar $\alpha$ balances the two terms in \eqref{eq:primal-dual-gap}. In contrast to the compromise goal \eqref{eq:imitate-sb}, a more ambitious goal is to seek a variable selection rule $\cV$ that minimizes the measure across a set of MILP instances $\sM$:
\begin{equation}
    \label{eq:minimize-gap}
    \min_{\cV} \sum_{(\vc,\sX_{\mathrm{MILP}}) \in \sM} I(\cV, (\vc, \sX_{\mathrm{MILP}}), \alpha). 
\end{equation}

\subsubsection{Learning to branch}
\label{sec:l2b}

Utilizing modern machine learning tools and optimization solvers, it is possible to quickly implement a pipeline targeting \eqref{eq:imitate-sb}. This pipeline can be universally applicable to any data set $\sM$, regardless of the need for specialized knowledge related to the specific form of the MILP instances in $\sM$. Although incorporating domain knowledge into machine learning may enhance performance, we are suggesting here the feasibility of an elementary model that operates without such knowledge. The pipeline consists of the following critical stages:
\begin{itemize}
    \item (Parameterization) Initially, one needs to define the input $\cdot$ of $\cV$ and establish the explicit form of $\cV$. There should be a sufficient number of parameters in the form of $\cV$ to allow a search for parameters that best fit \eqref{eq:imitate-sb}. Neural networks, due to their strong representation capacity, are a good choice to express $\cV$. (For further comprehension on parameterizing a mapping with a neural network, refer to Section \ref{subsec:dnn}) One has the flexibility to choose from various neural network structures or sizes (in terms of the number of parameters). Once the neural network's structure and size have been chosen, the complexity of invoking the neural network is limited by an upper bound, which satisfies the constraint in \eqref{eq:imitate-sb}.
    \item (Data collection) Following this, it is necessary to collect data for training during the BnB process. At each BnB node, we collect \emph{``(feature, label)" pairs}. The \emph{feature} refers to the specific values of the input $\cdot$ under the current circumstances, while the \emph{label} corresponds to the target we want the neural network's output to match. In this case, we can use the SB score as label.
    \item (Training and validation) This stage follows a standard supervised learning procedure. Given those ``(feature, label)" pairs, we can train a neural network to fit the strong branching. In the training process, there might be some hyperparameters that should be tuned by hand. We can use a separate validation set to help us choose the best hyperparameters.
    \item (BnB with trained models) Upon completion of the training, the trained neural network can be applied during the BnB process. Specifically, one can invoke the neural network to assign a score to each branching candidate in $\mathrm{frac}(\underline{\vx})$ in Algorithm \ref{algo:bnb}, and select the candidate with the highest score to branch on, just as in \eqref{eq:branch-with-score}.
\end{itemize}
The entire pipeline develops branching rules with machine learning, and is named \emph{Learning to Branch (L2B)}. In the subsequent paragraphs, we will individually discuss the components of this pipeline.

\paragraph{Parameterization and complexity.} In this section, we discuss several neural network-based scoring rules and compare their complexity to traditional score calculations like SB and PC.

\subparagraph{Approach 1: Separate models} 
The most straightforward parameterization method treats each variable as an independent entity. Features of variables are gathered individually and processed separately using a shared neural network. More specifically, the feature vector for each variable may include elements like $\vf_j = (c_j,l_j,u_j,z_j, \underline{x}_j,s_{j,\mathrm{PC}})$, where $c_j,l_j,u_j$ are as defined in \eqref{eq:milp}, $z_j$ signifies if $x_j$ is mandated to be an integer (where $z_j=1$ if $j\in\sI$ and $z_j=0$ otherwise), $\underline{x}_j$ is the $j$-th component of $\underline{\vx}$ as defined in Line 13 of Algorithm \ref{algo:bnb}, and $s_{j,\mathrm{PC}}$ represents the PC score. The ultimate goal is to build a machine learning model, such as an MLP \eqref{eq:mlp-simple}, that can map this feature vector to the strong branching score: 
\begin{equation}
    \label{eq:SB-imitation-approach1}
    \cV_{\mathrm{Sepa}}(j,\vc,\vl,\vu,\vz,\underline{\vx},\vs_{\mathrm{PC}}; \vtheta) := \mathrm{MLP}(\vf_j;\vtheta) \approx s_{j,\mathrm{SB}},
\end{equation}
where $\vtheta$ denotes the parameters in that neural network.

\subparagraph{Graph representations} To enrich the previous model, we will consider the relationship between variables and constraints in Approaches 2 and 3. To clearly describe this point, we adopt the language of graphs. We start by demonstrating how to model an MILP instance using a \emph{bipartite} graph, as shown in Figure \ref{fig:bipartite}. Variables are represented as nodes in the ``$\sV$" group, while constraints are represented as nodes in the ``$\sC$" group. No edges exist within each group. An edge with weight $A_{i,j}$ connects the $i$-th constraint to the $j$-th variable if $A_{i,j} \neq 0$. Each variable node carries a feature vector $\vf_j$, and each constraint node carries a feature vector $\vg_i$ that may include quantities related to that specific constraint, like $\vg_i = (\circ_i,b_i)$. (Recall \eqref{eq:milp} for the definition of ``$\circ,\vb$") For each variable node $j\in\{1,\cdots,n\}$, $\cN(j)$ is denoted as the set of its adjacent constraints. Conversely, $\cN(i)$ denotes the set of variables involved in the $i$-th constraint ``$\vA_{i,:}\vx \circ_i b_i$."

\begin{figure}[t]
  \begin{minipage}[b]{0.3\textwidth}
\[
\begin{aligned}
\min_{\vx} & ~ x_1 + 2x_2 + 3x_3 \\
\text{s.t.} & ~ x_1 + x_2 \geq 1\\
& ~ x_2 + 2x_3 \leq 2\\
& ~ x_1,x_2,x_3\in\{0,1\}
\end{aligned}
\]


\hspace{1cm}
\input{bipartite.tikz}

\vspace{0.5cm}

\captionof{figure}{An MILP instance represented by a bipartite graph.}
\label{fig:bipartite}
\end{minipage}
\hfill
\begin{minipage}[b]{0.68\textwidth}
\input{sepa.tikz}

\vspace{2mm}

\input{graph-aug.tikz}

\vspace{2mm}

\input{gnn.tikz}
\captionof{figure}{Illustration of three methods to represent $\cV(j,\cdot)$.}
\label{fig:gnn}
\end{minipage}
\end{figure}

\subparagraph{Approach 2: Features augmented with graph information} We can enrich $\vf_j$ by incorporating the relationships between variables and constraints. Rather than independently gathering features for each variable, we collect additional features that represent the neighborhood of $x_j$. These features can include statistical information about $\{\vg_i\}_{i\in\cN(j)}$, such as the distribution of constraints with $\circ_i=``\geq",``=",``\leq"$ respectively, and the mean and standard deviation of the set $\{b_i\}_{i\in\cN(j)}$. We can also consider more complex features including information about variables that are adjacent to $x_j$'s neighbors, like the mean and standard deviation of these variables' objective coefficients, lower bounds, and upper bounds. These values can be appended to $\vf_j$ to create an enhanced feature vector $\vf_{j,\mathrm{aug}}$, which provides a more comprehensive view of a variable. Using these enhanced features, we can more accurately predict the SB score: 
\begin{equation}
    \label{eq:SB-imitation-approach2}
    \cV_{\mathrm{Aug}}(j,\vc,\vl,\vu,\vz,\underline{\vx},\vs_{\mathrm{PC}}, \vA, \circ, \vb; \vtheta) := \mathrm{MLP}(\vf_{j,\mathrm{aug}};\vtheta) \approx s_{j,\mathrm{SB}}.
\end{equation}

\subparagraph{Approach 3: Graph neural network} Graph Neural Networks (GNNs) offer a way to represent and capture the essential information about the relationships between variables and constraints, without manual intervention. Specifically, we alternate between passing information from variables to constraints and vice versa, iteratively updating the feature vectors:
\begin{equation}
    \label{eq:gnn-message-passing}
    \begin{aligned}
        \vg_{i} \leftarrow & \mathrm{MLP}\left( \vg_{i}, \sum_{j \in \cN(i)} \mathrm{MLP}( \vg_{i}, \vf_{j}, A_{i,j}; \vtheta_1 ) ; \vtheta_2 \right), \\
\vf_{j} \leftarrow & \mathrm{MLP}\left( \vf_{j}, \sum_{i \in \cN(j)} \mathrm{MLP}( \vg_{i}, \vf_{j}, A_{i,j}; \vtheta_3 );\vtheta_4 \right)
    \end{aligned}
\end{equation}
After several iterations, the final resulting variable and node features are denoted as $\{\vf_{j,\mathrm{gnn}}\}_{j=1}^{n}$ and $\{\vg_{i,\mathrm{gnn}}\}_{i=1}^{m}$. An independent model is applied to each $\vf_{j,\mathrm{gnn}}$ to predict the SB score: $\mathrm{MLP}(\vf_{j,\mathrm{gnn}};\vtheta_5)\approx s_{j,\mathrm{SB}}$. For more information about the merits of GNNs, such as scalability and permutation invariance, please refer to Section \ref{sec:other-neural-nets}. The complete model can be abbreviated as
\begin{equation}
    \label{eq:SB-imitation-appraoch3}
    \cV_{\mathrm{GNN}}(j,\vc,\vl,\vu,\vz,\underline{\vx},\vs_{\mathrm{PC}}, \vA, \circ, \vb; \vtheta) := \mathrm{GNN}(j, \cG;\vtheta) \approx s_{j,\mathrm{SB}}. 
\end{equation}
Here $\cG = \big(\vA,\{\vf_j\}_{j=1}^n,\{\vg_i\}_{i=1}^{m}\big)$ represents the entire input to a GNN and $\vtheta = (\vtheta_1,\vtheta_2,\vtheta_3,\vtheta_4,\vtheta_5)$ includes all parameters in the GNN. A comparison between these three methods can be found in Figure \ref{fig:gnn}. Note that this diagram only displays one iteration of message passing (from variables to constraints and back), while in practice, multiple rounds might be employed.

\subparagraph{Complexity analysis} Initially, we will discuss the complexity of SB and PC. In SB, the primary computational overhead lies in solving $2 n_{\mathrm{frac}}$ LPs. There is still no polynomial-time algorithm for finding an exact LP solution. The simplex method, often used for LPs in BnB, may exhibit exponential complexity $\cO(2^n)$ in the worst case~\cite{bertsimas1997lp}. Despite typically not reaching this upper limit, the LP solving overhead is still substantial in real-world applications, with its complexity being somewhat unpredictable due to these potential worst-case outcomes. In summary, the cost of calculating the SB score is $\cO(\cC_{\mathrm{LP}} ~ n_{\mathrm{frac}})$, where $\cC_{\mathrm{LP}}$ denotes the complexity of solving a single LP. 
In contrast, the cost of PC is significantly lower. It consists of two parts: updating $\delta f_{j,\mathrm{up}}, \delta f_{j,\mathrm{down}}$ for the chosen $j$ and calculating $\delta_{j,\mathrm{up}}, \delta_{j,\mathrm{down}}$ for all $j \in \mathrm{frac}(\underline{\vx})$. The first part is dominated by standard BnB computations, and the second part is of $\cO(n_{\mathrm{frac}})$. Hence, the overall complexity of PC is simply $\cO(n_{\mathrm{frac}})$.
To simplify the analysis,  we will use $\cC_{\mathrm{MLP}}$ to denote the complexity of each neural network used in ML-based strategies. 
In Approach 1, we use an individual model for each candidate in $\mathrm{frac}(\underline{\vx})$, resulting in a complexity of $\cO(\cC_{\mathrm{MLP}}~n_{\mathrm{frac}})$. 
In Approach 2, when passing information from variables to constraints, the complexity for each constraint is proportional to its number of neighbors. Adding all these complexities, the total is $\cO(\mathrm{nnz}(\vA))$. The complexity of passing information from constraints to variables is the same. Thus, the complexity of Approach 2 is $\cO(L ~\mathrm{nnz}(\vA)) + \cO(\cC_{\mathrm{MLP}}~n_{\mathrm{frac}})$, where $L$ is the number of message-passing rounds. 
Using similar logic, the complexity of Approach 3 is $\cO(L ~\cC_{\mathrm{MLP}} ~\mathrm{nnz}(\vA)) + \cO(\cC_{\mathrm{MLP}}~n_{\mathrm{frac}})$. Here, the first term can obviously dominate the second one due to $\mathrm{nnz}(\vA)) \geq n \geq n_{\mathrm{frac}}$, so it can be expressed as $\cO(L ~\cC_{\mathrm{MLP}} ~\mathrm{nnz}(\vA))$. 
A summary of all complexities is provided in Table \ref{tab:complexity}. 

\begin{table}[H]
\centering
\caption{Summary of the complexities of branching rules}
\label{tab:complexity}
\begin{tabular}{l|l|l|l|l}
\hline
SB & PC & Approach 1 & Approach 2 & Approach 3 \\ \hline
$\cO(\cC_{\mathrm{LP}} ~ n_{\mathrm{frac}})$ & $\cO(n_{\mathrm{frac}})$ & $\cO(\cC_{\mathrm{MLP}}~n_{\mathrm{frac}})$ & $\cO(L ~\mathrm{nnz}(\vA)) + \cO(\cC_{\mathrm{MLP}}~n_{\mathrm{frac}})$ & $\cO(L ~\cC_{\mathrm{MLP}} ~\mathrm{nnz}(\vA))$ \\ \hline
\end{tabular}
\end{table}

Note that $\cC_{\mathrm{MLP}}$ is solely determined by the architecture and size of the neural network and can be fully controlled by the user. On the other hand, $\cC_{\mathrm{LP}}$ depends on the properties of the LP instance and the solver used, making it less straightforward for the user to control. \emph{This complexity explains our claim that neural networks can fulfill the constraint presented in \eqref{eq:imitate-sb}.} As we progress from Approach 1 to Approach 3, both the computational load and the power to fit the SB score increase. In practice, the choice of method should make a balance between efficiency and effectiveness.

\paragraph{Data collection.} Regardless of the approach chosen, it is necessary to collect sufficient data to train the chosen ML models. The collected data must correspond with the input structure of the ML model. Now let us take Approach 1 as an example. We maintain the data set $\sD_{\mathrm{Sepa}}$ in memory and initialize it as empty. Assuming we have a set of MILP instances $\sM$ available. During data collection, we iterate over all instances in $\sM$ and execute BnB (Algorithm \ref{algo:bnb}) for each instance. When Algorithm \ref{algo:bnb} invokes the variable selection rule, we perform the following steps:
\begin{itemize}
    \item Compute $\vf_j$ for all branching candidates $j \in \mathrm{frac}(\underline{\vx})$.
    \item Run SB and calculate the SB score $s_{j,\mathrm{SB}}$ for $j \in \mathrm{frac}(\underline{\vx})$
    \item Append these newly gathered data to $\sD_{\mathrm{Sepa}}$: 
    \[\sD_{\mathrm{Sepa}} \leftarrow \sD_{\mathrm{Sepa}} \cup \Big\{\big(\vf_j, s_{j,\mathrm{SB}}\big)\Big\}_{j\in \mathrm{frac}(\underline{\vx})}\]
\end{itemize}
Finally, the set $\sD_{\mathrm{Sepa}}$ becomes the training data set. The data collection procedure for Approach 2 aligns with that of Approach 1, with the only modification being the replacement of $\vf_j$ with $\vf_{j,\mathrm{aug}}$ and the re-naming of $\sD_{\mathrm{Sepa}}$ to $\sD_{\mathrm{Aug}}$. The procedure for Approach 3 slightly differs as the GNN's input requires features from the entire MILP, not just those of variables with fractional values. For Approach 3, the data set is labeled as $\sD_{\mathrm{GNN}}$, and the following steps are taken when BnB calls the variable selection rule:
\begin{itemize}
    \item Compute $\vf_j$ for all variables in the current MILP $j \in \{1,\cdots,n\}$.
    \item Compute $\vg_i$ for all constraints in the current MILP $i \in \{1,\cdots,m\}$.
    \item Run SB and compute the SB score $s_{j,\mathrm{SB}}$ for $j \in \mathrm{frac}(\underline{\vx})$
    \item Define $\cG = \big(\vA,\{\vf_j\}_{j=1}^n,\{\vg_i\}_{i=1}^{m}\big)$, where $\vA$ corresponds to the current MILP. Append these newly collected data to $\sD_{\mathrm{GNN}}$:
    \[\sD_{\mathrm{GNN}} \leftarrow \sD_{\mathrm{GNN}} \cup \Big\{ \big(\cG, \mathrm{frac}(\underline{\vx}), \vs_{\mathrm{SB}}\big)\Big\}.\]
\end{itemize}
Note that the data form presented here are simply the basic format. \emph{Users can enhance this framework by adding more features.} For example, in variable-related features $\vf_j$, the reduced cost (an essential concept in LP) can be included; in constraint-related features $\vg_i$, the dual variable can be added to indicate if the corresponding constraint is active. Such advanced feature engineering can be found in the relevant literature \cite{zarpellon2021param}. Another important trick is \emph{invoking SB with a probability}. Triggering SB at each BnB node might be extremely time-consuming. As such, we can call SB and gather data with a small likelihood, say $10\%$, and in other instances, resort to less expensive branching rules like PC. This tactic is discussed in \cite{ecole}.

\paragraph{Model training and validation.} Once we have the datasets $\sD_{\mathrm{Sepa}}$, $\sD_{\mathrm{Aug}}$, or $\sD_{\mathrm{GNN}}$ ready, we can proceed to the model training phase using standard machine learning procedures. For the first two approaches, the machine learning models can be trained by minimizing:
\begin{align}
    \min_{\vtheta} \ccL_{\mathrm{Sepa}}(\sD_{\mathrm{Sepa}},\vtheta) := \sum_{(\vf, s) \in \sD_{\mathrm{Sepa}}} \ell \left( s , \mathrm{MLP}(\vf; \vtheta) \right) \label{eq:training-1} \\
    \min_{\vtheta} \ccL_{\mathrm{Aug}}(\sD_{\mathrm{Aug}},\vtheta) := \sum_{(\vf, s) \in \sD_{\mathrm{Aug}}} \ell \left( s , \mathrm{MLP}(\vf; \vtheta) \right) \label{eq:training-2}
\end{align}
In these equations, $\ell$ represents the loss function, which should be manually chosen. A typical loss function might be the squared error: $\ell(a,b) = (a-b)^2$. For the GNN model, the loss function for training can be formulated as
\begin{equation}
\label{eq:training-3}
    \min_{\vtheta} \ccL_{\mathrm{GNN}}(\sD_{\mathrm{GNN}},\vtheta) := \sum_{(\cG, \sJ, \vs) \in \sD_{\mathrm{GNN}}} \sum_{j\in\sJ} \ell \left( s_j , \mathrm{GNN}(j, \cG; \vtheta) \right) 
\end{equation}
All these minimization formulas can be seen as specific forms of \eqref{eq:imitate-sb}. Given the complexity of neural networks can be directly controlled as demonstrated in Table \ref{tab:complexity}, the constraints in \eqref{eq:imitate-sb} can be met. Furthermore, minimization over the parameter space can be executed with ``back-propagation" which is built into many modern machine learning tools. See Section \ref{subsec:dnn-training} for a comprehensive introduction on training machine learning models. Note that certain hyperparameters, like the size of the neural networks and the learning rate, should be manually chosen during the training procedure. These hyperparameters are typically chosen through the \emph{validation} process. In addition to the MILP set used for training $\sM$, we prepare an independent MILP set $\sM'$, where we execute the same data collection procedure to result in the \emph{validation set} ${\sD'}_{\mathrm{Sepa}}$, ${\sD'}_{\mathrm{Aug}}$ or ${\sD'}_{\mathrm{GNN}}$. Using the validation set, we compute the loss function values $\ccL{\cdot}(\cdot, \vtheta_\ast)$ with $\vtheta_\ast$ as the parameters trained via \eqref{eq:training-1}, \eqref{eq:training-2} or \eqref{eq:training-3} with varying hyperparameters. The hyperparameters that yield the lowest loss function value on the validation set are selected. This is a standard validation procedure in machine learning. The model validated with a separate set $\sM'$ typically exhibits robust generalization ability, meaning it performs well on instances not included in the training set $\sM$.

\paragraph{Branching with training models.}  Once the model has been trained, it is straightforward to apply it. Given an MILP instance, we run BnB and trigger the trained models when BnB calls for the variable selection rule. Subsequently, each branching candidate is assigned a score, and we choose the candidate with the highest score for branching, as indicated by \eqref{eq:branch-with-score}.

\paragraph{Exploring new rules with reinforcement learning.} Instead of trying to imitate the strong branching, a more ambitious question arises: \textit{Can we develop a branching rule that directly aims to minimize the primal-dual gap integral \eqref{eq:minimize-gap}?} In comparison to fitting the strong branching, this method offers two potential advantages:
\begin{itemize}
    \item The possibility to outperform Strong Branching. Despite the thoughtful design behind the strong branching rule, it may not be optimal in terms of \eqref{eq:minimize-gap}. Directly devising branching rules based on \eqref{eq:minimize-gap} may lead to a more powerful and innovative rule.
    \item Freedom from reliance on the resource-intensive strong branching during the data collection phase. The strong branching can be exceedingly time-consuming on large-scale instances. Even a single BnB node might take longer than fully solving the MILP using economic branching rules like the pseudocost. This can result in inefficient data collection, while a training process that is independent of strong branching can remedy this issue.
\end{itemize}
Although the supervised learning pipelines discussed earlier in this section are not applicable to solve a minimization problem like \eqref{eq:minimize-gap} due to complex dependencies between the branching rule $\cV$ and the objective function $I(\cdot)$, this problem fits perfectly within the scope of \emph{reinforcement learning} (RL). First, we demonstrate that the Branch and Bound method (Algorithm \ref{algo:bnb}) is a well-defined, albeit complex, Markov Decision Process (MDP), assuming we use a fixed node selection rule like Depth-First Search (DFS). The key elements of the MDP can be defined as follows: (For introduction of MDP and RL, refer to Section \ref{sec:rl}.)
\begin{itemize}
\item State: All the data available when BnB calls for the branching rule: $\vc,\sX_{\mathrm{MILP}},\sX,\underline{\vx}, \underline{f}, \overline{\vx},\overline{f},\sS$.
\item Action: The action space consists of $\mathrm{frac}(\underline{\vx})$ and the action is the variable selected for branching.
\item Transition: Given the current state and action, execute BnB until it invokes the branching rule again. The resulting new state is the transitioned state and thus the transition probability is defined.
\item Reward $r$: If BnB does not stop, the reward is the negative of the primal-dual gap $r = -(\overline{f}-\underline{f})$ at that point; if BnB stops, the reward is the negative of the total solving time $r = -T_{\mathrm{end}}$. This definition will accurately reconstruct \eqref{eq:primal-dual-gap} when integrated over time.
\end{itemize}
Using these definitions, we can deploy RL to achieve the minimization outlined in \eqref{eq:minimize-gap}. Please refer to Section \ref{sec:rl} for an introduction to RL. Contrary to the supervised learning methods \eqref{eq:training-1},\eqref{eq:training-2},\eqref{eq:training-3}, where we first gather data using the strong branching rule in the BnB framework and then train ML models to fit the SB, RL approaches interact with the BnB: we obtain the state from BnB, make a decision, and send that action to BnB. BnB then returns a reward, indicating whether the lower bound has become tighter or whether BnB has stopped. Using this reward as a guide, we train the ML model. Both basic RL methods, value-based methods like Q-learning and policy-based methods like policy gradient, can be applied here, and we will delve into how these two methods can be utilized in learning to branch.

\subparagraph{Q-Learning} The Q-Learning algorithm aims at learning a scoring function, also known as the \emph{Q function}, which maps state-action pairs to scores. Given our definition of GNN model from \eqref{eq:SB-imitation-appraoch3}, we can consider the input $\cG$ as the state and $j$ as the action. Thus, $\mathrm{GNN}(j,\cG;\vtheta)$ can be seen as a parameterized Q function. Instead of imitating strong branching as in \eqref{eq:training-3}, Q-Learning seeks to make the Q function conform to the \emph{Bellman optimality equation}, which is an optimality condition for \eqref{eq:minimize-gap}:
\begin{equation}
    \label{eq:branching-q-learning}
    \min_{\vtheta} \sum_{(\cG,j,r,\cG')} \big( y - \mathrm{GNN}(j,\cG;\vtheta) \big)^2,~\text{where}~y = r + \max_{j'}\mathrm{GNN}(j',\cG';\vtheta_{\mathrm{prev}})
\end{equation}
Here, $\cG'$ represents the subsequent state given the current state $\cG$ and action $j$, and $\vtheta_{\mathrm{prev}}$ denotes the parameters established before the current training round. For instance, one can randomly initialize $\vtheta_{\mathrm{prev}}$, execute the BnB, and compile a collection of 4-tuples $\{(\cG,j,r,\cG')\}$. These are then trained using \eqref{eq:branching-q-learning} to derive $\vtheta$. This new $\vtheta$ then serves as $\vtheta_{\mathrm{prev}}$ in the next training round, eventually yielding $\vtheta_{\mathrm{next}}$, and so on. Upon obtaining the final parameters $\vtheta_{\ast}$, the corresponding GNN model $\mathrm{GNN}(\cdot,\cdot;\vtheta_{\ast})$ can be deployed as a branching score within BnB. When the BnB invokes the branching rule, the action $j_\ast = \max_{j \in \mathrm{frac}(\underline{\vx})} \mathrm{GNN}(j,\cG;\vtheta_{\ast})$ can be returned as the selected index to branch on.

\subparagraph{Policy Gradient} In contrast to Q-Learning, policy-based methods utilize a conditional probability distribution as the policy. This policy distribution is state-dependent and proposes a probability for all potential actions. The goal is to train the policy distribution to suggest actions reliably. By making a small modification, the GNN model from \eqref{eq:SB-imitation-appraoch3} can be treated as such a policy distribution: If we denote the outcomes of the GNN as $s_j$ for $j\in\{1,\cdots,n\}$, we can combine the scores of all branching candidates into a vector $(s_j)_{j\in\mathrm{frac}(\underline{\vx})}$ and apply a softmax function to that vector. The resulting output vector can then be seen as a conditional probability distribution and each entry in that vector represents the probability to select a specific index for branching. Such a probabilistic policy based on GNN is denoted as $p_{\mathrm{GNN}}(j|\cG;\vtheta)$.
Now, we use this probabilistic policy as the branching rule and execute the BnB on multiple MILP instances, resulting in a collection of trajectories. We gather the states, actions, and rewards from these trajectories. With this data at hand, we can compute the gradient of the function \eqref{eq:primal-dual-gap} with respect to the parameters in the policy and use this gradient to update the parameters, just like Section \ref{sec:rl}. After training, suppose the final parameters are $\vtheta_{\ast}$. We then use $p_{\mathrm{GNN}}(j|\cG;\vtheta_{\ast})$ as the branching rule: sample an index from this distribution and branch along that index.

\paragraph{Bibliographical notes.} The methodology of approximating SB via machine learning models can be traced back to 2014 \cite{alvarez2014supervised}. Initially, research on this topic yielded Approaches 1 and 2 presented in this paper \cite{alvarez2017machine,khalil2016learning}. From 2019 onward, innovative representations of MILP, including bipartite \cite{gasse2019exact} or tripartite representations \cite{ding2020accelerating}, as well as their related GNNs, have been introduced. Since then, such GNNs (described as Approach 3 in this paper) have been a fundamental model for ML applications in MILP. Recent research has explored alternative neural network architectures, such as transformers \cite{lin2022learning} and hierarchical structures \cite{zarpellon2021param}. More recent advances can be found in \cite{nair2020solving,gupta2020hybrid,gupta2022lookback}.
The exploration of reinforcement learning for crafting branching rules emerged later. Preliminary studies showcased that the performance of rules based on RL faced challenges in matching up to those derived by imitating SB. However, recent developments have significantly propelled this domain forward, making RL not only viable but also often outpacing imitation learning in certain contexts. Note that the methods presented in this section are only foundational setups that illustrate the potential of reinforcement learning for branching rules. To achieve competitive results, we recommend referring to recent literature \cite{etheve2020reinforcement,cappart2021combining,qu2022improved,scavuzzo2022learning,zhang2022deep,parsonson2023reinforcement}, which provides more advanced techniques.

\subsubsection{Node selection and learning to search}
\label{sec:l2s}
Revisiting Algorithm \ref{algo:bnb} for the BnB method, we can see that at Line 5, the BnB calls the node selection rule (also referred to as the \emph{searching rule}). This rule selects a node from the set of all active nodes $\sS$ (active nodes being those that have not been solved or pruned) and pops it from $\sS$. As illustrated in Figure \ref{fig:bnb-node-example}, this searching rule can impact the efficiency of BnB. 
We outline some standard searching strategies used in MILP solvers below. 
\begin{itemize}
    \item Best-First Search. This strategy selects the node with the smallest local lower bound. Looking at the global lower bound formula in Line 9 of Algorithm \ref{algo:bnb}: the global lower bound equals the smallest of all local lower bounds from nodes in $\sS$. By removing nodes with smaller local lower bounds from $\sS$ earlier, the global lower bound can be improved quickly.
    \item Breadth-First Search. This method selects the node with the least depth. Intuitively, nodes with lower depth usually possess a smaller local lower bound. Therefore, Breadth-First Search can also improve the global lower bound. However, Best-First Search is a more direct approach, and the Breadth-First Search technique offers little benefit over it.
    \item Depth-First Search. This method selects the next node candidate at the maximum depth in the tree. Unlike the best-first search, the depth-first search is more likely to produce feasible MILP solutions, often early in the search process, allowing for a faster improvement of the upper bound.
    \item Hybrid Search. This strategy involves switching between different search methods. For instance, one might initially employ the depth-first search to find a feasible solution quickly, and then transition to the best-first search to improve the dual bound.
\end{itemize}
To summarize, a searching rule in BnB is formulated as
\begin{equation}
    \label{eq:node-selection}
    (d, f, \sX) = \argmax_{(d', f', \sX') \in \sS} \cS\big((d', f', \sX'), \cdot \big)
\end{equation}
In this equation, the mapping $\cS$ assigns a score to each node, relying on its data $(d', f', \sX')$ and other beneficial information ``$\cdot$", such as the original MILP $(\vc,\sX_{\mathrm{MILP}})$. Then the node with the highest score $(d, f, \sX)$ is chosen for exploration, as outlined in Line 5 of Algorithm \ref{algo:bnb}.

\paragraph{Target analysis.}
The ultimate goal of the node selection rule is the same with the variable selection rule: to minimize the primal-dual gap as quickly as possible and reduce the size of the BnB tree. Formally, assuming a fixed variable selection rule such as the pseudocost rule, we seek the optimal node selection rule over a set of MILP instances $\sM$:
\begin{equation}
    \label{eq:node-selection-minimizing}
      \min_{\vS} \sum_{(\vc, \sX_{\mathrm{MILP}}) \in \sM} I(\cS, (\vc, \sX_{\mathrm{MILP}}), \alpha) := \alpha \int_{0}^{T_{\mathrm{end}}} \overline{f}(t) - \underline{f}(t) \mathrm{d}t + T_{\mathrm{end}}.
\end{equation}
Here, $\underline{f}(t),\overline{f}(t),T_{\mathrm{end}}$ have the same definitions as those in \eqref{eq:primal-dual-gap}. 
To further clarify the ambitious objective in \eqref{eq:node-selection-minimizing}, we discuss ``lowering $\overline{f}(t)$" and ``boosting $\underline{f}(t)$" separately.
\begin{itemize}
\item Improving the dual bound. Initially, we discuss how to improve the dual bound as rapidly as possible by choosing nodes in a suitable order. It can be demonstrated that \emph{the best-first search, previously discussed, is indeed optimal under this criterion}. Assume we have two nodes $(d, f,\sX)$ and $(d', f',\sX')$ in the set $\sS$ with $f < f'$. The best-first search proposes to explore $(d, f,\sX)$ first. If we decide to explore the other node first, even when $(d', f',\sX')$ yields an incumbent $\overline{\vx}$, the node $(d, f,\sX)$ cannot be pruned because $f < f' \leq \vc^\top \overline{\vx}$. This implies that investigating $(d, f,\sX)$ cannot be skipped under any circumstances. The dual bound can only be enhanced after the exploration of $(d, f,\sX)$, thereby demonstrating that the best-first search is indeed optimal under this measure.
\item Decreasing the primal bound. Contrary to improving the dual bound, no node selection rule that guarantees optimal primal bound reduction exists. None of the strategies introduced previously can assure the production of a feasible solution to the MILP and a consequent reduction of the primal bound. However, this presents us with an opportunity for improvement: one can design an effective strategy based on the domain-specific knowledge of a certain type of MILP, or develop a rule from data using machine learning techniques. 
\end{itemize}
Inspired by the above analysis, we will develop a strategy to reduce the primal bound rapidly. 
\emph{Intuitively, if we can identify the node in $\sS$ containing the optimal solution of MILP $\vx^\ast$ (i.e., finding the node $(d,f,\sX)\in\sS$ such that $\vx^\ast \in \sX$), selecting that node will be a logical choice.} Branching along this node may yield incumbents $\overline{\vx}$ with a lower objective function and will ultimately provide the best primal bound because that node comprises the feasible solution with the smallest objective. Therefore, an ideal node scorer can be defined as:
\begin{equation}
    \label{eq:node-score-target}
    \cS_{\mathrm{ref}}\Big((d,f,\sX),\vx^\ast \Big) = \begin{cases}
        1,~~\text{if }\vx^\ast \in \sX \\
        0,~~\text{Otherwise}
    \end{cases}
\end{equation}
However, this mapping is not directly applicable as a node scorer in BnB due to its dependence on the optimal solution $\vx^\ast$. Therefore, it can be advantageous to develop an efficient mapping that approximates $\cS_{\mathrm{ref}}$ without relying on $\vx^\ast$. Formally, this can be written as:
\begin{equation}
    \label{eq:find-node-scorer}
        \min_{\cS} \sum_{(\vc,\sX_{\mathrm{MILP}})\in\sM} \sum_{(d,f,\sX) \in \sS} \left\| \cS\Big((d,f,\sX),(\vc,\sX_{\mathrm{MILP}})\Big) - \cS_{\mathrm{ref}}\Big((d,f,\sX),\vx^\ast \Big) \right\| 
        ~\text{s.t. Complexity}(\cS) \leq B
\end{equation}
Here, $\sS$ denotes the active node set generated during the BnB process on the MILP instance $(\vc,\sX_{\mathrm{MILP}})$ and $\vx^\ast$ is one of the optimal solutions to $(\vc,\sX_{\mathrm{MILP}})$.

Upon examining the formula in \eqref{eq:find-node-scorer}, one might ask: \emph{is \eqref{eq:find-node-scorer} essentially the same as solving an MILP?} We argue that finding the mapping $\cS$ in \eqref{eq:find-node-scorer} is significantly less challenging than directly finding a solution to an MILP, especially in the early stage of BnB. 
To illustrate this, consider the scenario where we have just explored the root node of the BnB tree. Based on Algorithm \ref{algo:bnb}, there will be two child nodes in $\sS$: $(1,\vc^\top\underline{\vx},\sX_{\mathrm{up}})$ and $(1,\vc^\top\underline{\vx},\sX_{\mathrm{down}})$. Now we must decide which node to explore first. As discussed earlier, it would be beneficial to explore the node containing the final optimal solution $\vx^\ast$.  
The goal of $\cS$ in \eqref{eq:find-node-scorer} is merely to predict whether $\vx^\ast$ is in $\sX_{\mathrm{up}}$ or $\sX_{\mathrm{down}}$, without needing to determine the exact value of $\vx^\ast$. In common scenarios, both $\sX_{\mathrm{up}}$ and $\sX_{\mathrm{down}}$ are large sets containing many feasible solutions because BnB is in the early stage, \emph{so determining $\vx^\ast$ in either of these sets can only provide a rough estimate of $\vx^\ast$, far from determining the precise location.} Thus, the mapping $\cS$ in \eqref{eq:find-node-scorer} is much less complex than solving an MILP, yet useful in BnB node selection. 
In the following paragraph, we will discuss how \eqref{eq:find-node-scorer} can be effectively implemented using machine learning.

\paragraph{Learning to search.} Developing a machine learning model to act as a searching rule is referred to as \emph{learning to search (L2S)}. Similar to L2B, as discussed in Section \ref{sec:l2b}, L2S follows the same pipeline: parameterization, data collection, model training and validation, and the application of the trained model within the BnB process.
\subparagraph{Parameterization} Here we present a method to express all the inputs of $\cS$ ``$(d,f,\sX),(\vc,\sX_{\mathrm{MILP}})$" in a format compatible with a GNN. This enables us to parameterize $\cS$ with a GNN. 
To start, consider the original MILP $(\vc,\sX_{\mathrm{MILP}})$. As depicted in Figure \ref{fig:bipartite}, each variable and each constraint can be regarded as a vertex, thus conceptualizing the MILP as a bipartite graph. 
The feasible set $\sX_{\mathrm{MILP}}$ is determined by $(\vA,\circ,\vb,\vl,\vu,\vz)$ where $\vz \in \{0,1\}^{n}$ denotes the variable type (i.e., $z_j=1$ for $j\in\sI$ and $z_j=0$ otherwise).  
Here, $(\vc,\vl,\vu,\vz)$ can be viewed as features to characterize a variable, and $\circ,\vb$ characterizes a constraint. The connections between variables and constraints are expressed by the matrix $\vA$. 
Compared to the original feasible set $\sX_{\mathrm{MILP}}$, the feasible set corresponding to a BnB node $\sX$ only differs in terms of $\vl,\vu$, as BnB only tightens the lower or upper bound of a variable (see Line 28 in Algorithm \ref{algo:bnb}). We denote the new variable bounds as $\vl'$ and $\vu'$.
Finally, the variable feature is defined with $\vf_j = (c_j,l_j,u_j,z_j,{l'}_j,{u'}_j)$ and the constraint feature is defined with $\vg_i = (\circ_i,b_i)$. We then apply the message passing procedure defined in \eqref{eq:gnn-message-passing} for several rounds to obtain updated feature vectors $\{\vf_{j,\mathrm{gnn}}\}_{j=1}^n$ and $\{\vg_{i,\mathrm{gnn}}\}_{i=1}^m$ and propagate them together, applying a MLP as $\mathrm{MLP}(\sum_j \vf_{j,\mathrm{gnn}}, \sum_i \vg_{i,\mathrm{gnn}}, d, f; \vtheta_{6})$, of which the output is a real number. This result can be regarded as the probability of $\vx^\ast\in\sX$. We use $\cG = (\vA,\{\vf_j\}_{j=1}^n, \{\vg_i\}_{i=1}^m,d,f)$ to represent the input of the GNN. The entire model can be denoted as:
\begin{equation}
    \label{eq:gnn-node-selection}
    \cS_{\mathrm{GNN}}(\cG; \vtheta ) := \mathrm{MLP}\left(\sum_{j=1}^n \vf_{j,\mathrm{gnn}}, \sum_{i=1}^m \vg_{i,\mathrm{gnn}}, d, f; \vtheta_{6}\right) \approx \cS_{\mathrm{ref}}\Big((d,f,\sX),\vx^\ast \Big)
\end{equation}
In the equation above, $\vtheta = (\vtheta_1,\vtheta_2,\vtheta_3,\vtheta_4,\vtheta_6)$, with $\vtheta_1,\vtheta_2,\vtheta_3,\vtheta_4$ being defined in \eqref{eq:gnn-message-passing}. This is a graph-level GNN rather than a node-level GNN. (For definition of types of GNN, refer to Section \ref{sec:gnn}.)

\subparagraph{Data collection and model training}  The processes of data collection, model training, and validation closely follow the pipeline outlined for L2B in Section \ref{sec:l2b}. For each MILP instance in the training set, we run BnB, and when BnB calls on the node selection rule, we gather data in the form of $\cG$. This data is then appended to the dataset $\sD_{\mathrm{Input}}$, which is initialized as an empty set:
\[ \sD_{\mathrm{Input}} \leftarrow \sD_{\mathrm{Input}} \cup \{\cG\} \]
In $\sD_{\mathrm{Input}}$, each element corresponds to a node in the BnB tree. Once BnB has been completed for an MILP instance and an optimal solution $\vx^\ast$ is obtained, we can determine whether $\vx^\ast$ is included in each BnB node, enabling us to assign labels to each node. If the node corresponding to $\cG$ includes $\vx^\ast$ (i.e., $\cS_{\mathrm{ref}}\Big((d,f,\sX),\vx^\ast \Big)=1$), we set $s_{\cG}=1$, otherwise $s_{\cG}=0$. The dataset for training a search strategy is then defined as:
\[ \sD_{\mathrm{Search}} = \{(\cG,s_{\cG})\}_{\cG \in \sD_{\mathrm{Input}}}. \]
With the above data set, we train the GNN with
\[ \min_{\vtheta} \sum_{(\cG,s_{\cG})\in\sD_{\mathrm{Search}}} \left| \cS_{\mathrm{GNN}}(\cG; \vtheta ) - s_{\cG} \right| \]
For information regarding model validation and hyperparameter tuning, please refer to Section \ref{sec:l2b}.

\subparagraph{ML-based node selection} Once the training phase is complete, the trained GNN model, denoted as $\cS_{\mathrm{GNN}}(\cG; \vtheta )$, assigns a score to each BnB node. When BnB calls upon the node selection rule, we select the node with the highest score. If $\cS_{\mathrm{GNN}}$ is a good approximation of $\cS_{\mathrm{ref}}$, as previously discussed, this node selection strategy is advantageous for reducing the primal bound and can be employed during the early stages of BnB. After several iterations of BnB, it may be beneficial to transition to the best-first search strategy, as previously mentioned.

\paragraph{Bibliographical notes.} The presented L2S framework serves as a basic configuration to demonstrate its potential applicability. In practice, to achieve more competitive outcomes, neural networks are typically employed to compare two given BnB nodes rather than assigning a score to a single node~\cite{he2014learning,khalil2022mip}. For other efforts on ML-based node selection, please refer to the designated reference \cite{labassi2022learning,yilmaz2021study}.

\subsection{Cutting Plane Methods}
\label{sec:cut}
Besides the branch-and-bound method, the \emph{cutting plane method} is another important algorithm in solving MILP. Rather than employing the BnB approach that generates and solves multiple sub-MILPs, the cutting plane method progressively incorporates linear constraints into the root LP relaxation. We first introduce some basic concepts in the cutting plane method and then present how to utilize machine learning in it.

\subsubsection{Introduction to cutting plane methods}
\begin{figure}[ht]
    \centering
    \input{cuts.tikz}
    \caption{An example of cuts. The MILP feasible set is defined with $\sX_{\mathrm{MILP}} = \{(x_1,x_2): x_1,x_2 \in \sZ, 0 \leq x_1 \leq 2.5, 0 \leq x_2 \leq 2, x_1 + x_2 \leq 3.5\}$. Suppose we run the simplex algorithm and obtain an LP optimal solution $\underline{\vx}=(1.5,2)$. To tighten the LP relaxation, one can introduce an extra constraint. In this example, both Cut 1 and Cut 2 can fulfill this purpose, as they both separate $\underline{\vx}$ from all feasible solutions $\vx \in \sX_{\mathrm{MILP}}$. In other words, both cuts will successfully eliminate $\underline{\vx}$ while preserving the feasible set $\sX_{\mathrm{MILP}}$. However, the efficiency of these two cuts differs. Using Cut 1, the LP solution becomes $(2,1.5)$, leading to a lower bound of $-3.5$. With Cut 2, an LP solution is $(1,2)$ and the lower bound is $-3$. Clearly, Cut 2 is better as it provides a greater enhancement to the lower bound.}
    \label{fig:cut}
\end{figure}
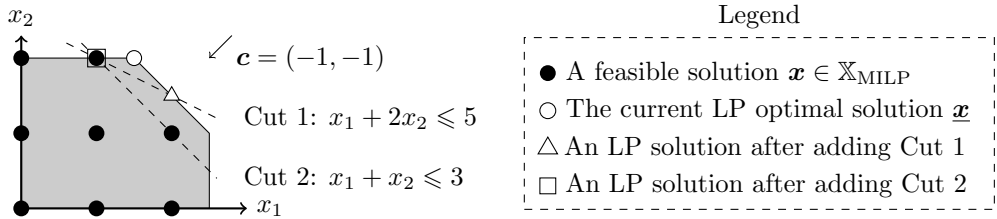

\paragraph{Cuts and separators.} Here we adopt the notation outlined in Algorithm \ref{algo:bnb}. For a given BnB node $(d,f,\sX)$, we solve the LP relaxation $(\vc, \mathrm{Relax}(\sX))$. For simplicity, \emph{we assume such an LP is feasible and bounded}, and hence we can achieve its optimal solution $\underline{\vx}$. Suppose $\underline{\vx}$ violates the integer constraint: $\mathrm{frac}(\underline{\vx}) \neq \emptyset$ and therefore $\underline{\vx} \not \in \sX$. If there exist $\valpha = (\alpha_1,\cdots,\alpha_n) \in \sQ^{n}$ and $\beta \in \sQ$ such that
\[\valpha^\top \vx \leq \beta,~\forall \vx \in \sX, \qquad \valpha^\top \underline{\vx} > \beta \]
then we term $(\valpha,\beta)$ as a \emph{cut} or a \emph{separator} of this BnB node. Since $\underline{\vx} \not \in \sX$, adding the additional constraint $\valpha^\top\vx \leq \beta$ to $\sX$ will cut off certain points (including at least $\underline{\vx}$) from $\mathrm{Relax}(\sX)$, while preserving the MILP feasible set $\sX$:
\[ \sX = \sX \cap \{\vx: \valpha^\top\vx \leq \beta\}, \qquad  \mathrm{Relax}(\sX) \subsetneq \mathrm{Relax}(\sX) \cap \{\vx: \valpha^\top\vx \leq \beta\} \]
Let us denote $\underline{\vx}^\prime$ as the optimal solution to the LP relaxation after adding the cut. Surely it holds that $\vc^\top \underline{\vx}^\prime \geq \vc^\top \underline{\vx}$. Therefore, such a cut can tighten the LP relaxation and hopefully improve the lower bound. Intuitively, a cut separates $\underline{\vx}$ from $\sX$ and Figure \ref{fig:cut} visually demonstrates this concept.

\paragraph{Generating cuts.} Indeed, the literature reveals an extensive array of potential cut candidates applicable to an MILP. In this context, we introduce the Chvátal-Gomory cut (CG cut) as a representative example. Consider a scenario in which we deal with a pure integer program where all variables are required to be non-negative integers. In the context of equation \eqref{eq:milp}, we assume $\sI = \{1,\cdots,n\}$ and $\vl = \vzero, \vu \in \sZ_{+}^n$. Let us assume that the $i$-th constraint of the MILP takes the form $\va_{i}^\top \vx \leq b_i$, with $\va_i$ representing the $i$-th row of the matrix $\vA$ and $b_i$ being the $i$-th element in $\vb$. From this constraint, we can derive a cut through the following rounding technique:
\begin{equation}
    \label{eq:cg-cut}
    \lfloor \va_{i} \rfloor^\top \vx \leq \lfloor b_i \rfloor
\end{equation}
Here, $\lfloor \va_{i} \rfloor$ signifies the element-wise floor operation. Clearly any feasible integer solution $\vx \in \sX$ must satisfy \eqref{eq:cg-cut}: $\lfloor \va_{i} \rfloor^\top \vx \leq \va_{i}^\top \vx \leq b_i$ due to the non-negativity of $\vx$, and $\lfloor \va_{i} \rfloor^\top \vx$ must be an integer, which implies \eqref{eq:cg-cut}. Thus, adding the constraint \eqref{eq:cg-cut} to $\sX$ will not cut off any desirable solutions. Given $\underline{\vx}$, it is straightforward to examine if $\lfloor \va_{i} \rfloor^\top \underline{\vx} > \lfloor b_i \rfloor$ is true. If this condition is met, \eqref{eq:cg-cut} will serve as a valid cut that separates $\underline{\vx}$ from $\sX$. Through this methodology, one may potentially create as many as $m$ cuts (equivalent to the number of constraints), a significant amount. Further, there exists a wide variety of practical techniques for generating cuts, such as strong CG-cut, MIR-cut, lift-and-project cut, and so on.\cite{wolsey2020integer} Each of these can give rise to a substantial quantity of cuts applicable to an MILP.  While this tutorial refrains from diving into the vast scope of papers on this subject, it may offer a typical number of cuts as a guide for newcomers. Given an MILP of size $\vA \in \sQ^{500\times1000}$, the number of all the potential cut candidates might exceed $10^4$.

\paragraph{Selecting cuts.} Given the large number of potential cuts, it is impractical to incorporate all potential cutting plane candidates. The reason is two-fold: (1) the resulting LP relaxation might be significantly enlarged and subsequently challenging to solve; (2) many cut candidates cannot significantly improve the lower bound. The first reason is relatively straightforward, while the second one can be illustrated by Figure \ref{fig:cut}, where Cut 1 is obviously undesirable as it cannot improve the lower bound. Therefore, it becomes vital to select part of the cut candidates and then resolve the subsequent LP.

\paragraph{Cutting plane methods.} The procedure detailed above can be executed iteratively. Specifically, we first generate cuts, then  select from among them, add the chosen cuts to $\sX$, and resolve the LP relaxation. This sequence of actions is repeated continuously until either a satisfactory solution is discovered or the designated budget is reached. This entire iterative approach is known as the \emph{cutting plane method} or \emph{separation method}, which is summarized in Algorithm \ref{algo:cut}. 

\begin{algorithm}[ht]
\renewcommand{\algorithmicrequire}{\textbf{Input:}}
\renewcommand{\algorithmicensure}{\textbf{Output:}}
\caption{Cutting Plane Algorithm}
\label{algo:cut}
\begin{algorithmic}[1]
    \REQUIRE A feasible and bounded (sub-)MILP $(\vc, \sX)$; Cut selection rule;
    \ENSURE An optimal solution to $(\vc, \sX)$ or a tightened feasible set $\sX^\prime$;
    \STATE Initialize $\sX^\prime \leftarrow \sX$;
    \WHILE{The time budget is not reached}
        \STATE Solve $(\vc, \mathrm{Relax}(\sX^\prime))$ to obtain its optimal solution $\underline{\vx}^\prime$;
        \IF{$\underline{\vx}^\prime \in \sX$ (i.e., all the integer constraints are satisfied)}
            \RETURN $\underline{\vx}^\prime$ as the optimal solution;
        \ENDIF
        \STATE Generate a bag of cuts $\{(\valpha_i,\beta_i)\}_{i=0}^{m'}$ such that $\valpha_i^\top \vx \leq \beta_i$ for all $\vx \in \sX$, and $\valpha_i^\top \underline{\vx} > \beta_i$;
        \STATE Select a subset of cuts from $\{(\valpha_i,\beta_i)\}_{i=0}^{m'}$ according to the cut selection rule; \COMMENT{The selected cuts form a new set $\{(\valpha_i,\beta_i)\}_{i\in \sI'}$.}
        \STATE Update $\sX^\prime$: $\sX^\prime \leftarrow \sX^\prime \cap \left(\cap_{i\in \sI'} \{\vx: \valpha_i^\top \vx \leq \beta_i\}\right)$;
    \ENDWHILE
    \RETURN $\sX^\prime$
\end{algorithmic}
\end{algorithm}

If this algorithm is directly applied to the original MILP $(\vc,\sX_{\mathrm{MILP}})$, and an ample computational budget is allocated, integer programs with rational data can theoretically be solved within a finite number of steps, provided that the cuts are appropriately selected (as demonstrated by Gomory~\cite{gomory1958algorithm,gomory1960solving}). However, a pure cutting plane method usually requires an extremely large number of iterations to produce the optimal solution in practice, making it far from efficient compared with the BnB method. Such a situation remained unchanged until Balas et al.~\cite{balas1996gomory} showed that the cutting plane method is efficient when combined with the BnB. State-of-the-art MIP solvers typically combine BnB and cutting planes with one of the following approaches~\cite{achterberg2007constraint}:
\begin{itemize}
    \item Cut and Branch. We first apply Algorithm \ref{algo:cut} on the original MILP $(\vc,\sX_{\mathrm{MILP}})$ with a controlled budget and obtain a tightened feasible set $\sX^\prime$. Afterwards, the problem $(\vc,\sX^\prime)$ is solved with BnB (Algorithm \ref{algo:bnb}). With the global lower bound already improved through the cutting planes, this cut-and-branch procedure often concludes earlier than the standard BnB method, without the necessity of traversing the entire BnB tree.
    \item Branch and Cut. The MILP is solved with BnB, while all BnB nodes might be tightened by the cutting plane method. The cuts included at the root node (the original MILP) are considered global cuts, whereas those integrated into a sub-MILP are viewed as local cuts. Local cuts, valid solely for the sub-MILP and its subsequent child nodes, must be discarded when the exploration extends beyond that specific subtree.
\end{itemize}
In this paper, we only consider cut-and-branch for the sake of simplicity. However, this narrowed scope does not significantly deviate from the practices commonly employed in real-world applications. Modern MILP solvers utilizing the branch-and-cut method typically allocate more computational resources for cuts in the original MILP (or the root node) compared to those in a sub-MILP, due to the global applicability of the cuts added at the root node.

\subsubsection{Learning to select cuts}
\label{sec:l2c}

Similar to the rules for selecting branching variables and BnB nodes, the rule for cut selection can be stated as: assign a score to every potential cut and choose the topmost $\gamma\%$ of them. Generally, the score assigned to a cut is influenced not only by the cut $(\valpha_i,\beta_i)$ itself, but also by the MILP data $(\vc,\sX)$. The scoring function can be expressed as:
\begin{equation}
    \label{eq:cut-score}
    s_{i,\mathrm{cut}} = \cC((\vc,\sX), (\valpha_i,\beta_i)).
\end{equation}
Based on this equation, we want to discover a mapping $\cC$ that assigns high ratings to effective cuts and low ratings to ineffective ones. It is worth noting that the cut selection proportion $\gamma\%$ is another critical aspect that can be optimized using machine learning. However, in this context, we treat it as a manually tunable hyperparameter for simplicity's sake. Like the variable selection rule stated in Section \ref{sec:l2b}, the cut selection strategy can also be parameterized by neural networks and trained via both supervised and reinforcement learning approaches.

\paragraph{Supervised learning.} We present an exemplary application of supervised learning for the selection of cuts, drawing inspiration from the principles of strong branching. In SB, rather than immediately proceeding with a branch, the solver simulates branching on various candidates to assess their potential effect on the objective function. Similarly, in cut selection, the most impactful cut, in terms of objective improvement, is a good choice. By incorporating a cut of the form $\valpha_i^\top \vx \leq \beta_i$ into the set $\sX$ and evaluating the modified LP relaxation, the resulting shift in the objective value provides a valuable measure of the cut's effectiveness. Specifically, this idea can be formulated as:
\begin{equation}
    \label{eq:look-ahead-cut}
    s_{i,\mathrm{cut}} := \frac{\Delta_i}{\Delta_{i_\ast}},\quad \text{where }\Delta_i = \min_{\vx \in \mathrm{Relax}(\sX) \cap \{\vx: \valpha_i^\top \vx \leq \beta_i\} } \vc^\top \vx - \vc^\top \underline{\vx}  ~~\text{and }i_\ast = \argmax \Delta_i
\end{equation}
Here, $s_{i,\mathrm{cut}}$ represents a score assigned to the $i$-th cut $\valpha_i^\top \vx \leq \beta_i$, quantifying its effectiveness. Such a ``look-ahead" cut-scoring function maximizes the dual bound improvement greedily, but it is computationally demanding as it requires solving multiple LPs - one for each proposed cut. Considering that the number of potential cuts even might exceed the count of linear constraints in the original MILP, the computational load due to LP evaluations is notably expensive. Still, there's potential to \emph{leverage neural networks as fast approximators of the scoring function} delineated in \eqref{eq:look-ahead-cut}. This mirrors the strategy where neural networks have been utilized to imitate SB, as expounded in Section \ref{sec:l2b}.

\subparagraph{Parameterization} Here we present a GNN-based method to parameterize the scoring mapping, $\cC$, that maps $(\vc,\sX),(\valpha_i,\beta_i)$ to $s_{i,\mathrm{cut}}$. Drawing inspiration from Figure \ref{fig:bipartite} in Section \ref{sec:l2b}, we consider each variable of the MILP as a node. Similarly, every constraint, including the additional constraints $\{\valpha_i^\top \vx \leq \beta_i\}_{i=1}^{m'}$, is also treated as a node. For variable nodes, we define their feature vector as $\vf_j = (c_j,l_j,u_j,z_j,\underline{x}_j)$. For constraint nodes, the associated feature vector is captured by $\vg_i = (\circ_i,b_i)$ for the range $1 \leq i \leq m$ and $\vg_i = (\leq, \beta_{i-m})$ for the range $m+1 \leq i \leq m+m'$. 
Given the additional constraints, we must expand the matrix $\vA\in\sQ^{m \times n}$ to $\vA_{\mathrm{cut}}\in\sQ^{(m+m')\times n}$. The initial $m$ rows of $\vA_{\mathrm{cut}}$ are identical to $\vA$. The succeeding $m'$ rows are defined by $\vA_{i,:,\mathrm{cut}}:=\valpha_{i-m}$ for all $m+1 \leq i \leq m+m'$. 
With the concepts defined above, we can apply the message passing rule \eqref{eq:gnn-message-passing} for several rounds and obtain the outcomes $\{\vf_{j,\mathrm{gnn}}\}_{j=1}^n$ and $\{\vg_{i,\mathrm{gnn}}\}_{i=1}^{m+m'}$. Finally, a dense neural network is employed individually on $\{\vg_{i,\mathrm{gnn}}\}_{i=m+1}^{m+m'}$ to fit the cut scores in \eqref{eq:look-ahead-cut}: specifically, $\mathrm{MLP}(\vg_{i,\mathrm{gnn}}; \vtheta_8) \approx s_{i-m, \mathrm{cut}}$. Summarizing, the entire GNN model can be represented as:
\[ \cC_{\mathrm{GNN}}\left( i , \cG_{\mathrm{cut}}; \vtheta \right) := \mathrm{MLP}(\vg_{i+m,\mathrm{gnn}}; \vtheta_8) \approx s_{i, \mathrm{cut}},\quad 1 \leq i \leq m'. \]
Here, $\cG_{\mathrm{cut}}=(\vA_{\mathrm{cut}}, \{\vf_j\}_{j=1}^n,\{\vg_i\}_{i=1}^{m+m'})$ serves as the GNN input while $\vtheta = (\vtheta_1,\vtheta_2,\vtheta_3,\vtheta_4,\vtheta_8)$ represents all the parameters in the GNN, with $\vtheta_1,\vtheta_2,\vtheta_3,\vtheta_4$ being defined in \eqref{eq:gnn-message-passing}.

\subparagraph{Data collection and training} For each MILP instance in the training set, we collect $\cG_{\mathrm{cut}}$ and calculate $\{s_{i, \mathrm{cut}}\}_{i=1}^{m'}$ based on \eqref{eq:look-ahead-cut}, each tuple $(\cG_{\mathrm{cut}}, \{s_{i, \mathrm{cut}}\}_{i=1}^{m'})$ is treated as an individual sample. 
Following this, we identify the top $\gamma\%$ of cuts based on their scores, integrate them into the MILP, and thereby generate a tightened MILP instance. For this new instance, another sample is collected. By iterating this procedure across multiple cycles, we accumulate a series of training samples. After processing all the MILP instances in the training dataset and consolidating the gathered samples, our compiled training data is denoted as: $\sD_{\mathrm{cut}}$. Utilizing this dataset, our aim is to optimize the GNN in order to best predict the cut scores:
\[ \min_{\vtheta} \sum_{(\cG_{\mathrm{cut}}, \{s_{i, \mathrm{cut}}\}_{i=1}^{m'}) \in \sD_{\mathrm{cut}}} \sum_{i=1}^{m'} \ell\left( \cC(i,\cG_{\mathrm{cut}};\vtheta) , s_{i, \mathrm{cut}} \right) \]
Within this context, the loss function, denoted as $\ell$, might be the mean squared error or perhaps the cross-entropy as seen in \ref{eq:gnn-cross-entropy}. Other procedures, like model validation and the tuning of hyperparameters, align with conventional machine learning methodologies. Please refer to Section \ref{sec:l2b} for details. Upon successful training, the trained GNN serves as the surrogate for the mapping $\cC$ as presented in \eqref{eq:cut-score}, equipping us with a GNN-driven strategy for cut selection.

\paragraph{Reinforcement learning.} Here we present an alternate method utilizing reinforcement learning (RL) for cut selection. We start by outlining the motivation behind the application of RL. While the cut selector, derived from \eqref{eq:look-ahead-cut}, might appear optimal for a single round of cutting planes, it is actually a local optimization strategy. Such a greedy approach might not necessarily produce the globally optimal cut selection strategy as defined by:
\begin{equation}
    \label{eq:rl-goal-cut}
    \max_{\cC} \sum_{(\vc,\sX) \in \sM} I(\cC, (\vc,\sX) ) := \int_{0}^{T_{\mathrm{end}}} \underline{f}(t) \mathrm{d}t
\end{equation}
Here $\sM$ aggregates all conceivable (sub-)MILPs encountered during cutting plane methods, with $\underline{f}(t)$ representing the lower bound at a specific time $t$, and $T_{\mathrm{end}}$ marking the termination of Algorithm \ref{algo:cut}. Termination can result from achieving optimality or due to budget control. To obtain an effective cut selection rule in the sense of \eqref{eq:rl-goal-cut}, one can adopt reinforcement learning. By presuming the selection of the single best cut in every cycle, Algorithm \ref{algo:cut} can naturally be formulated as an MDP: 
\begin{itemize}
    \item State: The previously defined $\cG_{\mathrm{cut}}$ describes the current state. 
    \item Action: The cut to be selected. We use $i \in \{1,\cdots,m'\}$ to denote its index.
    \item Transition: Upon selecting a cut, it is added into the MILP, resulting in the tightened feasible set $\sX^\prime$ as Line 9 in Algorithm \ref{algo:cut}.
    \item Reward: This is dictated by the dual bound, with higher values being more desirable.
\end{itemize}
With these definitions, methods such as Q-learning and policy gradient methods, as highlighted in Section \ref{sec:l2b}, can be conveniently adopted for enhancing cut selection. For example, when using Q-learning, one can simply replace $j$ with $i$ and the policy $\mathrm{GNN}(j,\cG;\vtheta)$ with $\cC_{\mathrm{GNN}}(i,\cG_{\mathrm{cut}};\vtheta)$ in equation \eqref{eq:branching-q-learning}. To adopt the policy gradient, one may treat $\cC_{\mathrm{GNN}}$ as a stochastic policy rule, where the output corresponds to the likelihood of selecting a particular cut candidate. However, it is worth noting that \emph{the action space for cut selection is much larger than that of variable selection} so that its training is also much more challenging. To provide context: for an MILP instance described by $\vA\in\sQ^{500\times1000}$, the action space for variable selection consists of all the fractional variables, typically of size $20\sim50$. While in the cut selection, the action space involves all the cut candidates, typically of size that may exceed $10^4$ as we mentioned before. Expanding the number of cuts in a given round further inflates this action space. For example, the size of the action space will be $\binom{10^4}{10}$ if we select $10$ cuts each time.

\paragraph{Bibliographical notes.} The incorporation of machine learning into the cutting plane method emerged relatively late compared to the evolution of learning to branch. Several factors contribute to this delay. Firstly, unlike the straightforward criterion of strong branching for variable selection, cut selection lacks a universally accepted standard. Secondly, as previously discussed, the action space for cut selection is considerably larger. Nevertheless, with the development of deep reinforcement learning, there's been a growing interest in learning for cut selection. The discussions on supervised learning in this section are mainly inspired by \cite{paulus2022learning,huang2022learning}, notable for their analogy to strong-branching imitation. In parallel, the reinforcement learning aspects are informed by studies such as \cite{tang2020reinforcement,wang2023learning}. Other noteworthy contributions to this field include \cite{jia2021benders,turner2023adaptive,berthold2022learning}. For a comprehensive review on this topic, readers may refer to \cite{deze2023machine}.

\subsection{Primal Heuristics}
\label{sec:heur}
In addition to branch-and-bound and cutting plane methods, \emph{primal heuristics} are also an important technique in solving MILP. Although heuristics lack the optimality guarantees that BnB offers, they are typically more computationally efficient and able to promptly reduce the primal bound. Modern solvers usually first resort to heuristic methods before implementing BnB, treating BnB as a final solution. This section will discuss primal heuristics and how machine learning can enhance them. 

\subsubsection{Introduction to heuristics}
Any procedure aimed at swiftly creating a feasible solution $\overline{\vx} \in \sX_{\mathrm{MILP}}$ without promising optimality is termed a \emph{primal heuristic}. The motivation behind primal heuristics is to rapidly reduce the primal bound $\overline{f}=\vc^\top\overline{\vx}$. If a desirable primal bound can be obtained in the early stage of BnB, it is plausible to prune a substantial portion of nodes in the BnB tree (Refer to Lines 6 and 14 in Algorithm \ref{algo:bnb} for pruning). The following are some typical primal heuristics:
\paragraph{Rounding methods.} The aim of rounding methods is to construct a solution that satisfies the integer constraint from a fractional solution, such as an optimal solution to the root LP relaxation $\underline{\vx}$. If $\mathrm{frac}(\underline{\vx}) \neq \emptyset$, some fractional values may be rounded to an integer value. A simple rounding method is as follows:
    \begin{equation}
        \label{eq:heur-rounding}
        x_j^\prime = \begin{cases}
        \lfloor \underline{x}_j \rfloor, \text{if }c_j > 0\\
        \lceil \underline{x}_j \rceil, \text{if }c_j < 0\\
        [\underline{x}_j], \text{if }c_j = 0
    \end{cases} 
    \end{equation}
    where $[x]$ denotes the nearest integer to $x$. This rounding method may generate a solution ${\vx}^\prime$ that violates the linear constraint $\vA\vx\circ\vb$. However, if the generated solution ${\vx}^\prime$ satisfies all the constraints, we consider a feasible solution has been discovered and mark ${\vx}^\prime$ as $\overline{\vx}$ and update $\overline{f}=\vc^\top\overline{\vx}$.
\paragraph{Diving methods.} Essentially, diving methods are \emph{incomplete BnB methods}. In BnB branching, only one child node may be generated instead of two, as prescribed in Algorithm \ref{algo:bnb}. A common diving method is the \emph{fractionality diving}, running BnB with a branching strategy that selects the fractional element with the smallest fractionality $j_\ast = \argmin_{j \in \mathrm{frac}(\underline{\vx})} \max(\underline{x}_j - \lfloor \underline{x}_j \rfloor, \lceil \underline{x}_j \rceil - \underline{x}_j)$ and generates a single child node defined as the following rule until either infeasibility is reached or a feasible solution $\overline{\vx}$ is obtained:
    \begin{equation}
        \label{eq:heur-diving}
        \sX_{\mathrm{child}} = \sX \cap \{ \vx: x_j = [\underline{x}_j] \}
    \end{equation}
    This differs from the typical most fractionality rule that we might adopt in a complete BnB; a diving heuristic does not aim to find the optimal solution and improve the dual bound, but instead seeks to quickly identify a feasible solution by selecting the element with the smallest fractionality.
    Given that such a method resembles diving into the BnB tree, it is aptly termed as a diving method.
\paragraph{Neighborhood search.} Given an estimated solution $\hat{\vx}$, one may explore integer solutions within the neighborhood of $\hat{\vx}$. For instance, a subset of the values in $\vx$ can be fixed to match $\hat{\vx}$, forming a ``sub-MIP" with the remaining unfixed variables. The sub-MIP typically has fewer variables than the original MILP and is hence likely to be solved in less time. Formally, an index set $\sJ$ is determined to create the sub-MIP as follows:
    \begin{equation}
        \label{eq:neighborhood-search}
        \min_{\vx } \vc^\top \vx ~~\text{s.t. }\vx \in \sX_{\mathrm{MILP}}, ~~ x_j = \hat{x}_j \text{ for }j \in \sJ. 
    \end{equation}
    If $\hat{\vx}$ is not feasible, a feasible solution in its neighborhood is sought; if it is already feasible, the goal is to find a feasible solution better than $\hat{\vx}$. If this goal is achieved, the optimal solution to \eqref{eq:neighborhood-search} is designated as $\overline{\vx}$, and $\overline{f}$ is accordingly updated.

Note that many heuristics do not guarantee the generation of a feasible solution within a given time limit. One strategy may be to allocate a time budget for each heuristic and try different heuristics either prior to or during the BnB process. If a feasible solution with a superior primal bound is discovered, BnB will update the incumbent $\overline{\vx}$ and the primal bound $\overline{f}$. Refer to \cite{berthold2006primal} for more practical heuristics.

\subsubsection{Learning to predict a solution}
\label{sec:learn-sol-predictor}
Given that the primary aim of primal heuristics is to discover a feasible solution with as low a primal bound as possible, \emph{it is always beneficial if we can train a neural network to predict the MILP solution.}
Even though the predicted solution might not be fully accurate, it can be immediately employed in rounding methods and neighborhood searches. Moreover, it can support diving methods by suggesting the branching variable. Applications of a machine learning-based solution predictor in primal heuristics will be discussed in the following section. This section will focus on how to train a solution predictor. The training process adheres to the previously mentioned pipeline: parameterization, data collection, model training, and validation. 

\paragraph{Assumptions.} To streamline notation and highlight the principal concept, we presume all the integer variables to be binary: $l_j=0,u_j=1$ for all $j \in \sI$, and our solution predictor only predicts the values of these binary variables. \textbf{This assumption is applicable for Sections \ref{sec:learn-sol-predictor}, \ref{sec:ml-based-heuristics}, and \ref{sec:improving-heur}.}

\paragraph{Parameterization.} Here our goal is to build a parameterized mapping that maps $(\vc,\sX_{\mathrm{MILP}})$ to the predicted solution $\hat{\vx}$. 
The three parameterization approaches, \eqref{eq:SB-imitation-approach1}, \eqref{eq:SB-imitation-approach2} and \eqref{eq:SB-imitation-appraoch3} detailed in Section \ref{sec:l2b}, are all applicable for solution prediction. We will now illustrate the structure of the GNN for solution prediction, bearing in mind that the other two apporaches can be adapted similarly. The solution predictor's input is $(\vc,\sX_{\mathrm{MILP}})$, which can be represented as $\cG = (\vA, \{\vf_j\}_{j=1}^n, \{\vg_i\}_{i=1}^m)$ to suit the GNN, where $\vf_j = (c_j,l_j,u_j,z_j)$ and $\vg_i = (\circ_i,b_i)$. We then apply \eqref{eq:gnn-message-passing} on $\cG$ for multiple rounds to obtain the outcomes $\{\vf_{j,\mathrm{gnn}}\}_{j=1}^n$ and $\{\vg_{i,\mathrm{gnn}}\}_{i=1}^m$. Finally, a binary classifier is applied on $\{\vf_{j,\mathrm{gnn}}\}_{j\in\sI}$ to predict a probability distribution of the binary variable values.
\begin{equation}
    \label{eq:gnn-sol}
    \mathrm{GNN}(j,\cG;\vtheta):=\mathrm{MLP}\left( \vf_{j,\mathrm{gnn}}; \vtheta_7 \right) = p_j,
\end{equation}
where $p_j$ signifies the probability that $\hat{x}_j = 1$, and $1-p_j$ represents the probability that $\hat{x}_j = 0$. Given the above GNN, to derive a predicted solution from the GNN outcome $p_j$, one can either sample $\hat{x}_j$ from the Bernoulli distribution with probability $p_j$ or decide the values directly by setting $\hat{x}_j = 0$ if $p_j \leq 0.5$ and $\hat{x}_j = 1$ otherwise. The parameters of the GNN, denoted by $\vtheta$, include $\vtheta = (\vtheta_1,\vtheta_2,\vtheta_3,\vtheta_4,\vtheta_7)$, where $\vtheta_1,\vtheta_2,\vtheta_3,\vtheta_4$ are as defined in \eqref{eq:gnn-message-passing}.

\paragraph{Data collection.} The process of data collection in this context is straightforward: collecting pairs of MILPs and their corresponding solutions in the format $(\cG,\vx^\ast)$. Given a set of MILP instances for training, one should execute BnB on each instance to obtain its optimal solution. By aggregating these pairs, a dataset denoted as $\sD_{\mathrm{Sol}}$ can be established.

\paragraph{Model training and validation.} This phase follows the standard end-to-end training procedure and employs the binary cross-entropy loss function to determine the parameters in GNN:
\begin{equation}
    \label{eq:gnn-cross-entropy}
    \min_{\vtheta} \sum_{(\cG,\vx^\ast) \in \sD_{\mathrm{Sol}}} ~\sum_{j \in \sI} ~~ - x^\ast_j \log(p_j) - (1-x^\ast_j)\log(1-p_j)
\end{equation}
This type of cross-entropy loss is typically employed in binary classification scenarios. To understand cross-entropy, consider that a perfect fit implies $p_j=1$ for all ${j\in\sI:x^\ast_j=1}$ and $p_j=0$ for all ${j\in\sI:x^\ast_j=0}$. The tuning of hyperparameters and model validation follows the same process outlined in Section \ref{sec:l2b}.

\subsubsection{Incorporating solution predictors into heuristics}
\label{sec:ml-based-heuristics}

Having trained a GNN solution predictor as per \eqref{eq:gnn-sol}, the question arises: How can it be effectively utilized to assist with MILP solving? Using the predicted solution $\hat{\vx}$ directly as $\overline{\vx}$ might entail risks, given that \emph{such a solution predicted by the GNN might be not just suboptimal, but also infeasible.} Nonetheless, incorporating this predicted solution with the three categories of heuristics previously introduced can yield significant benefits. In this section, we outline some methods for merging neural network predictions with primal heuristics.

\paragraph{Rounding.}  A straightforward method to incorporate the GNN in rounding is to directly use the outcome of the GNN, $p_j$, as the $\underline{x}_j$ in \eqref{eq:heur-rounding}. A safer alternative would be to set an acceptance threshold for $p_j$. Let us say the threshold is denoted as $\xi\in(0,1)$, only variables with a probability of being $0$ or $1$ that exceeds $\xi$ are accepted: $x^\prime_j = 1$ if $p_j > \xi$ and $x^\prime_j = 0$ if $1-p_j > \xi$. For all other cases, we use \eqref{eq:heur-rounding} to create values of $\vx^\prime$. Please note that neither of these rounding methods can guarantee feasibility. However, due to their minimal computational cost, it is still worth attempting these methods. If a new feasible solution is discovered, it will be highly beneficial; if not, the rounding results can be discarded.

\paragraph{Diving.} The output of the GNN can potentially enhance a diving method. Given the GNN outputs $\{p_j\}_{j\in\sI}$, one can simply use $p_j$ as the $\underline{x}_j$ in the fractionality diving \eqref{eq:heur-diving}. This leads to the branching rule for diving being defined as: 
\[ j_\ast = \argmin_{j} \max(p_j,1-p_j) \]
The other aspects of diving remain unchanged and follow the fractionality diving approach. If the GNN's output is closer to the optimal solution of the MILP $\vx^\ast$, such an ML-based diving method can potentially yield a better solution compared to the conventional fractionality diving.

\paragraph{Neighborhood search.} Based on the idea of ``fixing some variables and solving the sub-MIP" introduced in \eqref{eq:neighborhood-search}, an ML-enhanced neighborhood search approach can be formulated: variables exceeding a certain acceptance threshold $\xi$ are fixed, with the remaining variables left unfixed, and the resulting sub-MILP is then solved: 
\begin{equation}
    \label{eq:neighborhood-search-gnn}
    \min_{\vx \in \sX_{\mathrm{MILP}}} \vc^\top \vx ~~\text{s.t. } x_j = 1 \text{ for }j \in \{j\in\sI:p_j > \xi\},~ x_j=0 \text{ for }j \in \{j\in\sI:p_j < 1- \xi\}.
\end{equation}
Let $\sJ(\xi) = \{j\in\sI:p_j > \xi\} \cup \{j\in\sI:p_j < 1- \xi\}$ represent the set of variables to be fixed. If $\sJ(\xi)$ is large, the sub-MIP is easier to solve due to its smaller size but the possibility of causing infeasibility increases due to the GNN-output not guaranteeing feasibility. If $\sJ(\xi)$ is small, the opposite scenario is expected. By adjusting $\xi$, the size of $\sJ(\xi)$ can be controlled. However, to ensure feasibility when high-confidence prediction errors are present, one might have to use a very small fixing set and cannot fully leverage the results of the GNN. Figure \ref{fig:local-branching} illustrates an example of this, a typical situation in practice.
\begin{figure}
    \centering
    \input{local-branching.tikz}
    \caption{Determining the threshold in neighborhood search. In this example, we are working with 16 binary variables. For simplicity, we have sorted these variables according to the confidence measure $\max(p_j,1-p_j)$, accepting and fixing only those variables that exhibit sufficient confidence. With the methods in \eqref{eq:neighborhood-search-gnn}, we must set $\xi > 0.8$ to maintain correctness, consequently fixing only the first three variables. However, if we utilize \eqref{eq:local-branching-gnn} with $\tau=2$, thereby allowing a maximum of two elements to be incorrect, we can set $\xi = 0.65$ and fix an additional five variables.}
    \label{fig:local-branching}
\end{figure}
To address this, the constraints in \eqref{eq:neighborhood-search-gnn} may be relaxed, permitting a small error ratio within the set $\sJ(\xi)$: 
\begin{equation}
    \label{eq:local-branching-gnn}
    \min_{\vx \in \sX_{\mathrm{MILP}}} \vc^\top \vx \quad \text{s.t.}  \left( \sum_{j \in \{j\in\sI:p_j > \xi\}} (1-x_j) + \sum_{j \in \{j\in\sI:p_j < 1- \xi\}} x_j \right) \leq \tau.
\end{equation}
In this equation, $\tau$ represents the tolerance for errors in $\sJ(\xi)$. If $\tau=0$, \eqref{eq:local-branching-gnn} becomes \eqref{eq:neighborhood-search-gnn}. If $\tau=1$, the constraint in \eqref{eq:local-branching-gnn} allows for a single element in $\sJ(\xi)$ to differ from the GNN's prediction. This method, known as \emph{local branching}~\cite{fischetti2003local} in the MILP community, offers considerably more flexibility than \eqref{eq:neighborhood-search-gnn}.

\subsubsection{Improvement Heuristics} 
\label{sec:improving-heur}

The methods described earlier aim to produce a feasible solution $\overline{\vx} \in \sX_{\mathrm{MILP}}$ for a MILP from its potentially infeasible estimate $\hat{\vx}$. Additionally, some heuristics can improve an existing feasible solution $\overline{\vx}$, which are usually named \emph{improvement heuristics}. For instance, the technique of local branching, as outlined previously, can be employed to discover an improved feasible solution $\overline{\vx}^\prime$ from a given $\overline{\vx}$ through
\begin{equation}
    \label{eq:local-branching-2}
   \overline{\vx}^\prime = \argmin_{\vx \in \sX_{\mathrm{MILP}}} \vc^\top \vx \quad \text{s.t.}  \left( \sum_{j \in \{j\in\sI:\overline{x}_j = 1\}} (1-x_j) + \sum_{j \in \{j\in\sI:\overline{x}_j = 0\}} x_j \right) \leq \tau.
\end{equation}
This constraint ensures that $\overline{\vx}$ and $\overline{\vx}^\prime$ differ in at most $\tau$ coordinates, effectively limiting the search to a ``neighborhood" around $\overline{\vx}$. 
If $\vc^\top \overline{\vx}^\prime = \vc^\top \overline{\vx}$, it indicates that $\overline{\vx}$ is the optimal solution within this neighborhood. In this case, one may enlarge $\tau$ to explore potentially better feasible solutions within a larger neighborhood. Otherwise, if $\vc^\top \overline{\vx}^\prime < \vc^\top \overline{\vx}$, an improved solution has been found. In such cases, the process repeats with $\overline{\vx}^\prime$ as the new starting point, potentially leading to a sequence of increasingly better solutions $\overline{\vx}^{\prime\prime}, \overline{\vx}^{\prime\prime\prime}, \cdots$. This iterative process is subject to computational constraints, emphasizing the heuristic nature of quickly finding feasible solutions without guarantees of optimality.

While local branching confines the search to a neighborhood, the actual scope of this neighborhood can be vast due to the unknown specific coordinates where $\overline{\vx}$ and $\overline{\vx}^\prime$ differ. Totally there are $\binom{n}{\tau}$ variations, where $n$ is the number of variables. The complexity provides a chance for machine-learning approaches to \textit{predict these specific coordinates, thereby simplifying the MILP by reducing the variable number to $\tau$}, with the remaining variables fixed. In other words, the goal of this approach is to \textit{imitate local branching} through ML models.

A straightforward implementation involves the following steps:
\begin{itemize}
    \item Collecting Data. Iteratively apply local branching to obtain a series of solutions $\overline{\vx}^\prime, \overline{\vx}^{\prime\prime}, \overline{\vx}^{\prime\prime\prime}, \cdots$ until a computational limit is reached or no further improvements are found. Compute difference vectors to indicate changes between successive solutions:
    \[ \vy^\prime = |\overline{\vx} - \overline{\vx}^\prime|, ~~~ \vy^{\prime\prime} = |\overline{\vx}^\prime - \overline{\vx}^{\prime\prime}|, ~~~ \vy^{\prime\prime\prime} = |\overline{\vx}^{\prime\prime} - \overline{\vx}^{\prime\prime\prime}|, ~~~\cdots \]
    Here $|\cdot|$ denotes the coordinate-wise absolute value, and $|\overline{\vx} - \overline{\vx}^\prime| \in \{0,1\}^n$ describes an indicator vector that measures whether $\overline{\vx}$ and $\overline{\vx}^\prime$ are different. 
    We want to train a neural network that maps $(\vc,\sX_{\mathrm{MILP}},\overline{\vx})$ to $\vy^\prime$, and maps $(\vc,\sX_{\mathrm{MILP}},\overline{\vx}^\prime)$ to $\vy^{\prime\prime}$, and maps $(\vc,\sX_{\mathrm{MILP}},\overline{\vx}^{\prime\prime})$ to $\vy^{\prime\prime\prime}$, etc. The data is collected in the form of 
    \[  
    \begin{aligned}
    \cG =& (\vA, \{\vf_j^{\prime}\}_{j=1}^n, \{\vg_i\}_{i=1}^m),~~ \vf_j = (c_j,l_j,u_j,z_j,\overline{x}_j)\\
    \cG^{\prime} =& (\vA, \{\vf_j^{\prime}\}_{j=1}^n, \{\vg_i\}_{i=1}^m) ,~~\vf_j^{\prime} = (c_j,l_j,u_j,z_j,\overline{x}_j^{\prime})\\
    \cG^{\prime\prime} =& (\vA, \{\vf_j^{\prime\prime}\}_{j=1}^n, \{\vg_i\}_{i=1}^m),~~\vf_j^{\prime\prime} = (c_j,l_j,u_j,z_j,\overline{x}_j^{\prime\prime})\\
    & \cdots 
    \end{aligned}
    \]
    Repeat the process for all MILPs in the training set and collect all the data together.
    
    \item Training. The GNN structure follows \eqref{eq:gnn-sol} and training follows \eqref{eq:gnn-cross-entropy}.
    
    \item Inference. After training, the model is then capable of identifying which variables should remain unfixed. Given a new MILP with a feasible solution, one can structure the MILP data into the previously mentioned format $\cG$ and apply the ML model, which outcomes a $\vy$. For each variable $j$, $y_j$ might be a continuous value, with higher values suggesting a greater likelihood that the corresponding variable $x_j$ would change to achieve a better feasible solution $\overline{\vx}^\prime$. Consequently, by ranking the elements of $\vy$ and selecting the top $\tau$ variables, denoted as $\sJ_\tau$, to remain unfixed, and fixing the values of the remaining variables to their current values in $\overline{\vx}$. A sub-MILP is then formulated: 
    \[\min_{\vx \in \sX_{\mathrm{MILP}}} \vc^\top \vx ~~~\mathrm{s.t.}~~ x_j = \overline{x}_j~~\mathrm{for}~j \not\in \sJ_\tau. \] 
    Solving this reduced MILP may lead to a better feasible solution. If no better solution is obtained, one may enlarge the neighborhood size $\tau$ and retry. Repeat the above process until the time budget is reached.
\end{itemize}
This pipeline represents a basic approach, and achieving competitive results may require incorporating advanced techniques from recent research \cite{sonnerat2021learning,wu2021learning,huang2023searching}. We provide some references for other learning heuristics in the next paragraph.

\paragraph{Bibliographical notes.} The success of AlphaGo~\cite{alphago} showed the promising combination of machine learning techniques and exact searching methods. Machine learning excels at swiftly delivering intuitive high-level guidance. Conversely, search methods, like the Monte Carlo tree search utilized in AlphaGo, offer a refined, precise solution grounded on the foresight provided by machine learning. This framework closely mirrors the ML-enhanced heuristics described earlier: the ML model yields a rough solution, while exact mathematical procedures refine and produce an MILP feasible solution based on ML predictions. Within ML4CO (Machine Learning for Combinatorial Optimization), heuristic learning has become one of the mainstream approaches. While the ML-guided rounding explored in this section is straightforward, the inspirations for ML-assisted diving and neighborhood search in this paper are drawn from the works of \cite{shen2021learning} and \cite{huang2022improving,ding2020accelerating} respectively. The following references are highlighted for other aspects on this subject: \cite{nair2020solving,song2020general,hottung2020neural,chmiela2021learning,sonnerat2021learning,wu2021learning,hendel2022adaptive,huang2022anytime,ma2022efficient,falkner2022large,liu2022learning,huang2023searching,paulus2023learning}.

\subsection{Configurations}
\label{sec:config-mip}

Modern MILP solvers integrate the aforementioned BnB, cutting plane, and heuristics, along with other techniques not detailed in this note, like presolving. These solvers provide a collection of hyperparameters to allow users flexibility in controlling specific components. Here, we highlight some of these hyperparameters:
\begin{itemize}
    \item (LP) Which algorithm should be employed for the LP relaxation: simplex or the interior-point method?
    \item (Cut) The budget for the cutting plane method and the order of trying different kinds of cuts.
    \item (Heuristics) The budget for heuristics and the order of trying different heuristic approaches.
    \item (BnB) The chosen branching and searching strategies. It is worth noting that numerous manually-designed strategies are not covered in this note.
\end{itemize}
Let us represent the $k$-th hyperparameter as $h_k$. By consolidating all the hyperparameters into a vector $\vh = (h_1,\cdots,h_K)$, we define a \emph{configuration}. The performance of an MILP solver greatly depends on this configuration.
Typically, an MILP solver offers a default configuration $\vh_{\mathrm{default}}$. Developers derive this configuration from diverse and general-purpose datasets to ensure its reliability. For datasets catering to specific MILP instances, $\vh_{\mathrm{default}}$ might be conservative and not the optimal choice. To obtain the optimal configuration for a specific MILP set $\sM$, one would need to solve:
\begin{equation}
    \label{eq:uniform-config}
  \min_{\vh} \sum_{(\vc,\sX_{\mathrm{MILP}}) \in \sM} \cM\Big( (\vc,\sX_{\mathrm{MILP}}), \vh \Big)
\end{equation}
where $\cM$ measures the solver's performance on a specific MILP instance $(\vc,\sX_{\mathrm{MILP}})$ under a given configuration $\vh$. This performance can be assessed based on time consumption, the BnB tree size, among other metrics. While in this note we only consider the sum over all MILP instances in $\sM$, other possible methods include the maximum (for robust hyperparameter optimization), the geometric mean, etc. Compared with the learning methods discussed in Sections \ref{sec:bnb-ml}, \ref{sec:cut}, and \ref{sec:heur}, the hyperparameter optimization approach treats the MILP solver as a black box, without diving into specific components like heuristics, cuts, and branching. This methodology is also referred to as \emph{black-box optimization} or \emph{algorithm configuration}. 

\paragraph{Hyperparameter tuning.} Solving \eqref{eq:uniform-config} is usually challenging due to several factors: the non-differentiability of $\cM$, its non-convex nature, and the fact that some hyperparameters are categorical rather than continuous. Yet, while pinpointing the global optimum of \eqref{eq:uniform-config} can be challenging, finding reasonably good sub-optimal hyperparameters is tractable, with ample research available in the field. Strategies like random search~\cite{bergstra2012random}, grid search, Bayesian optimization techniques~\cite{snoek2012practical}, and meta-heuristic methods~\cite{gogna2013meta} have been explored. For a comprehensive overview, see \cite{yang2020hyperparameter}. For hyperparameter tuning in the context of MILP challenges, refer to \cite{hutter2009param,hutter2011smac,balcan2018learning,himmich2023mpils,zhang2023mindopt}. This paper will not delve into specifics of any single hyperparameter optimization technique. Instead, we will treat problems in the form of \eqref{eq:uniform-config} as solvable modules and focus on leveraging machine learning techniques to enhance configuration.

\paragraph{Adaptive configuration.} Assigning a uniform configuration to all instances in $\sM$ as per \eqref{eq:uniform-config} lacks flexibility and typically does not yield the optimal outcome. For instance, let us imagine $\sM$ consists of two distinct MILP types. For one type, allocating a substantial computational budget to heuristics might drastically shrink the primal bound and BnB tree size. However, for the other MILP variety, heuristics might prove ineffective, warranting a minimal allocation. Using a strategy like \eqref{eq:uniform-config} would struggle to produce a universally efficient configuration for these diverse MILP types. Rather than settling for a one-size-fits-all solution, adaptive configuration offers a more dynamic alternative: \emph{configuring each instance uniquely based on its needs}. One direct method to achieve this is to tune the hyperparameters for each specific MILP $(\vc,\sX_{\mathrm{MILP}})$ to pinpoint its optimal configuration $\vh_\ast$. Subsequently, a neural network can be trained to map each $(\vc,\sX_{\mathrm{MILP}})$ to its respective $\vh_\ast$ within all instances in $\sM$. However, this method might be too ambitious and rarely adopted in practice. A more practical strategy is partitioning the instances in $\sM$ into several clusters, designating a unique configuration for each cluster, and then training a classification model based on these clusters. Once set up, when we meet a new MILP instance beyond $\sM$, we can invoke the classifier and assign a configuration to this MILP based on its classification result. Here, we delve deeper into this workflow.

\subparagraph{Clustering} The purpose of clustering is to divide the entire set $\sM$ into several subsets, represented by $\sM = \cap_{c=1}^{C} \sM_c$, with each subset referred to as a \emph{cluster}. The goal is to make the instances within each cluster similar enough that a universally effective configuration can be applied to them. The criteria for clustering can be quite diverse and may include:
\begin{itemize}
    \item (Data source). In real-world scenarios, $\sM$ might be gathered from various origins, such as different teams across various businesses. Even though all these problems can be formulated as MILPs, treating instances from distinct sources separately can be more effective.
    \item (Domain knowledge). Utilizing background information about the MILPs in $\sM$ can facilitate clustering. Take the vehicle routing problem as an example; MILPs derived from weekday data and weekend data might be clustered and configured separately, based on such insights.
    \item (Trial and error). A more generic clustering approach is trial-and-error. Imagine a scenario where only one binary hyperparameter $h$ is considered: the LP relaxation algorithm. The value of $h$ represents which method to use, with $0$ for the simplex method and $1$ for the interior-point method. By testing both options across all instances in $\sM$, it is possible to divide $\sM$ into two clusters: $\sM_0$, where the simplex method outperforms the interior-point method, and $\sM_1$, encompassing all remaining instances. For cases involving multiple hyperparameters, this standard trial-and-error approach might be expensive, and practitioners may choose to eliminate certain undesirable configurations based on experience.
\end{itemize}
For each identified cluster $\sM_c$, hyperparameter tuning akin to \eqref{eq:uniform-config} can be employed to derive the optimal configuration $\vh_c$ tied to that cluster. Every MILP instance is tagged with a label to signify its cluster affiliation. If an MILP is part of $\sM_c$, the corresponding label is $\vp_\ast \in \sR^C$, with the $c$-th component set to $1$ and the rest marked as $0$.

\subparagraph{Training} With the defined clusters, our next step is to establish a mapping from a given MILP instance, represented as $(\vc,\sX_{\mathrm{MILP}})$, to its corresponding cluster indicator, $\vp_\ast$. To accomplish this, one can utilize a parameterized model like Approach 2 (Graph-augmented neural networks) or Approach 3 (GNNs) as detailed in Section \ref{sec:l2b}. Here, for demonstration, we use Approach 3. By employing the method outlined in Section \ref{sec:l2b}, data in the form of $(\vc,\sX_{\mathrm{MILP}})$ can be converted into GNN-friendly inputs denoted as $\cG = (\vA, \{\vf_j\}_{j=1}^n, \{\vg_i\}_{i=1}^m)$. Following this, we apply the message passing mechanism \eqref{eq:gnn-message-passing} to $\cG$. The results are labeled as $\{\vf_{j,\mathrm{gnn}}\}_{j=1}^n$ and $\{\vg_{j,\mathrm{gnn}}\}_{j=1}^m$. By aggregating these results, we derive two feature vectors, $\sum_j \vf_{j,\mathrm{gnn}}$ and $\sum_{i} \vg_{j,\mathrm{gnn}}$, that characterize the nature of the MILP instance $(\vc,\sX_{\mathrm{MILP}})$. An MLP is then applied to these vectors. The neural network's output is formatted as a $C$-dimensional vector, depicting the likelihood of the MILP instance's association with a particular cluster. This complete model can be represented as:
\[\ccH_{\mathrm{GNN}}(\cG;\vtheta):= \mathrm{MLP}\left( \sum_{j=1}^n \vf_{j,\mathrm{gnn}}, \sum_{i=1}^m \vg_{j,\mathrm{gnn}}; \vtheta_9 \right) = (p_1,p_2,\cdots,p_C) \in \sR^C\]
In this context, $p_c$ indicates the probability of the MILP aligning with cluster $\sM_c$, and $\vtheta = (\vtheta_1,\vtheta_2,\vtheta_3,\vtheta_4,\vtheta_9)$ embodies all the GNN's parameters, with the definitions of $\vtheta_1,\vtheta_2,\vtheta_3,\vtheta_4$ presented in \eqref{eq:gnn-message-passing}. To effectively train this GNN, we collect data represented as pairs $(\cG,\vp_\ast)$, compiling them into $\sD_{\mathrm{config}}$. Then we fit the GNN's output with $\vp_\ast$. Formally, this is described as:
\[ \min_{\vtheta} \sum_{(\cG,\vp_\ast) \in \sD_{\mathrm{config}}} \ell\big(\vp_\ast, \ccH_{\mathrm{GNN}}(\cG;\vtheta)\big) \]
where $\ell(\vp,\vq)=\sum_{c=1}^C p_c \log(q_c)$ is taken as the cross entropy, a generalized form of the binary cross entropy defined in \eqref{eq:gnn-cross-entropy}.

\subparagraph{Inference} Once the GNN has been trained, it can be utilized to categorize a specific MILP, even if it is outside the training set $\sM$. Assuming the output from the GNN is $\hat{\vp}$, and the classification aligns the MILP to cluster $c_\ast$ (where $c_\ast = \argmax_c \hat{p}_c$), we can then apply the configuration $\vh_{c_\ast}$ to the presented MILP based on this classification.

\paragraph{Bibliographical notes.} The concept of instance-adaptive or instance-aware algorithm configuration is far from new, with initial implementations tracing back to works such as \cite{kadioglu2010isac,xu2011hydramip}, to the best of our knowledge. This direction has gained arising interest in recent years \cite{ansotegui2016maxsat,bonami2022classifier,hosny2023auto}. The methodology and pipeline presented within this paper draw inspirations from two state-of-the-art MIP configuration works \cite{song2023instance,valentin2022instance}. For a comprehensive review, please refer to \cite{malitsky2014isac}.

\subsection{Numerical Results}

In this subsection, we present numerical results of learning to branch, learning to select cuts, learning heuristics, and learning configurations on a common type of MILP problem, namely \textit{Set Covering}. This classic problem in combinatorial optimization aims to find the minimum number of subsets needed to cover a given set of elements. Formally, let \(U = \{1, 2, \ldots, m\}\) denote a set of \(m\) elements, and let \(S = \{S_1, S_2, \ldots, S_n\}\) be a collection of \(n\) subsets of \(U\). The objective is to select the smallest number of these subsets such that every element in \(U\) is included in at least one of the selected subsets. This can be formulated as the following MILP problem:
\[
\min_{\vx \in \{0,1\}^n} ~~~ \sum_{j=1}^{n} x_j, ~~~ \textup{subject to: }\sum_{j: i \in S_j} x_j \geq 1, \quad \forall i \in U,
\]
where \(x_j\) is a binary decision variable that equals 1 if subset \(S_j\) is selected and 0 otherwise. The objective function minimizes the total number of subsets selected, while the constraints ensure that every element in the universe \(U\) is covered by at least one of the chosen subsets.

In our experiments, we randomly generate Set Covering instances with a Python software package named Ecole~\cite{ecole}. The generation scheme follows the convention in~\cite{balas1980set}, which has been recognized and followed by recent work such as~\cite{gasse2019exact}. Set Covering instances for training, validation and testing are generated independently.

Our experiment serves as a numerical showcase on a common ground of the effectiveness that can be expected out of different L2O methodologies for mixed-integer optimization. We select works to report based on consideration of both representativeness and reproducibility.

\paragraph{Environments.} The experiments were conducted on a high-performance workstation equipped with an Intel(R) Xeon(R) Platinum 8163 CPU and 8 NVIDIA Tesla V100 GPUs, utilizing SCIP~8.0.1~\cite{scip8} as the baseline MILP solver and PyTorch 1.8.1 for neural network training. Additionally, Ecole 0.8.1 was deployed to generate problem instances and to facilitate the extraction of bipartite graph representations from MILP instances. This setup also included the use of Ecole's provided demonstrations for learning branching strategies and solver configurations.

\paragraph{ML-based branching, cut selection, and configuration.}

We reproduced numerical results of three existing methods, one for each direction. We implemented \cite{gasse2019exact} as a representative for ML-based branching methods, \cite{wang2023learning} for ML-based cut selection and \cite{valentin2022instance} for ML-based solver configuration. \cite{gasse2019exact,valentin2022instance} are based on Graph Neural Networks (GNNs) that take bipartite representations of MILP instances as inputs. \cite{wang2023learning} exploits reinforcement learning to learn a cut selection policy. Note that here we reduce the configuration space of~\cite{valentin2022instance} to only 16 possible candidates to make the label generation time manageable.\footnote{We set heuristics and separating parameters among [`DEFAULT', `AGGRESSIVE', `FAST', `OFF'], resulting in a configuration space of $4\times4=16$ elements. Other solver settings are set to default.}

All three methods are trained and evaluated on a common dataset of Set Covering problems released by the authors of~\cite{wang2023learning}.\footnote{Please refer to \url{https://github.com/MIRALab-USTC/L2O-HEM-Torch}.} The dataset contains a training partition of 10,000 instances of Set Covering problem, a validation partition of 2,000 instances and a testing partition of 100 instances. During training, we follow the exact hyperparameter settings (e.g., learning rate) used in the original implementation.

Results of the average solving times and primal-dual integrals are reported in Table~\ref{tab:branch-cut-config-experiment}, which clearly show that all three ML-based methods outperform SCIP~8.0.1 with default configuration by significant margins in both metrics. Moreover,
\begin{itemize}
    \item While the ML-based branching method~\cite{gasse2019exact} achieves the best performance, it requires an expensive process of collecting high-quality branching examples to facilitate supervised learning to imitate the full strong branching policy. This characteristic limits its scalability to larger problems as the complexity of example collection becomes prohibitive.
    \item The ML-based configuration method~\cite{valentin2022instance} runs a similar process to collect training labels for each instance-configuration pairs (primal-dual integrals within a fixed time limit in its case). The time cost grows linearly with respect to the size of configuration space and exponentially with respect to the number of solver parameters considered. Therefore, it is essential to exclude non-promising configurations early, which is one of main contributions of~\cite{valentin2022instance}.
    \item In contrast, the RL-based cut selection method~\cite{wang2023learning} is free of such data collection process and thus is more scalable to larger-scale mixed-integer programming problems.
\end{itemize}

\begin{table}[t]
\centering
\caption{Results of ML-based branching, cut selection and configuration on Set Covering problems with 500 constraints and 1,000 variables. Numbers reported here are the average on 100 testing instances.}
\label{tab:branch-cut-config-experiment}
\begin{tabular}{lcccc}
\toprule
 \multirow{2}{*}{Method} & \multicolumn{2}{c}{Time (seconds)} & \multicolumn{2}{c}{Primal Dual Integral} \\
                         \cmidrule(l){2-3} \cmidrule(l){4-5}
    &   Value   & Improvement  &  Value  & Improvement \\ \hline
SCIP Default~\cite{scip8}
    &  11.82    &  -           & 103.13  &  -          \\ \hline
ML-based Branching~\cite{gasse2019exact}
    &  2.03     &  +82.83\%    &  24.76  &  +75.99\%   \\ \hline
ML-based Configuration~\cite{valentin2022instance}
    &  4.24     &  +64.13\%    &  44.18  &  +57.16\%   \\ \hline
ML-based Cut Selection~\cite{wang2023learning}
    &  3.12     &  +73.60\%    &  59.66  &  +42.15\%   \\
\bottomrule
\end{tabular}
\end{table}

\paragraph{Learning heuristics.}
For the heuristic learning component, we incorporated methodologies from the research presented in \cite{huang2023searching}, available at \url{https://github.com/facebookresearch/CL-LNS.git}. This approach is grounded in the strategies discussed in Section \ref{sec:improving-heur}, but it further integrates sophisticated ML models such as Graph Attention Networks (GAT) along with an enhanced training protocol known as contrastive learning. In this experiment, we leveraged Ecole to create two distinct sets of problem instances: one set involving smaller problems with dimensions $n=1000,m=500$, and another set consisting of larger problems with dimensions $n=4000,m=5000$. Each set involved $1000$ MILP instances for training and validation, alongside additional $100$ instances for testing.

In contrast to the evaluation metrics applied to ML-based branching, cut selection, and configuration—where the emphasis is on the duration taken to achieve a certain primal-dual gap—our heuristic performance assessment employs a different criterion. We evaluate the effectiveness of heuristics based on \textit{the quality of the best solutions they can generate within a specified time limit}. This alternative metric is commonly used in relevant literature.\footnote{The rationale behind this metric lies in the nature of heuristics as somewhat independent components within an MILP solver. Evaluating heuristics based on the optimal objective achievable within a fixed time frame aligns with their fundamental purpose: to quickly generate feasible solutions with potentially lower objective values. Conversely, branching strategies, cut selection, and solver configuration are deeply tied with specific solvers, making it impractical to isolate and evaluate these modules as standalone entities. Consequently, their performance is inherently linked to the overall execution time of the solver where they are integrated.} Specifically, we employ the following two metrics to measure the heuristics's performance:
\begin{itemize}
    \item The first metric is a curve that tracks the best achievable objective value over time, $\vc^\top \overline{\vx}(t)$, with $\overline{\vx}(t)$ representing the optimal feasible solution at time $t$. This value is computed as an average across all $100$ test instances.
    \item The second metric assesses the comparative performance at a specific time $t$ across the test instances, detailing the number of instances where the ML-based heuristics outperform, the instances where SCIP's provided heuristics win, and the instances where both approaches result in a tie.
\end{itemize}
While the first metric provides a broad overview of average performance, the second offers a deeper analysis, examining outcomes in an instance-by-instance manner. In our analysis, we benchmark our findings against the heuristics provided by SCIP. It is important to note that SCIP employs various manually crafted heuristic strategies rather than relying on a single method. However, using SCIP's default configuration is unfair here. This is because SCIP's default mode is not solely focused on primal heuristics; it also dedicates significant computational effort to branching and cutting processes. Therefore, we use the aggressive heuristic mode in SCIP, encouraging SCIP to allocate more resources towards heuristic processes and minimizing the focus on other computations. 

\begin{figure}
    \centering
    \begin{minipage}[b]{0.49\textwidth}
        \includegraphics[width=\textwidth]{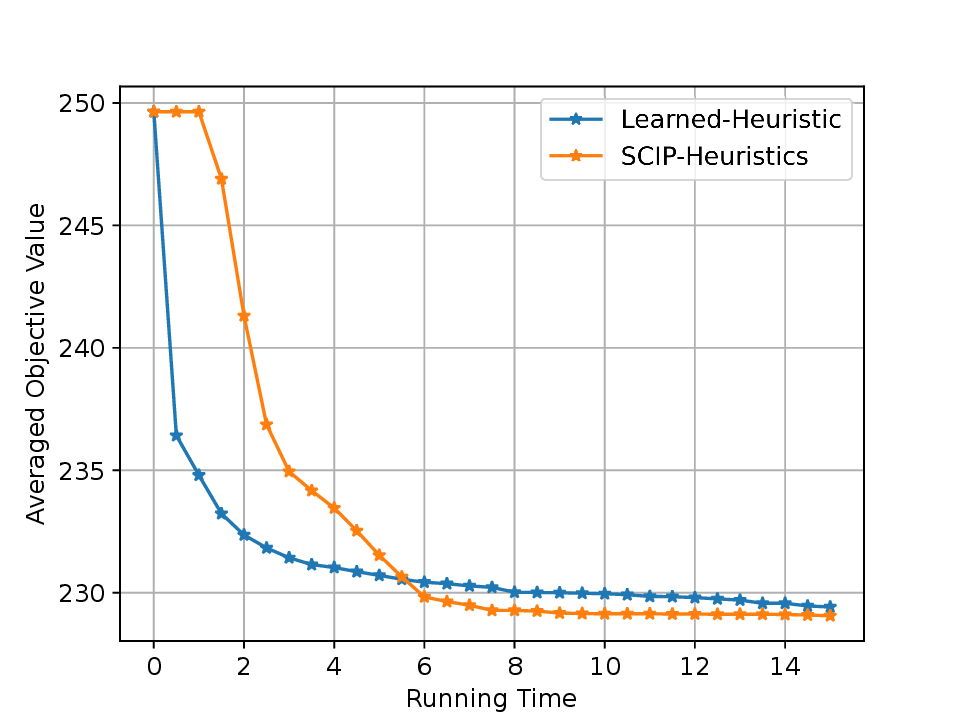}
        \caption{Set Covering problems with size $500 \times 1000$.}
        \label{fig:heuristic-small}
    \end{minipage}
    \hfill 
    \begin{minipage}[b]{0.49\textwidth}
        \includegraphics[width=\textwidth]{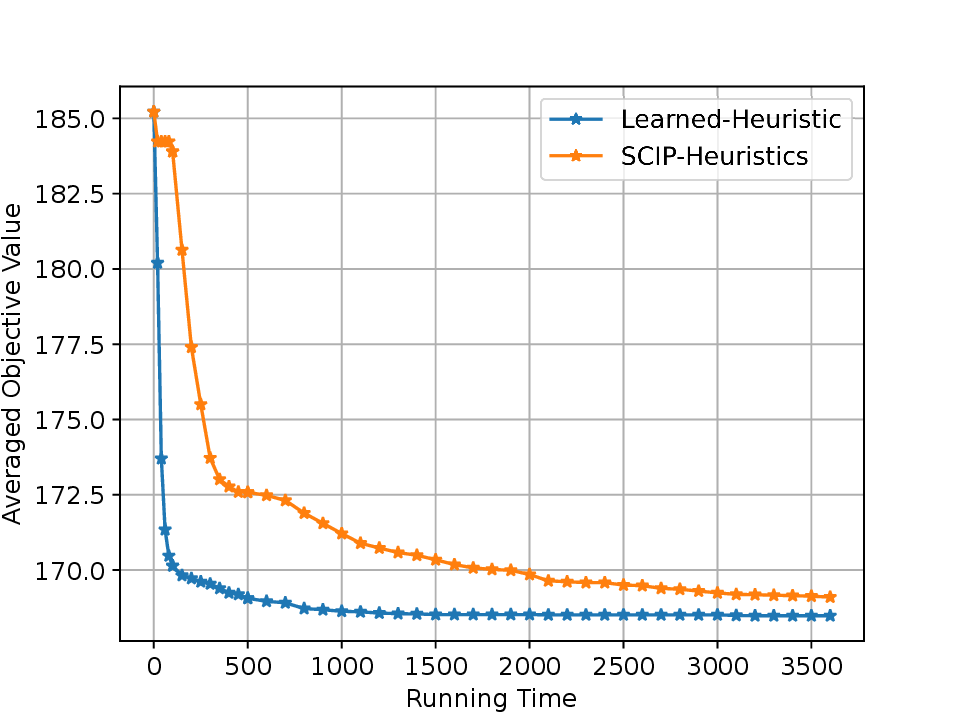}
        \caption{Set Covering problems with size $5000 \times 4000$.}
        \label{fig:heuristic-large}
    \end{minipage}
\end{figure}

\begin{table}[t]
\centering
\caption{Comparative Victory Count Across Test Instances for Various Heuristic Approaches.}
\label{tab:heuristics-experiment}
\begin{tabular}{c|c|c|c|c|c|c}
\hline\hline
\multicolumn{7}{c}{Performance on small-scale set-covering problems (size: $500 \times 1000$)}\\ \hline
Running time (secs) & 1  & 2  & 4  & 6  & 8  & 12 \\ \hline
Learned-Heuristic wins & 98 & 72 & 24 & 5  & 4  & 4  \\ \hline
SCIP-Heuristics win & 0  & 17 & 36 & 38 & 33 & 28 \\ \hline
Tied & 2  & 11 & 40 & 57 & 63 & 68 \\ \hline\hline
\multicolumn{7}{c}{Performance on large-scale set-covering problems (size: $5000 \times 4000$)}\\ \hline
Running time (secs) & 100 & 200 & 400 & 800 & 1600 & 3200 \\ \hline
Learned-Heuristic wins & 99  & 93  & 88  & 79  & 67   & 50   \\ \hline
SCIP-Heuristics win & 0   & 3   & 5   & 8   & 17   & 26   \\ \hline
Tied & 1   & 4   & 7   & 13  & 16   & 24   \\ \hline\hline
\end{tabular}
\end{table}

The results using the first evaluation metric can be seen in Figures \ref{fig:heuristic-small} and \ref{fig:heuristic-large}, while those corresponding to the second metric are detailed in Table \ref{tab:heuristics-experiment}. It is worth mentioning that our implementation involved a scaled-down version of the GAT network and a reduced batch size during training, diverging from the configurations recommended in \cite{huang2023searching}, due to the memory limit of our hardware.\footnote{Specifically, we set the embedding size in GAT as $32$ rather than $64$, and set the batch size as $8$ rather than $32$.} The numerical results clearly show that: \textit{Initially, the ML-based heuristic demonstrates an overwhelming advantage; however, as the computation progresses, this advantage diminishes and might even fall behind the heuristics provided by SCIP.} This phenomenon can be attributed to the \textit{intuitive nature} of machine learning models in contrast to the precision of mathematical algorithms employed by SCIP. Consequently, leveraging ML techniques to swiftly derive a feasible solution is promising, particularly in outperforming SCIP in terms of quickly finding superior feasible solutions within a \textit{limited timeframe}.

\subsection{Summaries}

\begin{figure}
    \centering
    \input{ml4mip-summary.tikz}
    \caption{A recipe of machine learning in solving MILP. Here each ML module interacts with a certain component in the MILP solver: obtaining the training data from the solver and returning actions to guide the solver. Note that this is a simplified diagram. In modern MILP solvers, some heuristics might be employed before the root LP relaxation or following the cutting plane stage, and specific cuts might be considered before heuristics or during the BnB phase.}
    \label{fig:ml4mip-summary}
\end{figure}
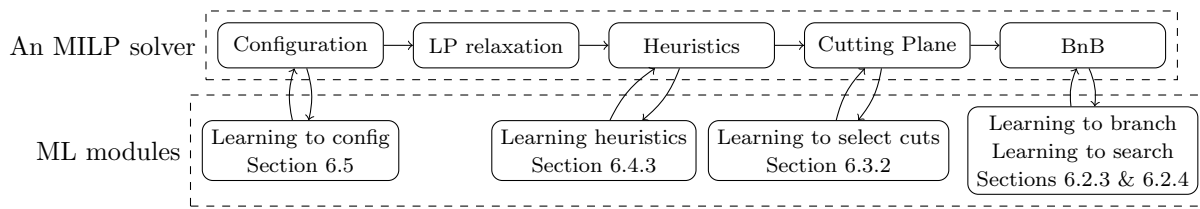

The above-introduced techniques may leave readers questioning how these methods are integrated within an MILP solver, and what connections exist between the ML-based approaches and the MILP solver itself. Figure \ref{fig:ml4mip-summary} offers an illustrative answer. Here, we highlight the motivations of the aforementioned ML-based modules:
\begin{itemize}
    \item (Learning to branch). Branching rules have a significant impact on the size of the BnB tree. While the SB rule is effective in controlling the BnB tree, it can be computationally expensive. Training a neural network to quickly approximate SB has proven effective, and recent research suggests that novel branching rules beyond SB can be discovered.
    \item (Learning to search). Quick discovery of BnB nodes containing the optimal solution may lead to better feasible solutions and a substantial reduction in the primal bound. Neural networks can be trained to predict whether a BnB node contains the optimal solution, thereby guiding the order in which nodes are processed.
    \item (Learning to select cuts). Properly adding cuts to the LP relaxation can greatly enhance the dual bound. However, indiscriminately adding all possible cuts can lead to computational challenges. Neural networks can be used to predict the impact of potential cuts on the dual-bound, guiding the selection of priority cut candidates.
    \item (Learning heuristics). Before initiating BnB, swiftly finding a feasible solution with a low objective value can reduce both the primal bound and the overall BnB tree. Neural networks can help predict the solution values for MILPs. While the predicted solutions might be infeasible, they can still serve as a foundation for constructing feasible solutions.
    \item (Learning to config). Since MILP instances may not all require the same optimal configuration, clustering them into distinct categories allows for specialized hyperparameter tuning. Neural networks can predict the category an MILP instance fits into, facilitating the assignment of adaptive configurations.
\end{itemize}
While this overview may inspire the start of hands-on projects, questions about the selection of the right techniques or theoretical foundations might still arise. Others may wish to delve into the open challenges in this field. In the following sections, we will provide insights into these aspects.

\paragraph{Practical questions.} This paragraph aims to guide the reader on how to choose the appropriate method from the many techniques detailed previously. Below are specific questions and answers to assist in this decision-making process:
\begin{itemize}
    \item \textbf{Which solver component should be improved using machine learning?} The answer depends on the specific MILP instances. We recommend employing a ``\emph{ground truth observation}" strategy. For instance, when considering heuristics, to identify if the MILPs of interest are suitable for learning heuristics, one might first solve an MILP instance to retrieve its optimal solution. If we hypothetically had this optimal solution (ground truth) prior to running the solver, it can be used as the predicted solution $\hat{\vx}$ in the heuristics detailed in Section \ref{sec:ml-based-heuristics}. The aim is to see if this optimal-solution-driven heuristic can reduce the solving time. Here are potential outcomes:
    \begin{itemize}
        \item If the solving time is not notably shortened using these heuristics, it indicates learning heuristics for this instance is useless because a learned model cannot surpass the performance of the optimal solution.
        \item On the other hand, if there's a noticeable benefit, it demonstrates the possibility of learning heuristics, assuming the neural network can predict the solution accurately. 
    \end{itemize}
    This same approach can be applied to other modules. To check the viability of learning to config, one can tune hyperparameters for each MILP instance separately and observe the effect of this optimal configuration (the ground truth). In learning to branch, one can use SB score as the ground truth; in learning to search, one can use \eqref{eq:node-score-target}; in learning to select cuts, one can use \eqref{eq:look-ahead-cut}. Using this ground truth observation method allows for quick decisions on which components should be boosted using ML. If the ground truth does not significantly improve the solver, pursuing that approach must be useless. Otherwise, it merits consideration.
    \item \textbf{What sequence should one follow when trying different ML modules?} We suggest the following sequence: start with learning to configure, then move to learning heuristics, followed by learning to select cuts, and finally, delve into learning to branch and search in BnB. The reason for prioritizing configuration is its minimal need for interaction with the solver and its relative independence. If a simpler method is already efficient, it negates the need for more complicated ones. The sequence of other modules is then determined by their dependency on solvers and frequency of interaction.
    \item \textbf{Among the three parameterization methods in Section \ref{sec:l2b}, which should one opt for?} We recommend the following sequence: Approach 1, followed by Approach 2, and lastly, Approach 3. The reasoning aligns with previous points: starting with the simplest to implement and progressing to the most complex.
    \item \textbf{Should one utilize supervised learning or reinforcement learning?} Supervised learning is suggested to begin with, owing to its straightforward training process. Moreover, recent studies \cite{alphago,qu2022improved} indicate that expert strategies in supervised learning can initialize the policies in reinforcement learning. Reinforcement learning becomes an option when: 
    (1) Supervised learning has been executed, but there's a quest for an even superior model. 
    (2) Determining the supervised learning ground truth becomes intractable, necessitating reinforcement learning. For instance, if certain MILP instances have more than $10^5$ variables, then acquiring full strong branching at a single BnB node can demand thousands of LP evaluations. Clearly, it is impractical to employ full strong branching to determine the supervised learning ground truth.
\end{itemize}

\paragraph{Theoretical questions.} In this paragraph, we introduce theoretical aspects concerning the application of ML to MILP, while highlighting recent literature contributions. Common theoretical concerns in general machine learning contexts include the expressive capacity of specific ML models (e.g., GNNs), the methodology behind training ML models and its theoretical convergence assurance, and how trained models generalize to instances absent from the training dataset. Let us discuss these within the context of MILP tasks.
\begin{itemize}
    \item Expressivity. The expressive power measures the ability of given ML models to fit certain mappings. For example, within the framework of learning heuristics (refer to Sections \ref{sec:learn-sol-predictor} and \ref{sec:ml-based-heuristics}), the GNN's expressiveness can be quantified using the metric:
    \begin{equation}
        \label{eq:expressive-gnn}
        \ccL(\vtheta, \sD_{\mathrm{Sol}}) := \sup_{(\cG,\vx^\ast) \in \sD_{\mathrm{Sol}}} \sum_{j\in\sI} |x^\ast_j - \lfloor p_j \rceil |,~~\text{where }p_j = \mathrm{GNN}(j,\cG;\vtheta)
    \end{equation}
    Here, $\mathrm{GNN}$ refers to the GNN detailed in \eqref{eq:gnn-sol}, and $\lfloor \cdot \rceil$ indicates nearest neighbor rounding. If $\inf_{\vtheta} \ccL(\vtheta, \sD_{\mathrm{Sol}}) = 0$, it implies that this class of GNNs has \emph{strong expressive power}.  
    For any $\varepsilon > 0$, a corresponding $\vtheta$ exists such that the GNN's predicted solution is proximate to the optimal one, within a distance of $\varepsilon$.  
    Otherwise, employing GNNs to approximate MILP solutions inherently introduces errors. In such cases, GNN may not be universally adept as a solution approximator for every instance in $\sD_{\mathrm{Sol}}$. 
    Existing research, such as \cite{chen2022lp,chen2022milp,chen2024rethinking}, has made strides in understanding the expressiveness of GNNs concerning MILPs. A key takeaway is that while standard GNNs are universally solution approximators for all LPs, they falter for all MILPs. Instances of MILP with pronounced symmetry may exhibit inherent representational inaccuracies, prompting practitioners to break this symmetry using methods like random features.
    \item Training. Possessing robust expressiveness merely confirms the existence of a $\vtheta$ that minimizes error. This arises a question: Are there methods to calculate such a $\vtheta$ within a time limit? While algorithms, such as Adam \cite{adam}, can address \eqref{eq:expressive-gnn} and often yield a satisfactory $\vtheta$ in practice, the theoretical foundation to ensure the derivation of a $\vtheta$ so that $\ccL(\vtheta, \sD_{\mathrm{Sol}}) < \varepsilon$ within a complexity limit remains elusive.
    \item Generalization. Even if one successfully identifies $\vtheta_\ast = \argmin_{\vtheta}\ccL(\vtheta, \sD_{\mathrm{Sol}})$, there still exists a significant question: Is $\vtheta_\ast$ applicable to an independent dataset $\sD_{\mathrm{Sol}}^\prime$?  Unfortunately, to the best of our knowledge, such theoretical results still lack, although the generalization ability of ML models is numerically tested in almost all papers focusing on machine learning for MILP.
\end{itemize}
Although we define the above concepts in the context of learning heuristics, they can be naturally extended to other ML modules presented in this paper.

\section{Conclusions}

This tutorial has provided an in-depth exploration of L2O, a promising frontier that integrates the strengths of machine learning with traditional optimization techniques. By examining numerous optimization cases, we have demonstrated that L2O has the potential to reshape the way optimization problems are approached and solved dramatically. Whether it is about accelerating established algorithms, directly producing solutions, or even reconfiguring the optimization problem itself, L2O offers a range of methodologies to cater to various application scenarios. While L2O thrives in environments with ample past experience, its application must be carefully selected, keeping in mind the nature of the optimization problems and the availability of representative data. 
Although most approaches introduced in this paper remain at a research stage and are not yet ready to be used in commercial solvers, we firmly believe that the incorporation between machine learning and optimization has just begun and will continue to evolve.






\end{document}

%% file: math_commands.tex
\usepackage{amsmath,amsfonts,bm}
\usepackage{amsthm,amssymb}

\def\1{\bm{1}}

\def\rvepsilon{{\mathbf{\epsilon}}}

\def\vzero{{\bm{0}}}
\def\vone{{\bm{1}}}

\def\vtheta{{\bm{\theta}}}
\def\vTheta{{\bm{\Theta}}}
\def\va{{\bm{a}}}
\def\vb{{\bm{b}}}
\def\vc{{\bm{c}}}
\def\vd{{\bm{d}}}
\def\ve{{\bm{e}}}
\def\vf{{\bm{f}}}
\def\vg{{\bm{g}}}
\def\vh{{\bm{h}}}

\def\vl{{\bm{l}}}

\def\vn{{\bm{n}}}

\def\vp{{\bm{p}}}
\def\vq{{\bm{q}}}

\def\vs{{\bm{s}}}

\def\vu{{\bm{u}}}
\def\vv{{\bm{v}}}
\def\vw{{\bm{w}}}
\def\vx{{\bm{x}}}
\def\vy{{\bm{y}}}
\def\vz{{\bm{z}}}

\def\mA{{\bm{A}}}

\def\mG{{\bm{G}}}

\def\mI{{\bm{I}}}

\def\mW{{\bm{W}}}

\DeclareMathAlphabet{\mathsfit}{\encodingdefault}{\sfdefault}{m}{sl}
\SetMathAlphabet{\mathsfit}{bold}{\encodingdefault}{\sfdefault}{bx}{n}

\def\gD{{\mathcal{D}}}

\def\sA{{\mathbb{A}}}

\def\sC{{\mathbb{C}}}
\def\sD{{\mathbb{D}}}

\def\sI{{\mathbb{I}}}
\def\sJ{{\mathbb{J}}}

\def\sM{{\mathbb{M}}}

\def\sQ{{\mathbb{Q}}}
\def\sR{{\mathbb{R}}}
\def\sS{{\mathbb{S}}}

\def\sV{{\mathbb{V}}}

\def\sX{{\mathbb{X}}}

\def\sZ{{\mathbb{Z}}}

\DeclareMathOperator*{\argmax}{arg\,max}
\DeclareMathOperator*{\argmin}{arg\,min}

\usepackage{amssymb}
\usepackage{amsbsy}
\usepackage{amsmath}
\usepackage{amsthm}
\usepackage{graphicx}
\usepackage{color}
\usepackage{multirow}
\usepackage{booktabs}
\usepackage{wasysym}
\usepackage{mathtools} 

\newcommand{\vA}{{\mathbf{A}}}

\newcommand{\vD}{{\mathbf{D}}}

\newcommand{\vI}{{\mathbf{I}}}

\newcommand{\vS}{{\mathbf{S}}}

\newcommand{\vW}{{\mathbf{W}}}

\newcommand{\vY}{{\mathbf{Y}}}

\newcommand{\valpha}{{\boldsymbol \alpha}}
\newcommand{\vbeta}{{\boldsymbol \beta }}

\newcommand{\vlambda}{{\boldsymbol \lambda}}

\newcommand{\vnu}{{\boldsymbol \nu}}

\newcommand{\vtau}{{\boldsymbol \tau}}

\newcommand{\cC}{{\mathcal{C}}}
\newcommand{\cE}{{\mathcal{E}}}

\newcommand{\cG}{{\mathcal{G}}}
\newcommand{\ccH}{{\mathcal{H}}} 

\newcommand{\ccL}{{\mathcal{L}}} 
\newcommand{\cM}{{\mathcal{M}}}
\newcommand{\cN}{{\mathcal{N}}}
\newcommand{\cO}{{\mathcal{O}}}

\newcommand{\cS}{{\mathcal{S}}}
\newcommand{\cT}{{\mathcal{T}}}

\newcommand{\cV}{{\mathcal{V}}}

\newcommand{\EE}{\mathbb{E}}
\newcommand{\RR}{\mathbb{R}}
\newcommand{\CC}{\mathbb{C}}
\newcommand{\ZZ}{\mathbb{Z}}
\newcommand{\QQ}{\mathbb{Q}}

\newcommand{\grad}{{\nabla}}    

\DeclareMathOperator*{\Min}{minimize}

\newcommand{\prox}{{\mathrm{prox}}}

\newcommand{\proj}{{\mathrm{proj}}}

\newcommand{\St}{{\quad\mathrm{subject~to}}}



%% file: unrolling-general.tikz
\begin{tikzpicture}
\definecolor{washedblue}{rgb}{0.85, 0.85, 1}
\tikzstyle{op}=[circle, draw, fill=washedblue, minimum size=1.8cm, align=center];
\tikzstyle{inter}=[circle, minimum size=1.8cm, align=center];

\node (input) {$\vd$};
\node[below = 2cm of input] (x0) {$\vx^{(0)}$};

\node[op, right = 1.5cm of x0] (layer1) {$\vg(\cdot; \vtheta^{(1)})$};
\node[inter, right = 2cm of layer1] (dots) {$\cdots\cdots$};
\node[op, right = 2cm of dots] (layerK) {$\vg(\cdot; \vtheta^{(T)})$};

\node[right=1.5cm of layerK] (output) {$\vx^{(T)}$};

\draw [>=latex,->] (input) -| (layer1);
\draw [>=latex,->] (input) -| (dots);
\draw [>=latex,->] (input) -| (layerK);

\draw [>=latex,->] (x0) edge (layer1);
\draw [>=latex,->] (layer1) edge node[auto] {$\vx^{(1)}$} (dots);
\draw [>=latex,->] (dots) edge node[auto] {$\vx^{(T-1)}$} (layerK);
\draw [>=latex,->] (layerK) edge (output);


  

\end{tikzpicture}

%% file: var-select-1.tikz
\begin{tikzpicture}[
level 1/.style={sibling distance=3.7cm},
  level 2/.style={sibling distance=1.85cm},
  every node/.style = {shape=rectangle, rounded corners,
    draw, align=center, font=\footnotesize,
    , minimum height=0.8cm
    }]]
  \node {$(0.95,1,0.55,0)$\\$29.6$}
    child { node {$(1,0.9,0.6,0)$\\$29.6$} 
        child { node {$(1,1,0.55,0)$\\$30.15$} 
            child { node {$(1,1,1,0)$\\$36$,~Incumbent} 
            edge from parent node[left, draw=none, fill=none] {$x_3 \geq 1$}
            }
            child { node {Infeasible} 
            edge from parent node[right, draw=none, fill=none] {$x_3 \leq 0$}
            }
            edge from parent node[left, draw=none, fill=none] {$x_2 \geq 1$}
        } 
        child { 
            node {Infeasible} 
            edge from parent node[right, draw=none, fill=none] {$x_2 \leq 0$}
        } 
        edge from parent node[left, draw=none, fill=none] {$x_1 \geq 1$}
    }
    child { node {$(0,1,1,0.5)$\\$32$} 
        child { node {$(0,0.9,0.6,1)$\\$32.6$} 
            child { node {$(0,1,0.55,1)$\\$33.15$} 
                child { node {$(0,1,1,1)$\\$39$} 
                    edge from parent node[left, draw=none, fill=none] {$x_3 \geq 1$}
                }
                child { node {Infeasible} 
                    edge from parent node[right, draw=none, fill=none] {$x_3 \leq 0$}
                }
                edge from parent node[left, draw=none, fill=none] {$x_2 \geq 1$}
            }
            child { 
                node {Infeasible} 
                edge from parent node[right, draw=none, fill=none] {$x_2 \leq 0$}
                }
            edge from parent node[left, draw=none, fill=none] {$x_4 \geq 1$}
        } 
        child { 
            node {Infeasible} 
            edge from parent node[right, draw=none, fill=none] {$x_4 \leq 0$}
            } 
        edge from parent node[right, draw=none, fill=none] {$x_1 \leq 0$}
    };
\end{tikzpicture}

%% file: var-select-2.tikz
\begin{tikzpicture}[
level 1/.style={sibling distance=3.7cm},
  level 2/.style={sibling distance=1.85cm},
   level 3/.style={sibling distance=2cm},
  every node/.style = {shape=rectangle, rounded corners,
    draw, align=center, font=\footnotesize,
    , minimum height=0.8cm
    }]]
  \node {$(0.95,1,0.55,0)$\\$29.6$}
    child { node {$(1,0.5,1,0)$\\$30$} 
        child { node {$(0.5,1,1,0)$\\$30.5$} 
            child { node {$(1,1,1,0)$\\$36$,~Incumbent} 
                edge from parent node[left, draw=none, fill=none] {$x_1 \geq 1$}
            }
            child { node {$(0,1,1,0.5)$\\$32$} 
                child { node {$(0,1,1,1)$\\$39$} 
                    edge from parent node[left, draw=none, fill=none] {$x_4 \geq 1$}}
                child { node {Infeasible} 
                    edge from parent node[right, draw=none, fill=none] {$x_4 \leq 0$}
                }
                edge from parent node[right, draw=none, fill=none] {$x_1 \leq 0$}
            }
            edge from parent node[left, draw=none, fill=none] {$x_2 \geq 1$}
        } 
        child { node {Infeasible} 
            edge from parent node[right, draw=none, fill=none] {$x_2 \leq 0$}
        } 
        edge from parent node[left, draw=none, fill=none] {$x_3 \geq 1$}
    }
    child { node {Infeasible} 
        edge from parent node[right, draw=none, fill=none] {$x_3 \leq 0$}
    };
\end{tikzpicture}

%% file: node-select-1.tikz
\begin{tikzpicture}[
level 1/.style={sibling distance=3cm},
  level 2/.style={sibling distance=2cm},
   level 3/.style={sibling distance=2.5cm},
    level 4/.style={sibling distance=2.1cm},
  every node/.style = {shape=rectangle, rounded corners,
    draw, align=center, font=\footnotesize,
    , minimum height=0.8cm
    }]]
  \node {$(0.95,1,0.55,0)$\\$4.6$,~\textcircled{0}}
    child { node {$(1,0.5,1,0)$\\$5$,~\textcircled{2}} 
        child { 
            node {$(0.5,1,1,0)$\\$5.5$,~\textcircled{4}}
            child { node {$(1,1,1,0)$\\$6$,~\textcircled{8},~Incumbent}
                edge from parent node[left, draw=none, fill=none] {$x_1 \geq 1$}
            }
            child { node {$(0,1,1,0.5)$\\$7$,~\textcircled{5}}
                child { node {$(0,1,1,1)$\\$9$,~\textcircled{7},~Incumbent}
                    edge from parent node[left, draw=none, fill=none] {$x_4 \geq 1$}
                }
                child { node {Infeasible\\~\textcircled{6}}
                    edge from parent node[right, draw=none, fill=none] {$x_4 \leq 0$}
                }
                edge from parent node[right, draw=none, fill=none] {$x_1 \leq 0$}
            }
            edge from parent node[left, draw=none, fill=none] {$x_2 \geq 1$}
        }
        child { node {Infeasible\\~\textcircled{3}} 
            edge from parent node[right, draw=none, fill=none] {$x_2 \leq 0$}
        }
        edge from parent node[left, draw=none, fill=none] {$x_3 \geq 1$}
    }
    child { node {Infeasible\\~\textcircled{1}} 
        edge from parent node[right, draw=none, fill=none] {$x_3 \leq 0$}
    };
\end{tikzpicture}

%% file: node-select-2.tikz
\begin{tikzpicture}[grow=down,
level 1/.style={sibling distance=3cm},
  level 2/.style={sibling distance=2cm},
  level 3/.style={sibling distance=2.5cm},
    level 4/.style={sibling distance=3cm},
  every node/.style = {shape=rectangle, rounded corners,
    draw, align=center, font=\footnotesize,
    , minimum height=0.8cm
    }]]
  \node {$(0.95,1,0.55,0)$\\$4.6$,~\textcircled{0}}
    child { node {$(1,0.5,1,0)$\\$5$,~\textcircled{2}} 
        child { 
            node {$(0.5,1,1,0)$\\$5.5$,~\textcircled{4}}
            child { node {$(1,1,1,0)$\\$6$,~\textcircled{6},~Incumbent}
                edge from parent node[left, draw=none, fill=none] {$x_1 \geq 1$}
            }
            child { node {$(0,1,1,0.5)$\\$7$,~\textcircled{5}}
                child { node{Pruned due to $7 > 6$}
                    edge from parent node[left, draw=none, fill=none] {$x_4 \geq 1$}
                } 
                child { node{Pruned due to $7 > 6$}
                    edge from parent node[right, draw=none, fill=none] {$x_4 \leq 0$}
                } 
                edge from parent node[right, draw=none, fill=none] {$x_1 \leq 0$}
            }
            edge from parent node[left, draw=none, fill=none] {$x_2 \geq 1$}
        }
        child { node {Infeasible\\~\textcircled{3}} 
            edge from parent node[right, draw=none, fill=none] {$x_2 \leq 0$}
        }
        edge from parent node[left, draw=none, fill=none] {$x_3 \geq 1$}
    }
    child { node {Infeasible\\~\textcircled{1}} 
        edge from parent node[right, draw=none, fill=none] {$x_3 \leq 0$}
    };
\end{tikzpicture}

%% file: bipartite.tikz
\begin{tikzpicture}[
    every node/.style={circle,draw},
    edge from parent/.style={->,draw},
    >=latex]
\node (A1) {$\vf_{1}$};
\node[below = 3mm of A1] (A2) {$\vf_{2}$};
\node[below = 3mm of A2] (A3) {$\vf_{3}$};

\node[right = 1cm of A1] (B1) {$\vg_{1}$};
\node[right = 1cm of A3] (B2) {$\vg_{2}$};
  
\draw (A1) -- node[pos=0.5,draw=none] {1} (B1); 
\draw (A2) -- node[pos=0.5,draw=none] {1} (B1);
\draw (A2) -- node[pos=0.5,draw=none] {1} (B2);
\draw (A3) -- node[pos=0.5,draw=none] {2} (B2);

\node[draw, dashed, fit=(A1) (A3), inner sep=1mm, rectangle, label=above:$\sV$] (V) {}; 
\node[draw, dashed, fit=(B1) (B2), inner sep=1mm, rectangle, label=above:$\sC$] (C) {}; 
\end{tikzpicture}

%% file: sepa.tikz
\begin{tikzpicture}
    \tikzstyle{whitebox}=[rectangle, draw, fill=none, minimum size=0.3cm];
    \tikzstyle{blackbox}=[rectangle, draw, fill=black, minimum size=0.3cm];
    \tikzstyle{textnode}=[fill=none, draw=none, minimum size=0cm];

    \newcommand\layersep{1cm}
    \newcommand\unitH{0.6cm}

    \foreach \name / \y in {1,...,3}
    \node[textnode] (V-1-\name) at (0,-\y*\unitH) {$\vf_{\y}$};

    \foreach \name / \y in {1,...,3}
    \node[blackbox] (N-4-\name) at (\layersep,-\y*\unitH) {};
    \foreach \name / \y in {1,...,3}
    \draw[->] (V-1-\name) -- (N-4-\name);

    \foreach \name / \y in {1,...,3}
    \node[textnode] (O-\name) at (2*\layersep,-\y*\unitH) {$s_{\y}$};
    \foreach \name / \y in {1,...,3}
    \draw[->] (N-4-\name) -- (O-\name);

    \node[draw, dashed, fit=(V-1-1) (O-3), inner sep=1mm, rectangle, label=left:Approach 1, font=\footnotesize] (V) {};

    \node[whitebox] (Legend1) at (3*\layersep, -2*\unitH) {};
    \node[blackbox] (Legend2) at (3*\layersep, -2.8*\unitH) {};
    \node[textnode,font=\footnotesize] (T1) at (4.9*\layersep, -2*\unitH) {White-box calculations};
    \node[textnode,font=\footnotesize] (T2) at (4.5*\layersep, -2.8*\unitH) {Neural networks};
    \node[draw, dashed, fit=(Legend1) (Legend2) (T1) (T2), inner sep=1mm, rectangle, label=above:Legend] {};
\end{tikzpicture}

%% file: graph-aug.tikz
\begin{tikzpicture}
    \tikzstyle{box}=[rectangle, draw, minimum size=0.3cm]
    \tikzstyle{whitebox}=[box, fill=none];
    \tikzstyle{blackbox}=[box, fill=black];
    \tikzstyle{edgebox}=[box, fill=white, pos=0.5]; 
    \tikzstyle{textnode}=[fill=none, draw=none, minimum size=0cm];

    \newcommand\layersep{1.2cm}
    \newcommand\unitH{0.9cm}

    \foreach \name / \y in {1,...,3}
    \node[whitebox] (N-1-\name) at (\layersep,-\y*\unitH) {};
    \foreach \name / \y in {1,...,3}
    \node[textnode] (V-1-\name) at (0.4*\layersep,-\y*\unitH) {$\vf_{\y}$};
    \foreach \name in {1,...,3}
    \draw[->] (V-1-\name) -- (N-1-\name);

    \node[whitebox] (N-2-1) at (2.5*\layersep, -1*\unitH) {};
    \node[whitebox] (N-2-2) at (2.5*\layersep, -3*\unitH) {};
    \node[textnode] (C-2-1) at (2.5*\layersep, -1.7*\unitH) {$\vg_1$};
    \node[textnode] (C-2-2) at (2.5*\layersep, -2.3*\unitH) {$\vg_2$};
    \foreach \name in {1,...,2}
    \draw[->] (C-2-\name) -- (N-2-\name);

    \draw[->] (N-1-1) -- node[edgebox]{} (N-2-1);
    \draw[->] (N-1-2) -- node[edgebox]{} (N-2-1);
    \draw[->] (N-1-2) -- node[edgebox]{} (N-2-2);
    \draw[->] (N-1-3) -- node[edgebox]{} (N-2-2);

    \foreach \name / \y in {1,...,3}
    \node[whitebox] (N-3-\name) at (4*\layersep,-\y*\unitH) {};

    \draw[->] (N-2-1) -- node[edgebox]{} (N-3-1);
    \draw[->] (N-2-1) -- node[edgebox]{} (N-3-2);
    \draw[->] (N-2-2) -- node[edgebox]{} (N-3-2);
    \draw[->] (N-2-2) -- node[edgebox]{} (N-3-3);

    \foreach \name / \y in {1,...,3}
    \node[blackbox] (N-4-\name) at (5.2*\layersep,-\y*\unitH) {};
    \foreach \name / \y in {1,...,3}
    \draw[->] (N-3-\name) -- node[textnode,above](A-\name){$\vf_{\y,\mathrm{aug}}$} (N-4-\name);

    \foreach \name / \y in {1,...,3}
    \node[textnode] (O-\name) at (5.8*\layersep,-\y*\unitH) {$s_{\y}$};
    \foreach \name / \y in {1,...,3}
    \draw[->] (N-4-\name) -- (O-\name);

    \node[draw, dashed, fit=(V-1-1) (O-3) (A-1), inner sep=1mm, rectangle, label=left:Approach 2] (V) {};
\end{tikzpicture}

%% file: gnn.tikz
\begin{tikzpicture}
    \tikzstyle{box}=[rectangle, draw, minimum size=0.3cm]
    \tikzstyle{whitebox}=[box, fill=black];
    \tikzstyle{blackbox}=[box, fill=black];
    \tikzstyle{edgebox}=[box, fill=black, pos=0.5]; 
    \tikzstyle{textnode}=[fill=none, draw=none, minimum size=0cm];

    \newcommand\layersep{1.2cm}
    \newcommand\unitH{0.9cm}

    \foreach \name / \y in {1,...,3}
    \node[whitebox] (N-1-\name) at (\layersep,-\y*\unitH) {};
    \foreach \name / \y in {1,...,3}
    \node[textnode] (V-1-\name) at (0.4*\layersep,-\y*\unitH) {$\vf_{\y}$};
    \foreach \name in {1,...,3}
    \draw[->] (V-1-\name) -- (N-1-\name);

    \node[whitebox] (N-2-1) at (2.5*\layersep, -1*\unitH) {};
    \node[whitebox] (N-2-2) at (2.5*\layersep, -3*\unitH) {};
    \node[textnode] (C-2-1) at (2.5*\layersep, -1.7*\unitH) {$\vg_1$};
    \node[textnode] (C-2-2) at (2.5*\layersep, -2.3*\unitH) {$\vg_2$};
    \foreach \name in {1,...,2}
    \draw[->] (C-2-\name) -- (N-2-\name);

    \draw[->] (N-1-1) -- node[edgebox]{} (N-2-1);
    \draw[->] (N-1-2) -- node[edgebox]{} (N-2-1);
    \draw[->] (N-1-2) -- node[edgebox]{} (N-2-2);
    \draw[->] (N-1-3) -- node[edgebox]{} (N-2-2);

    \foreach \name / \y in {1,...,3}
    \node[whitebox] (N-3-\name) at (4*\layersep,-\y*\unitH) {};

    \draw[->] (N-2-1) -- node[edgebox]{} (N-3-1);
    \draw[->] (N-2-1) -- node[edgebox]{} (N-3-2);
    \draw[->] (N-2-2) -- node[edgebox]{} (N-3-2);
    \draw[->] (N-2-2) -- node[edgebox]{} (N-3-3);

    \foreach \name / \y in {1,...,3}
    \node[blackbox] (N-4-\name) at (5.2*\layersep,-\y*\unitH) {};
    \foreach \name / \y in {1,...,3}
    \draw[->] (N-3-\name) -- node[textnode,above](A-\name){$\vf_{\y,\mathrm{gnn}}$}  (N-4-\name);

    \foreach \name / \y in {1,...,3}
    \node[textnode] (O-\name) at (5.8*\layersep,-\y*\unitH) {$s_{\y}$};
    \foreach \name / \y in {1,...,3}
    \draw[->] (N-4-\name) -- (O-\name);

    \node[draw, dashed, fit=(V-1-1) (O-3) (A-1), inner sep=1mm, rectangle, label=left:Approach 3] (V) {};
\end{tikzpicture}

%% file: cuts.tikz
\begin{tikzpicture}

    \fill[black,opacity=0.2] (0,2) -- (1.5,2) -- (2.5,1) -- (2.5,0) -- (0,0) -- cycle; 
    \draw (0,2) -- (1.5,2) -- (2.5,1) -- (2.5,0) -- (0,0) -- cycle; 

    \draw[thick,->] (0,0) -- (3,0) node[right] {$x_1$};
    \draw[thick,->] (0,0) -- (0,2.3) node[above] {$x_2$};

    \foreach \x in {0,1,2}
        \foreach \y in {0,1,2}
            \fill (\x,\y) circle (3pt);
    \fill[white] (2,2) circle (4pt);

    \draw[dashed] (0.8,2.2) -- (2.6,0.4) node[right] {~~Cut 2: $x_1 + x_2 \leq 3$};
    \draw[dashed] (0.6,2.2) -- (2.6,1.2) node[right] {~~Cut 1: $x_1 + 2x_2 \leq 5$};

    \node[draw, circle, fill=white, inner sep = 2pt] at (1.5,2) {};
    \node[draw, rectangle] at (1,2) {};
    \node[draw, regular polygon, regular polygon sides=3, inner sep=1.2pt,fill=white] at (2,1.5) {};

    \draw (7,1.3) circle (3pt) node[right](Leg1) {~The current LP optimal solution $\underline{\vx}$};
    \fill (7,1.8) circle (3pt) node[right](Leg2) {~A feasible solution $\vx \in \sX_{\mathrm{MILP}}$};
    \node[draw, rectangle, inner sep=3pt, label=right:An LP solution after adding Cut 2](Leg3) at (7,0.3) {};
    \node[draw, regular polygon, regular polygon sides=3, inner sep=1.5pt, label=right:An LP solution after adding Cut 1](Leg4) at (7,0.8) {};

    \node[draw, dashed, fit=(Leg1) (Leg2) (Leg3) (Leg4), inner sep=2mm, rectangle, label = above:Legend] {};

    \draw[->] (2.8,2.3) -- (2.5,2) node[right] {~~$\vc = (-1,-1)$};

\end{tikzpicture}

%% file: local-branching.tikz
\begin{tikzpicture}
\tikzstyle{box}=[rectangle, draw, inner sep=0.1mm, rounded corners, minimum width=10cm];
\tikzstyle{textnode}=[fill=none, draw=none, minimum size=0cm];

\newcommand\deltax{0.8}
\newcommand\deltay{0.6}

\def\xstar{{0 1 1 1 0 0 1 0 0 1 0 1}} 
\def\xhat{{0 1 1 0 1 0 1 0 0 1 - -}} 
\def\corr{{T T T F F T T T T T - -}} 
\def\p{{0.1  0.9  0.9  0.2  0.8  0.3  0.7  0.3  0.3  0.7  0.6  0.6}} 
\def\f{{0.1  0.1  0.1  0.2  0.2  0.3  0.3  0.3  0.3  0.3  0.4  0.4}} 

\node[textnode] (p-1) at (1*\deltax, 0) {0.1};
\node[textnode] (p-2) at (2*\deltax, 0) {0.9};
\node[textnode] (p-3) at (3*\deltax, 0) {0.9};
\node[textnode] (p-4) at (4*\deltax, 0) {0.2};

\node[textnode] (p-5) at (5*\deltax, 0) {0.8};
\node[textnode] (p-6) at (6*\deltax, 0) {0.3};
\node[textnode] (p-7) at (7*\deltax, 0) {0.7};
\node[textnode] (p-8) at (8*\deltax, 0) {0.3};

\node[textnode] (p-9) at (9*\deltax, 0) {0.3};
\node[textnode] (p-10) at (10*\deltax, 0) {0.7};
\node[textnode] (p-11) at (11*\deltax, 0) {0.6};
\node[textnode] (p-12) at (12*\deltax, 0) {0.6};

\node[textnode] (f-1) at (1*\deltax, -\deltay) {0.9};
\node[textnode] (f-2) at (2*\deltax, -\deltay) {0.9};
\node[textnode] (f-3) at (3*\deltax, -\deltay) {0.9};
\node[textnode] (f-4) at (4*\deltax, -\deltay) {0.8};

\node[textnode] (f-5) at (5*\deltax, -\deltay) {0.8};
\node[textnode] (f-6) at (6*\deltax, -\deltay) {0.7};
\node[textnode] (f-7) at (7*\deltax, -\deltay) {0.7};
\node[textnode] (f-8) at (8*\deltax, -\deltay) {0.7};

\node[textnode] (f-9) at (9*\deltax, -\deltay) {0.7};
\node[textnode] (f-10) at (10*\deltax, -\deltay) {0.7};
\node[textnode] (f-11) at (11*\deltax, -\deltay) {0.6};
\node[textnode] (f-12) at (12*\deltax, -\deltay) {0.6};

\node[textnode] (xhat-1) at (1*\deltax, -2*\deltay) {0};
\node[textnode] (xhat-2) at (2*\deltax, -2*\deltay) {1};
\node[textnode] (xhat-3) at (3*\deltax, -2*\deltay) {1};
\node[textnode] (xhat-4) at (4*\deltax, -2*\deltay) {0};

\node[textnode] (xhat-5) at (5*\deltax, -2*\deltay) {1};
\node[textnode] (xhat-6) at (6*\deltax, -2*\deltay) {0};
\node[textnode] (xhat-7) at (7*\deltax, -2*\deltay) {1};
\node[textnode] (xhat-8) at (8*\deltax, -2*\deltay) {0};

\node[textnode] (xhat-9) at (9*\deltax, -2*\deltay) {0};
\node[textnode] (xhat-10) at (10*\deltax, -2*\deltay) {1};
\node[textnode] (xhat-11) at (11*\deltax, -2*\deltay) {-};
\node[textnode] (xhat-12) at (12*\deltax, -2*\deltay) {-};

\node[textnode] (xstar-1) at (1*\deltax, -3*\deltay) {0};
\node[textnode] (xstar-2) at (2*\deltax, -3*\deltay) {1};
\node[textnode] (xstar-3) at (3*\deltax, -3*\deltay) {1};
\node[textnode] (xstar-4) at (4*\deltax, -3*\deltay) {1};

\node[textnode] (xstar-5) at (5*\deltax, -3*\deltay) {0};
\node[textnode] (xstar-6) at (6*\deltax, -3*\deltay) {0};
\node[textnode] (xstar-7) at (7*\deltax, -3*\deltay) {1};
\node[textnode] (xstar-8) at (8*\deltax, -3*\deltay) {0};

\node[textnode] (xstar-9) at (9*\deltax, -3*\deltay) {0};
\node[textnode] (xstar-10) at (10*\deltax, -3*\deltay) {1};
\node[textnode] (xstar-11) at (11*\deltax, -3*\deltay) {0};
\node[textnode] (xstar-12) at (12*\deltax, -3*\deltay) {1};

\node[textnode] (cor-1) at (1*\deltax, -4*\deltay) {T};
\node[textnode] (cor-2) at (2*\deltax, -4*\deltay) {T};
\node[textnode] (cor-3) at (3*\deltax, -4*\deltay) {T};
\node[textnode] (cor-4) at (4*\deltax, -4*\deltay) {F};

\node[textnode] (cor-5) at (5*\deltax, -4*\deltay) {F};
\node[textnode] (cor-6) at (6*\deltax, -4*\deltay) {T};
\node[textnode] (cor-7) at (7*\deltax, -4*\deltay) {T};
\node[textnode] (cor-8) at (8*\deltax, -4*\deltay) {T};

\node[textnode] (cor-9) at (9*\deltax, -4*\deltay) {T};
\node[textnode] (cor-10) at (10*\deltax, -4*\deltay) {T};
\node[textnode] (cor-11) at (11*\deltax, -4*\deltay) {-};
\node[textnode] (cor-12) at (12*\deltax, -4*\deltay) {-};

\node[box, fit=(p-1)(p-12), label=left:$p_j$] (p) {};
\node[box, fit=(f-1)(f-12), label=left:Confidence] (f) {};
\node[box, fit=(xhat-1)(xhat-12), label=left:$\hat{x}_j$] (xhat) {};
\node[box, fit=(xstar-1)(xstar-12), label=left:$x^\ast_j$] (xstar) {};
\node[box, fit=(cor-1)(cor-12), label=left:Correctness] (cor) {};

\draw [dashed] (3.5*\deltax, 1*\deltay) -- (3.5*\deltax, -5*\deltay) node [below] {$\xi = 0.85$};
\draw [dashed] (10.5*\deltax, 1*\deltay) -- (10.5*\deltax, -5*\deltay) node [below] {$\xi=0.65$};

\end{tikzpicture}

%% file: ml4mip-summary.tikz
\begin{tikzpicture}
\tikzstyle{box}=[rectangle, draw, rounded corners, minimum height=0.6cm, minimum width=2.2cm,align=center];
\tikzstyle{box2}=[box,font=\footnotesize];
\newcommand\deltax{2.6cm}
\newcommand\deltay{1.4cm}
    \node[box2] (A) at (0,0) {Configuration};
    \node[box2] (B) at (\deltax,0) {LP relaxation};
    \node[box2] (C) at (2*\deltax,0) {Heuristics};
    \node[box2] (D) at (3*\deltax,0) {Cutting Plane};
    \node[box2] (E) at (4*\deltax,0) {BnB};

    \draw[->] (A) -- (B);
    \draw[->] (B) -- (C);
    \draw[->] (C) -- (D);
    \draw[->] (D) -- (E);

    \node[box2] (A1) at (0,-\deltay) {Learning to config \\ Section \ref{sec:config-mip}};
    \node[box2] (C1) at (1.5*\deltax,-\deltay) {Learning heuristics \\ Section \ref{sec:ml-based-heuristics}};
    \node[box2] (D1) at (2.7*\deltax,-\deltay) {Learning to select cuts \\ Section \ref{sec:l2c}};
    \node[box2] (E1) at (4*\deltax,-\deltay) {Learning to branch \\ Learning to search \\ Sections \ref{sec:l2b} \& \ref{sec:l2s}};

    \draw[->] (A1) to[bend left=15] (A);
    \draw[->] (C1) to[bend left=15] (C);
    \draw[->] (D1) to[bend left=15] (D);
    \draw[->] (E1) to[bend left=15] (E);

    \draw[->] (A) to[bend left=15] (A1);
    \draw[->] (C) to[bend left=15] (C1);
    \draw[->] (D) to[bend left=15] (D1);
    \draw[->] (E) to[bend left=15] (E1);

    \node[draw, dashed, fit=(A) (E), inner sep=1.5mm, rectangle, label=left:An MILP solver] (solver) {};
    \node[draw, dashed, fit=(A1) (E1), inner sep=1.5mm, rectangle, label=left:ML modules] (ml) {};
\end{tikzpicture}

%% file: main.bbl
\begin{thebibliography}{99}


\bibitem{aberdam2021adalista} Aberdam A, Golts A, Elad M.
Ada-lista: Learned solvers adaptive to varying models. IEEE Transactions on Pattern Analysis and Machine Intelligence. 2021 Nov 4;44(12):9222-35.

\bibitem{ablin2019learning} Ablin P, Moreau T, Massias M, Gramfort A.
Learning step sizes for unfolded sparse coding. Advances in Neural Information Processing Systems 32, 2019.

\bibitem{achterberg2007constraint} Achterberg, Tobias. Constraint integer programming, 2007.

\bibitem{adler2018learned} Adler J, Öktem O.
Learned primal-dual reconstruction. IEEE transactions on medical imaging, 37(6), 1322-1332, 2018.

\bibitem{Agrawal2019differentiable} Agrawal A, Amos B, Barratt S, et al.
Differentiable convex optimization layers. Advances in neural information processing systems 32, 2019.

\bibitem{aharon2006ksvd} Aharon M, Elad M, Bruckstein A.
K-SVD: An algorithm for designing overcomplete dictionaries for sparse representation. IEEE Transactions on signal processing. 2006 Oct 16;54(11):4311-22.

\bibitem{alvarez2017machine} Alvarez A M, Louveaux Q, Wehenkel L.
A machine learning-based approximation of strong branching. INFORMS Journal on Computing 29.1 (2017): 185-195.

\bibitem{amos2017optnet} Amos B, Kolter J Z.
Optnet: Differentiable optimization as a layer in neural networks. International Conference on Machine Learning. PMLR, 2017.

\bibitem{ansotegui2016maxsat} Ansótegui, Carlos, et al. "MaxSAT by improved instance-specific algorithm configuration." Artificial Intelligence 235 (2016): 26-39.

\bibitem{bai2019deep} Bai S, Kolter J Z, Koltun V.
Deep equilibrium models. Advances in Neural Information Processing Systems 32, 2019.

\bibitem{balas1980set} Balas E, Ho A.
Set covering algorithms using cutting planes, heuristics, and subgradient optimization: a computational study. Springer Berlin Heidelberg; 1980.

\bibitem{balas1996gomory} Balas E, et al.
Gomory cuts revisited. Operations Research Letters 19.1 (1996): 1-9.

\bibitem{balatsoukas2019deep} Balatsoukas-Stimming A, Studer C.
Deep unfolding for communications systems: A survey and some new directions. In 2019 IEEE International Workshop on Signal Processing Systems (SiPS) 2019 Oct 20 (pp. 266-271). IEEE.

\bibitem{balcan2018learning} Balcan M-F, et al.
Learning to branch. International conference on machine learning. PMLR, 2018.

\bibitem{bansal2018can} Bansal N, Chen X, Wang Z.
Can we gain more from orthogonality regularizations in training deep networks?. Advances in Neural Information Processing Systems 31, 2018.

\bibitem{bartlett2017spectrally} Bartlett PL, Foster DJ, Telgarsky MJ. Spectrally-normalized margin bounds for neural networks. Advances in neural information processing systems 30, 2017.

\bibitem{fista} Beck A, Teboulle M.
A fast iterative shrinkage-thresholding algorithm with application to wavelet-based image deblurring. In2009 IEEE International Conference on Acoustics, Speech and Signal Processing 2009 Apr 19 (pp. 693-696). IEEE.

\bibitem{behboodi2020generalization} Behboodi A, Rauhut H, Schnoor E.
Compressive sensing and neural networks from a statistical learning perspective. In Compressed Sensing in Information Processing 2022 Oct 22 (pp. 247-277). Cham: Springer International Publishing.

\bibitem{behrens2020neurally} Behrens F, Sauder J, Jung P.
Neurally Augmented ALISTA. In International Conference on Learning Representations 2020 Oct 2.

\bibitem{bergstra2012random} Bergstra J, Bengio Y. Random search for hyper-parameter optimization. Journal of machine learning research. 2012 Feb 1;13(2).

\bibitem{berthet2020learning} Berthet Q, Blondel M, Teboul O, et al. Learning with differentiable pertubed optimizers. Advances in neural information processing systems. 2020;33:9508-19.

\bibitem{berthold2022learning} Berthold T, Francobaldi M, Hendel G. Learning to use local cuts. arXiv preprint arXiv:2206.11618. 2022 Jun 23.

\bibitem{berthold2006primal} Berthold T. Primal heuristics for mixed integer programs (Doctoral dissertation, Zuse Institute Berlin (ZIB)).

\bibitem{bertsimas1997lp} Bertsimas D, Tsitsiklis JN. Introduction to linear optimization. Belmont, MA: Athena Scientific; 1997 Jan.

\bibitem{Bertsimas2020predictive} Bertsimas D, Kallus N. From predictive to prescriptive analytics. Management Science. 2020 Mar;66(3):1025-44.

\bibitem{scip8} Bestuzheva K, Besançon M, Chen WK, et al. The SCIP optimization suite 8.0. arXiv preprint arXiv:2112.08872. 2021 Dec 16.

\bibitem{bonami2022classifier} Bonami P, Lodi A, Zarpellon G. A classifier to decide on the linearization of mixed-integer quadratic problems in CPLEX. Operations research. 2022 Nov;70(6):3303-20.

\bibitem{borgerding2017amp} Borgerding M, Schniter P, Rangan S.
AMP-inspired deep networks for sparse linear inverse problems. IEEE Transactions on Signal Processing. 2017 May 25;65(16):4293-308.

\bibitem{boyd2011distributed} Boyd S, Parikh N, Chu E, Peleato B, Eckstein J.
Distributed optimization and statistical learning via the alternating direction method of multipliers. Foundations and Trends® in Machine learning. 2011 Jul 25;3(1):1-22.

\bibitem{brauer2022learning} Brauer C, Breustedt N, De Wolff T, Lorenz DA. Learning variational models with unrolling and bilevel optimization. arXiv preprint arXiv:2209.12651. 2022 Sep 26.

\bibitem{brock2019large} Brock A, Donahue J, Simonyan K. Large Scale GAN Training for High Fidelity Natural Image Synthesis. In International Conference on Learning Representations 2018 Sep 27.

\bibitem{wilder2019melding} Wilder B, Dilkina B, Tambe M. Melding the data-decisions pipeline: Decision-focused learning for combinatorial optimization. In Proceedings of the AAAI Conference on Artificial Intelligence 2019 Jul 17 (Vol. 33, No. 01, pp. 1658-1665).

\bibitem{buades2005nonlocal} Buades A, Coll B, Morel JM. A non-local algorithm for image denoising. In 2005 IEEE computer society conference on computer vision and pattern recognition (CVPR'05) 2005 Jun 20 (Vol. 2, pp. 60-65). IEEE.

\bibitem{cai2021learned} Cai H, Liu J, Yin W. Learned robust pca: A scalable deep unfolding approach for high-dimensional outlier detection. Advances in Neural Information Processing Systems. 2021 Dec 6;34:16977-89.

\bibitem{cappart2021combining} Cappart Q, Moisan T, Rousseau LM, Prémont-Schwarz I, Cire AA. Combining reinforcement learning and constraint programming for combinatorial optimization. In Proceedings of the AAAI Conference on Artificial Intelligence 2021 May 18 (Vol. 35, No. 5, pp. 3677-3687).

\bibitem{chan2016plug} Chan SH, Wang X, Elgendy OA. Plug-and-play ADMM for image restoration: Fixed-point convergence and applications. IEEE Transactions on Computational Imaging. 2016 Nov 15;3(1):84-98.

\bibitem{chen2022learning} Chen T, Chen X, Chen W, et al. Learning to Optimize: A Primer and A Benchmark. Journal of Machine Learning Research. 2022;23:1-59.

\bibitem{chen2020learning} Chen X, Dai H, Li Y, et al. Learning to stop while learning to predict. In International conference on machine learning 2020 Nov 21 (pp. 1520-1530). PMLR.

\bibitem{chen2021hyperparameter} Chen X, Liu J, Wang Z, Yin W. Hyperparameter tuning is all you need for LISTA. Advances in Neural Information Processing Systems. 2021 Dec 6;34:11678-89.

\bibitem{chen2018theoretical} Chen X, Liu J, Wang Z, Yin W. Theoretical linear convergence of unfolded ISTA and its practical weights and thresholds. Advances in Neural Information Processing Systems. 2018;31.

\bibitem{chen2024rethinking} Chen Z, Liu J, Chen X, et al. Rethinking the Capacity of Graph Neural Networks for Branching Strategy. arXiv preprint arXiv:2402.07099 (2024).

\bibitem{chen2022lp} Chen Z, Liu J, Wang X, Yin W. On Representing Linear Programs by Graph Neural Networks. In The Eleventh International Conference on Learning Representations 2023.

\bibitem{chen2022milp} Chen Z, Liu J, Wang X, Yin W. On Representing Mixed-Integer Linear Programs by Graph Neural Networks. In The Eleventh International Conference on Learning Representations 2023.

\bibitem{chen2020understanding} Chen X, Zhang Y, Reisinger C, Song L. Understanding deep architecture with reasoning layer. Advances in Neural Information Processing Systems. 2020;33:1240-52.

\bibitem{chmiela2021learning} Chmiela A, Khalil E, Gleixner A, et al. Learning to schedule heuristics in branch and bound. Advances in Neural Information Processing Systems. 2021 Dec 6;34:24235-46.

\bibitem{cohen2021regularization} Cohen R, Elad M, Milanfar P. Regularization by Denoising via Fixed-Point Projection (RED-PRO). SIAM Journal on Imaging Sciences. 2021;14(3):1374-406.

\bibitem{condat} Condat L. A primal–dual splitting method for convex optimization involving Lipschitzian, proximable and linear composite terms. Journal of optimization theory and applications. 2013 Aug;158(2):460-79.

\bibitem{corbineau2019learned} Corbineau MC, Bertocchi C, Chouzenoux E, et al. Learned image deblurring by unfolding a proximal interior point algorithm. In2019 IEEE International Conference on Image Processing (ICIP) 2019 Sep 22 (pp. 4664-4668). IEEE.

\bibitem{dabov2017image} Dabov K, Foi A, Katkovnik V, Egiazarian K. Image denoising by sparse 3-D transform-domain collaborative filtering. IEEE Transactions on image processing. 2007 Jul 16;16(8):2080-95.

\bibitem{McKenzie2023faster} McKenzie D, Fung SW, Heaton H. Faster Predict-and-Optimize with Davis-Yin Splitting. arXiv preprint arXiv:2301.13395. 2023 Jan 31.

\bibitem{davis2017three} Davis D, Yin W. A three-operator splitting scheme and its optimization applications. Set-valued and variational analysis. 2017 Dec;25:829-58.

\bibitem{deze2023machine} Deza A, Khalil EB. Machine learning for cutting planes in integer programming: A survey. arXiv preprint arXiv:2302.09166. 2023 Feb 17.

\bibitem{ding2020accelerating} Ding JY, Zhang C, Shen L, et al. Accelerating primal solution findings for mixed integer programs based on solution prediction. In Proceedings of the AAAI conference on artificial intelligence 2020 Apr 3 (Vol. 34, No. 02, pp. 1452-1459).

\bibitem{Elmachtoub2022smart} Elmachtoub AN, Grigas P. Smart “predict, then optimize”. Management Science. 2022 Jan;68(1):9-26.

\bibitem{ryuyin2022} Ryu EK, Yin W. Large-scale convex optimization: algorithms \& analyses via monotone operators. Cambridge University Press; 2022 Dec 1.

\bibitem{etheve2020reinforcement} Etheve M, Alès Z, Bissuel C, et al. Reinforcement learning for variable selection in a branch and bound algorithm. In International Conference on Integration of Constraint Programming, Artificial Intelligence, and Operations Research 2020 Sep 19 (pp. 176-185). Cham: Springer International Publishing.

\bibitem{falkner2022large} Falkner JK, Thyssens D, Schmidt-Thieme L. Large neighborhood search based on neural construction heuristics. arXiv preprint arXiv:2205.00772. 2022 May 2.

\bibitem{fan2001variable} Fan J, Li R. Variable selection via nonconcave penalized likelihood and its oracle properties. Journal of the American statistical Association. 2001 Dec 1;96(456):1348-60.

\bibitem{fischetti2003local} Fischetti M, Lodi A. Local branching. Mathematical programming. 2003 Sep;98:23-47.

\bibitem{gasse2019exact} Gasse M, Chételat D, Ferroni N, et al.
Exact combinatorial optimization with graph convolutional neural networks. Advances in neural information processing systems. 2019;32.

\bibitem{gehring2016convolutional} Gehring J, Auli M, Grangier D, Dauphin YN.
A convolutional encoder model for neural machine translation. arXiv preprint arXiv:1611.02344. 2016 Nov 7.

\bibitem{geng2021training} Geng Z, Zhang XY, Bai S, et al. On training implicit models. Advances in Neural Information Processing Systems. 2021 Dec 6;34:24247-60.

\bibitem{giryes2018tradeoffs} Giryes R, Eldar YC, Bronstein AM, Sapiro G.
Tradeoffs between convergence speed and reconstruction accuracy in inverse problems. IEEE Transactions on Signal Processing. 2018 Jan 11;66(7):1676-90.

\bibitem{gogna2013meta} Gogna A, Tayal A. Metaheuristics: review and application. Journal of Experimental \& Theoretical Artificial Intelligence. 2013 Dec 1;25(4):503-26.

\bibitem{gomory1958algorithm} Gomory RE. An Algorithm for Integer Solutions to Lmear Programs. Princeton-IBM Mathematics Research Project Technical Report 1 (1958).

\bibitem{gomory1960solving} Gomory RE. Solving linear programming problems in integers. Combinatorial Analysis 10 (1960): 211-215.

\bibitem{deeplearning} Goodfellow I, Bengio Y, Courville A. Deep learning. MIT press; 2016 Nov 10.

\bibitem{lista} Gregor K, LeCun Y. Learning fast approximations of sparse coding.
In Proceedings of the 27th international conference on international conference on machine learning 2010 Jun 21 (pp. 399-406).

\bibitem{griewank2008evaluating} Griewank A, Walther A. Evaluating derivatives: principles and techniques of algorithmic differentiation. Society for industrial and applied mathematics; 2008 Jan 1.

\bibitem{gupta2018cnn} Gupta H, Jin KH, Nguyen HQ, et al. CNN-based projected gradient descent for consistent CT image reconstruction. IEEE transactions on medical imaging. 2018 May 3;37(6):1440-53.

\bibitem{gupta2020hybrid} Gupta P, Gasse M, Khalil E, et al. Hybrid models for learning to branch. Advances in neural information processing systems. 2020;33:18087-97.

\bibitem{gupta2022lookback} Gupta P, Khalil EB, Chetélat D, et al. Lookback for learning to branch. arXiv preprint arXiv:2206.14987, 2022.

\bibitem{han2013online} Han S, Fu R, Wang S, Wu X. Online adaptive dictionary learning and weighted sparse coding for abnormality detection. In 2013 IEEE International Conference on Image Processing 2013 Sep 15 (pp. 151-155). IEEE.

\bibitem{hauptmann2018model} Hauptmann A, Lucka F, Betcke M, et al. Model-based learning for accelerated, limited-view 3-D photoacoustic tomography. IEEE transactions on medical imaging. 2018 Mar 29;37(6):1382-93.

\bibitem{he2020model} He H, Wen CK, Jin S, Li GY. Model-driven deep learning for MIMO detection. IEEE Transactions on Signal Processing. 2020 Feb 28;68:1702-15.

\bibitem{he2016deep} He K, Zhang X, Ren S, Sun J. Deep residual learning for image recognition. In Proceedings of the IEEE conference on computer vision and pattern recognition 2016 (pp. 770-778).

\bibitem{he2014learning} He H, Daume III H, Eisner JM. Learning to search in branch and bound algorithms. Advances in neural information processing systems. 2014;27.

\bibitem{heaton2023safeguarded} Heaton H, Chen X, Wang Z, Yin W. Safeguarded learned convex optimization. In Proceedings of the AAAI Conference on Artificial Intelligence 2023 Jun 26 (Vol. 37, No. 6, pp. 7848-7855).

\bibitem{heaton2022wasserstein} Heaton H, Fung SW, Lin AT, Osher S, Yin W.
Wasserstein-based projections with applications to inverse problems. SIAM Journal on Mathematics of Data Science. 2022;4(2):581-603.

\bibitem{hendel2022adaptive} Hendel G. Adaptive large neighborhood search for mixed integer programming. Mathematical Programming Computation. 2022 Jun;14(2):185-221.

\bibitem{himmich2023mpils} Himmich I, El Hachemi N, El Hallaoui I, et al. MPILS: An automatic tuner for MILP solvers. Computers \& Operations Research. 2023 Nov 1;159:106344.

\bibitem{hornik1989multilayer} Hornik K, Stinchcombe M, White H. Multilayer feedforward networks are universal approximators. Neural networks. 1989 Jan 1;2(5):359-66.

\bibitem{hosny2023auto} Hosny A, Reda S. Automatic MILP solver configuration by learning problem similarities. Annals of Operations Research. 2023 Jul 14:1-28.

\bibitem{hottung2020neural} Hottung A, Tierney K. Neural large neighborhood search for the capacitated vehicle routing problem. arXiv preprint arXiv:1911.09539. 2019 Nov 21.

\bibitem{huang2022improving} Huang L, Chen X, Huo W, et al. Improving primal heuristics for mixed integer programming problems based on problem reduction: A learning-based approach. In 2022 17th International Conference on Control, Automation, Robotics and Vision (ICARCV) 2022 Dec 11 (pp. 181-186). IEEE.

\bibitem{huang2022anytime} Huang T, Li J, Koenig S, Dilkina B. Anytime multi-agent path finding via machine learning-guided large neighborhood search. In Proceedings of the AAAI Conference on Artificial Intelligence 2022 Jun 28 (Vol. 36, No. 9, pp. 9368-9376).

\bibitem{huang2023searching} Huang T, Ferber AM, Tian Y, et al. Searching large neighborhoods for integer linear programs with contrastive learning. In International Conference on Machine Learning 2023 Jul 3 (pp. 13869-13890). PMLR.

\bibitem{huang2022learning} Huang Z, Wang K, Liu F, et al. Learning to select cuts for efficient mixed-integer programming. Pattern Recognition. 2022 Mar 1;123:108353.

\bibitem{hutter2011smac} Hutter F, Hoos HH, Leyton-Brown K. Sequential model-based optimization for general algorithm configuration. InLearning and Intelligent Optimization: 5th International Conference, LION 5, Rome, Italy, January 17-21, 2011. Selected Papers 5 2011 (pp. 507-523). Springer Berlin Heidelberg.

\bibitem{hutter2009param} Hutter F, Hoos HH, Leyton-Brown K, Stützle T. ParamILS: an automatic algorithm configuration framework. Journal of Artificial Intelligence Research. 2009 Oct 30;36:267-306.

\bibitem{ioffe2015batch} Ioffe S, Szegedy C. Batch normalization: Accelerating deep network training by reducing internal covariate shift. In International conference on machine learning 2015 Jun 1 (pp. 448-456). PMLR.

\bibitem{mandi2020interior} Mandi J, Guns T. Interior point solving for lp-based prediction+ optimisation. Advances in Neural Information Processing Systems. 2020;33:7272-82.

\bibitem{jegelka2022gnn} Jegelka S. Theory of graph neural networks: Representation and learning. In The International Congress of Mathematicians 2022.

\bibitem{Bolte2023onestep} Bolte J, Pauwels E, Vaiter S. One-step differentiation of iterative algorithms. Advances in Neural Information Processing Systems. 2024 Feb 13;36.

\bibitem{jia2021benders} Jia H, Shen S. Benders cut classification via support vector machines for solving two-stage stochastic programs. INFORMS Journal on Optimization. 2021 Jul;3(3):278-97.

\bibitem{joukovsky2021generalization} Joukovsky B, Mukherjee T, Van Luong H, Deligiannis N. Generalization error bounds for deep unfolding RNNs. In Uncertainty in Artificial Intelligence 2021 Dec 1 (pp. 1515-1524). PMLR.

\bibitem{kadioglu2010isac} Kadioglu S, Malitsky Y, Sellmann M, Tierney K. ISAC–instance-specific algorithm configuration. In ECAI 2010 2010 (pp. 751-756). IOS Press.

\bibitem{kang2018deep} Kang E, Chang W, Yoo J, Ye JC. Deep convolutional framelet denosing for low-dose CT via wavelet residual network. IEEE transactions on medical imaging. 2018 Apr 6;37(6):1358-69.

\bibitem{khalil2022mip} Khalil EB, Morris C, Lodi A. Mip-gnn: A data-driven framework for guiding combinatorial solvers. In Proceedings of the AAAI Conference on Artificial Intelligence 2022 Jun 28 (Vol. 36, No. 9, pp. 10219-10227).

\bibitem{khalil2017learning} Khalil E, Dai H, Zhang Y, et al. Learning combinatorial optimization algorithms over graphs. Advances in neural information processing systems. 2017;30.

\bibitem{khalil2016learning} Khalil E, Le Bodic P, Song L, et al. Learning to branch in mixed integer programming. In Proceedings of the AAAI Conference on Artificial Intelligence 2016 Feb 21 (Vol. 30, No. 1).

\bibitem{adam} Kingma DP, Ba J. Adam: A method for stochastic optimization. arXiv preprint arXiv:1412.6980. 2014 Dec 22.

\bibitem{kouni2022generalization} Kouni V, Panagakis Y. DECONET: an Unfolding Network for Analysis-based Compressed Sensing with Generalization Error Bounds. IEEE Transactions on Signal Processing. 2023 May 3.

\bibitem{labassi2022learning} Labassi AG, Chételat D, Lodi A. Learning to compare nodes in branch and bound with graph neural networks. Advances in neural information processing systems. 2022 Dec 6;35:32000-10.

\bibitem{mnist} LeCun Y, Cortes C, Burges JC C. THE MNIST DATABASE of handwritten digits. \url{http://yann.lecun.com/exdb/mnist/}.

\bibitem{li2021deep} Li Y, Bar-Shira O, Monga V, Eldar YC. Deep algorithm unrolling for biomedical imaging. arXiv preprint arXiv:2108.06637. 2021 Aug 15.

\bibitem{lin2022learning} Lin J, Zhu J, Wang H, Zhang T. Learning to branch with Tree-aware Branching Transformers. Knowledge-Based Systems. 2022 Sep 27;252:109455.

\bibitem{liu2023towards} Liu J, Chen X, Wang Z, Yin W, Cai H. Towards Constituting Mathematical Structures for Learning to Optimize. arXiv preprint arXiv:2305.18577. 2023 May 29.

\bibitem{alista} Liu J, Chen X, Wang Z, Yin W. ALISTA: Analytic Weights Are As Good As Learned Weights in LISTA. In International Conference on Learning Representations 2018 Sep 27.

\bibitem{liu2022learning} Liu D, Fischetti M, Lodi A. Learning to search in local branching. In Proceedings of the AAAI conference on artificial intelligence 2022 Jun 28 (Vol. 36, No. 4, pp. 3796-3803).

\bibitem{unet} Long J, Shelhamer E, Darrell T. Fully convolutional networks for semantic segmentation. In Proceedings of the IEEE conference on computer vision and pattern recognition 2015 (pp. 3431-3440).

\bibitem{ma2022efficient} Ma Y, Li J, Cao Z, et al. Efficient neural neighborhood search for pickup and delivery problems. arXiv preprint arXiv:2204.11399. 2022 Apr 25.

\bibitem{malitsky2014isac} Malitsky Y. Instance-specific algorithm configuration. Springer International Publishing; 2014.

\bibitem{mao2016image} Mao X, Shen C, Yang YB. Image restoration using very deep convolutional encoder-decoder networks with symmetric skip connections. Advances in neural information processing systems. 2016;29.

\bibitem{alvarez2014supervised} Marcos Alvarez A, Louveaux Q, Wehenkel L. A supervised machine learning approach to variable branching in branch-and-bound. 2014.

\bibitem{mardani2018neural} Mardani M, Sun Q, Donoho D, et al. Neural proximal gradient descent for compressive imaging. Advances in Neural Information Processing Systems. 2018;31.

\bibitem{meinhardt2017learning} Meinhardt T, Moller M, Hazirbas C, Cremers D. Learning proximal operators: Using denoising networks for regularizing inverse imaging problems. In Proceedings of the IEEE International Conference on Computer Vision 2017 (pp. 1781-1790).

\bibitem{miyato2018spectral} Miyato T, Kataoka T, Koyama M, Yoshida Y. Spectral Normalization for Generative Adversarial Networks. In International Conference on Learning Representations 2018 Feb 15.

\bibitem{mnih2015human} Mnih V, Kavukcuoglu K, Silver D, et al. Human-level control through deep reinforcement learning. nature. 2015 Feb 26;518(7540):529-33.

\bibitem{monga2019algorithm} Monga V, Li Y, Eldar YC. Algorithm unrolling: Interpretable, efficient deep learning for signal and image processing. IEEE Signal Processing Magazine. 2021 Feb 25;38(2):18-44.

\bibitem{moreau2017understanding} Moreau T, Bruna J. Understanding neural sparse coding with matrix factorization. In International Conference on Learning Representation (ICLR) 2017 Apr.

\bibitem{nair2020solving} Nair V, Bartunov S, Gimeno F, et al. Solving mixed integer programs using neural networks. arXiv preprint arXiv:2012.13349. 2020 Dec 23.

\bibitem{oberman2018lipschitz} Oberman AM, Calder J. Lipschitz regularized deep neural networks converge and generalize. arXiv preprint arXiv:1808.09540. 2018 Aug 28.

\bibitem{parsonson2023reinforcement} Parsonson CW, Laterre A, Barrett TD. Reinforcement learning for branch-and-bound optimisation using retrospective trajectories. In Proceedings of the AAAI Conference on Artificial Intelligence 2023 Jun 26 (Vol. 37, No. 4, pp. 4061-4069).

\bibitem{pascanu2013on} Pascanu R, Mikolov T, Bengio Y. On the difficulty of training recurrent neural networks. In International conference on machine learning 2013 May 26 (pp. 1310-1318). PMLR.

\bibitem{paulus2023learning} Paulus M, Krause A. Learning to dive in branch and bound. Advances in Neural Information Processing Systems. 2024 Feb 13;36.

\bibitem{paulus2022learning} Paulus M, Zarpellon G, Krause A, et al. Learning to cut by looking ahead: Cutting plane selection via imitation learning. In International conference on machine learning 2022 Jun 28 (pp. 17584-17600). PMLR.

\bibitem{pramanik2020deep} Pramanik A, Aggarwal HK, Jacob M. Deep generalization of structured low-rank algorithms (Deep-SLR). IEEE transactions on medical imaging. 2020 Aug 5;39(12):4186-97.

\bibitem{Donti2017taskbased} Donti P, Amos B, Kolter JZ. Task-based end-to-end model learning in stochastic optimization. Advances in neural information processing systems. 2017;30.

\bibitem{ecole} Prouvost A, Dumouchelle J, Scavuzzo L, et al. Ecole: A Gym-like Library for Machine Learning in Combinatorial Optimization Solvers. In Learning Meets Combinatorial Algorithms at NeurIPS 2020, 2020 Nov 11.

\bibitem{qian2018l2} Qian H, Wegman MN. L2-Nonexpansive Neural Networks. In International Conference on Learning Representations 2018 Sep 27.

\bibitem{qu2022improved} Qu Q, Li X, Zhou Y, et al. An improved reinforcement learning algorithm for learning to branch. arXiv preprint arXiv:2201.06213. 2022 Jan 17.

\bibitem{rick2017one} Rick Chang JH, Li CL, Poczos B, Vijaya Kumar BV, Sankaranarayanan AC. One network to solve them all--solving linear inverse problems using deep projection models. In Proceedings of the IEEE International Conference on Computer Vision 2017 (pp. 5888-5897).

\bibitem{rudin1992nonlinear} Rudin LI, Osher S, Fatemi E. Nonlinear total variation based noise removal algorithms. Physica D: nonlinear phenomena. 1992 Nov 1;60(1-4):259-68.

\bibitem{ryu2019plug} Ryu E, Liu J, Wang S, et al. Plug-and-play methods provably converge with properly trained denoisers. In International Conference on Machine Learning 2019 May 24 (pp. 5546-5557). PMLR.

\bibitem{samuel2019learning} Samuel N, Diskin T, Wiesel A. Learning to detect. IEEE Transactions on Signal Processing. 2019 Feb 15;67(10):2554-64.

\bibitem{scarlett2022theoretical} Scarlett J, Heckel R, Rodrigues MR, et al. Theoretical perspectives on deep learning methods in inverse problems. IEEE journal on selected areas in information theory. 2022 Sep;3(3):433-53.

\bibitem{scavuzzo2022learning} Scavuzzo L, Chen F, Chételat D, et al. Learning to branch with tree mdps. Advances in Neural Information Processing Systems. 2022 Dec 6;35:18514-26.

\bibitem{schnoor2021generlization} Schnoor E, Behboodi A, Rauhut H. Generalization error bounds for iterative recovery algorithms unfolded as neural networks. Information and Inference: A Journal of the IMA. 2023 Sep;12(3):2267-99.

\bibitem{shen2021learning} Shen Y, Sun Y, Eberhard A, Li X. Learning primal heuristics for mixed integer programs. In 2021 international joint conference on neural networks (IJCNN) 2021 Jul 18 (pp. 1-8). IEEE.

\bibitem{alphago} Silver D, Huang A, Maddison CJ, et al. Mastering the game of Go with deep neural networks and tree search. nature. 2016 Jan;529(7587):484-9.

\bibitem{vgg} Simonyan K, Zisserman A. Very deep convolutional networks for large-scale image recognition. arXiv preprint arXiv:1409.1556. 2014 Sep 4.

\bibitem{snoek2012practical} Snoek J, Larochelle H, Adams RP. Practical bayesian optimization of machine learning algorithms. Advances in neural information processing systems. 2012;25.

\bibitem{solomon2019deep} Solomon O, Cohen R, Zhang Y, et al. Deep unfolded robust PCA with application to clutter suppression in ultrasound. IEEE transactions on medical imaging. 2019 Sep 13;39(4):1051-63.

\bibitem{song2020general} Song J, Yue Y, Dilkina B. A general large neighborhood search framework for solving integer linear programs. Advances in Neural Information Processing Systems. 2020;33:20012-23.

\bibitem{song2023instance} Song W, Liu Y, Cao Z, et al. Instance-specific algorithm configuration via unsupervised deep graph clustering. Engineering Applications of Artificial Intelligence. 2023 Oct 1;125:106740.

\bibitem{sonnerat2021learning} Sonnerat N, Wang P, Ktena I, et al. Learning a large neighborhood search algorithm for mixed integer programs. arXiv preprint arXiv:2107.10201. 2021 Jul 21.

\bibitem{sreehari2016plug} Sreehari S, Venkatakrishnan SV, Wohlberg B, et al. Plug-and-play priors for bright field electron tomography and sparse interpolation. IEEE Transactions on Computational Imaging. 2016 Aug 11;2(4):408-23.

\bibitem{sreter2018learned} Sreter H, Giryes R. Learned convolutional sparse coding. In 2018 IEEE International Conference on Acoustics, Speech and Signal Processing (ICASSP) 2018 Apr 15 (pp. 2191-2195). IEEE.

\bibitem{sutton2000policy} Sutton RS, McAllester D, Singh S, Mansour Y. Policy gradient methods for reinforcement learning with function approximation. Advances in neural information processing systems. 1999;12.

\bibitem{takabe2020complex} Takabe S, Wadayama T, Eldar YC. Complex trainable ista for linear and nonlinear inverse problems. In ICASSP 2020-2020 IEEE International Conference on Acoustics, Speech and Signal Processing (ICASSP) 2020 May 4 (pp. 5020-5024). IEEE.

\bibitem{takabe2020theoretical} Takabe S, Wadayama T. Theoretical interpretation of learned step size in deep-unfolded gradient descent. arXiv preprint arXiv:2001.05142. 2020 Jan 15.

\bibitem{tang2020reinforcement} Tang Y, Agrawal S, Faenza Y. Reinforcement learning for integer programming: Learning to cut. In International conference on machine learning 2020 Nov 21 (pp. 9367-9376). PMLR.

\bibitem{teerapittayanon2016branchynet} Teerapittayanon S, McDanel B, Kung HT. Branchynet: Fast inference via early exiting from deep neural networks. In 2016 23rd international conference on pattern recognition (ICPR) 2016 Dec 4 (pp. 2464-2469). IEEE.

\bibitem{matthieu2021enhanced} Terris M, Repetti A, Pesquet JC, Wiaux Y. Enhanced convergent pnp algorithms for image restoration. In 2021 IEEE International Conference on Image Processing (ICIP) 2021 Sep 19 (pp. 1684-1688). IEEE.

\bibitem{turner2023adaptive} Turner M, Koch T, Serrano F, Winkler M. Adaptive cut selection in mixed-integer linear programming. Open Journal of Mathematical Optimization. 2023;4:1-28.

\bibitem{ulyanov2018deep} Ulyanov D, Vedaldi A, Lempitsky V. Deep image prior. In Proceedings of the IEEE conference on computer vision and pattern recognition 2018 (pp. 9446-9454).

\bibitem{valentin2022instance} Valentin R, Ferrari C, Scheurer J, et al. Instance-wise algorithm configuration with graph neural networks. arXiv preprint arXiv:2202.04910. 2022 Feb 10.

\bibitem{venkatakrishnan2013plug} Venkatakrishnan SV, Bouman CA, Wohlberg B. Plug-and-play priors for model based reconstruction. In 2013 IEEE global conference on signal and information processing 2013 Dec 3 (pp. 945-948). IEEE.

\bibitem{vu} Vũ BC. A splitting algorithm for dual monotone inclusions involving cocoercive operators. Advances in Computational Mathematics. 2013 Apr;38(3):667-81.

\bibitem{wadayama2019deep} Wadayama T, Takabe S. Deep learning-aided trainable projected gradient decoding for LDPC codes. In 2019 IEEE International Symposium on Information Theory (ISIT) 2019 Jul 7 (pp. 2444-2448). IEEE.

\bibitem{wang2016d3} Wang Z, Liu D, Chang S, et al. D3: Deep dual-domain based fast restoration of JPEG-compressed images. In Proceedings of the IEEE Conference on Computer Vision and Pattern Recognition 2016 (pp. 2764-2772).

\bibitem{wang2023learning} Wang Z, Li X, Wang J, et al. Learning Cut Selection for Mixed-Integer Linear Programming via Hierarchical Sequence Model. In The Eleventh International Conference on Learning Representations 2022 Sep 29.

\bibitem{wei2020tuning} Wei K, Aviles-Rivero A, Liang J, et al. Tuning-free plug-and-play proximal algorithm for inverse imaging problems. In International Conference on Machine Learning 2020 Nov 21 (pp. 10158-10169). PMLR.

\bibitem{weng2018evaluating} Weng TW, Zhang H, Chen PY, et al. Evaluating the Robustness of Neural Networks: An Extreme Value Theory Approach. In International Conference on Learning Representations 2018 Feb 15.

\bibitem{Wilder2019Melding} Wilder B, Dilkina B, Tambe M. Melding the data-decisions pipeline: Decision-focused learning for combinatorial optimization. In Proceedings of the AAAI Conference on Artificial Intelligence 2019 Jul 17 (Vol. 33, No. 01, pp. 1658-1665).

\bibitem{wolpert1997no} Wolpert DH, Macready WG. No free lunch theorems for optimization. IEEE transactions on evolutionary computation. 1997 Apr;1(1):67-82.

\bibitem{wolsey2020integer} Wolsey LA. Integer programming. John Wiley \& Sons; 2020 Sep 10.

\bibitem{samy2022jfb} Fung SW, Heaton H, Li Q, et al. Jfb: Jacobian-free backpropagation for implicit networks. In Proceedings of the AAAI Conference on Artificial Intelligence 2022 Jun 28 (Vol. 36, No. 6, pp. 6648-6656).

\bibitem{wu2020sparse} Wu K, Guo Y, Li Z, Zhang C. Sparse coding with gated learned ISTA. In International conference on learning representations 2019 Sep 23.

\bibitem{wu2022gnn} Wu L, Cui P, Pei J, et al. Graph neural networks: foundation, frontiers and applications. In Proceedings of the 28th ACM SIGKDD Conference on Knowledge Discovery and Data Mining 2022 Aug 14 (pp. 4840-4841).

\bibitem{wu2021learning} Wu Y, Song W, Cao Z, Zhang J. Learning large neighborhood search policy for integer programming. Advances in Neural Information Processing Systems. 2021 Dec 6;34:30075-87.

\bibitem{wollmer2013lstm} Wöllmer M, Kaiser M, Eyben F, et al. LSTM-modeling of continuous emotions in an audiovisual affect recognition framework. Image and Vision Computing. 2013 Feb 1;31(2):153-63.

\bibitem{xie2019diff} Xie X, Wu J, Liu G, et al. Differentiable linearized ADMM. In International Conference on Machine Learning 2019 May 24 (pp. 6902-6911). PMLR.

\bibitem{xu2011hydramip} Xu L, Hutter F, Hoos HH, Leyton-Brown K. Hydra-MIP: Automated algorithm configuration and selection for mixed integer programming. In RCRA workshop on experimental evaluation of algorithms for solving problems with combinatorial explosion at the international joint conference on artificial intelligence (IJCAI) 2011 Jul (pp. 16-30).

\bibitem{yang2020learning} Yang C, Gu Y, Chen B, Ma H, So HC. Learning proximal operator methods for nonconvex sparse recovery with theoretical guarantee. IEEE Transactions on Signal Processing. 2020 Mar 5;68:5244-59.

\bibitem{yang2020hyperparameter} Yang L, Shami A. On hyperparameter optimization of machine learning algorithms: Theory and practice. Neurocomputing. 2020 Nov 20;415:295-316.

\bibitem{Kao2009directed} Kao YH, Roy B, Yan X. Directed regression. Advances in Neural Information Processing Systems. 2009;22.

\bibitem{yilmaz2021study} Yilmaz K, Yorke-Smith N. A study of learning search approximation in mixed integer branch and bound: Node selection in scip. AI. 2021 Apr 12;2(2):150-78.

\bibitem{yuan2020plug} Yuan X, Liu Y, Suo J, Dai Q. Plug-and-play algorithms for large-scale snapshot compressive imaging. In Proceedings of the IEEE/CVF Conference on Computer Vision and Pattern Recognition 2020 (pp. 1447-1457).

\bibitem{zarka2020deep} Zarka J, Thiry L, Angles T, Mallat S. Deep Network Classification by Scattering and Homotopy Dictionary Learning. In International Conference on Learning Representations 2020.

\bibitem{zarpellon2021param} Zarpellon G, Jo J, Lodi A, Bengio Y. Parameterizing branch-and-bound search trees to learn branching policies. In Proceedings of the AAAI conference on artificial intelligence 2021 May 18 (Vol. 35, No. 5, pp. 3931-3939).

\bibitem{zhang2018istanet} Zhang J, Ghanem B. ISTA-Net: Interpretable optimization-inspired deep network for image compressive sensing. In Proceedings of the IEEE conference on computer vision and pattern recognition 2018 (pp. 1828-1837).

\bibitem{zhang2021plug} Zhang K, Li Y, Zuo W, et al. Plug-and-play image restoration with deep denoiser prior. IEEE Transactions on Pattern Analysis and Machine Intelligence. 2021 Jun 14;44(10):6360-76.

\bibitem{zhang2017beyond} Zhang K, Zuo W, Chen Y, et al. Beyond a gaussian denoiser: Residual learning of deep cnn for image denoising. IEEE transactions on image processing. 2017 Feb 1;26(7):3142-55.

\bibitem{zhang2017leanring} Zhang K, Zuo W, Gu S, Zhang L. Learning deep CNN denoiser prior for image restoration. In Proceedings of the IEEE conference on computer vision and pattern recognition 2017 (pp. 3929-3938).

\bibitem{zhang2018ffdnet} Zhang K, Zuo W, Zhang L. FFDNet: Toward a fast and flexible solution for CNN-based image denoising. IEEE Transactions on Image Processing. 2018 May 25;27(9):4608-22.

\bibitem{zhang2018dynamically} Zhang X, Lu Y, Liu J, Dong B. Dynamically Unfolding Recurrent Restorer: A Moving Endpoint Control Method for Image Restoration. In International Conference on Learning Representations 2018 Sep 27.

\bibitem{zhang2023mindopt} Zhang M, Yin W, Wang M, et al. Mindopt tuner: Boost the performance of numerical software by automatic parameter tuning. arXiv preprint arXiv:2307.08085. 2023 Jul 16.

\bibitem{zhang2022deep} Zhang T, Banitalebi-Dehkordi A, Zhang Y. Deep reinforcement learning for exact combinatorial optimization: Learning to branch. In 26th international conference on pattern recognition (ICPR) 2022 Aug 21 (pp. 3105-3111). IEEE.

\bibitem{zhao2011online} Zhao B, Fei-Fei L, Xing EP. Online detection of unusual events in videos via dynamic sparse coding. In CVPR 2011 2011 Jun 20 (pp. 3313-3320). IEEE.

\bibitem{zou2022proximal} Zou Y, Zhou Y, Chen X, Eldar YC. Proximal Gradient-Based Unfolding for Massive Random Access in IoT Networks. arXiv preprint arXiv:2212.01839. 2022 Dec 4.

\end{thebibliography}
